\documentclass[pdftex, 12pt,a4paper,leqno, times]{report}
\usepackage[utf8]{inputenc}
\usepackage[T1]{fontenc}
\usepackage[english]{babel}
\usepackage{amsmath}
\usepackage{amssymb, textcomp}
\usepackage{enumerate}
\usepackage{times}
\usepackage[varg]{txfonts}
\usepackage{mathrsfs}
\usepackage{makeidx}
\usepackage[pdftex]{graphicx}
\usepackage{multicol}
\usepackage{color}
\newcommand{\newproof}[3]{
   \newenvironment{#1}[1][]%
   {%
      \begin{trivlist}%
         \item[\hspace{\labelsep}\textnormal{\textbf{#2%
            \def\op@@@arg{##1}%
            \ifx\op@@@arg\empty
            \else~ ##1\fi
         }}]%
   }%
   {%
      #3
      \end{trivlist}%
   }%
}
\newproof{biz}{Proof.}{\qed}
\def\qed{{\ifhmode\unskip\nobreak\hfil\penalty50 \hskip1em \else\nobreak\fi
   \mbox{}\nobreak\hfil\qedtext%
   \parfillskip=0pt \finalhyphendemerits=0 \par}}

\def\qedtext{\ensuremath{\square}}
\RequirePackage{amsfonts}
\DeclareMathSymbol{\square}       {\mathord}{AMSa}{"03}

\newtheorem{theorem}{Theorem}[section]
\newtheorem{lem}{Lemma}[section]
\newtheorem{cor}{Corollary}[section]
\newtheorem{defin}{Definition}[section]
\newtheorem{prop}{Proposition}[section]
\newtheorem{rem}{Remark}[section]
\newtheorem{Opp}{Open Problem}
\newtheorem{ex}{Example}[section]

\long\def\symbolfootnote[#1]#2{\begingroup
\def\thefootnote{\fnsymbol{footnote}}\footnote[#1]{#2}\endgroup}

\def\polhk#1{\setbox0=\hbox{#1}{\ooalign{\hidewidth
    \lower1.5ex\hbox{`}\hidewidth\crcr\unhbox0}}}

\title{Characterizations of derivations}

\author{Eszter Gselmann}

\makeindex

\begin{document}
\nocite{*}

\maketitle

\newpage
\thispagestyle{empty}
\tableofcontents

\newpage

\chapter*{Introduction}
The main purpose of this work is to characterize derivations through functional equations. 
Therefore, (besides the notion of derivations) it is natural to ask what a functional equation is. 
But there is no easy and satisfactory answer to this question. 
While such concepts as element, relation, mapping, operation, etc., 
are well defined in set theory, while such principal concept as set, is an undefined term.
As in set theory, we hope the reader will get a general insight of what this theory is
about.

Functional equations occur almost everywhere. 
Their influence and applications can be felt in every field, 
and all fields benefit from their contact, use and technique. 
The growth and development used to be influenced by their impact on other areas -- 
not only in mathematics but also in other disciplines. 
Applications can be found in a wide variety of fields e.g., analysis, behavioural and social science, 
biology, combinatorics, economics, engineering, geometry, inequalities, information theory, 
physics, psychology, statistics etc. 

Even though lots of mathematicians worked in this area, since the appearance of the famous monograph \cite{Acz66} of 
J.~Aczél\index{Aczél, J.} no systematic exposition existed. In this dissertation we will follow this monograph as well as that of M.~Kuczma\index{Kuczma, M.}, 
see \cite{Kuc09}. 

This work consists of five chapters. 
In the first one, we summarize the most important notions and results from 
the theory of functional equations that will be used afterwards. 
In Chapter 2 we collect all the definitions and results regarding derivations that are essential 
while studying this area. 
Let $Q$ be a ring and let $P$ be a subring of $Q$.
A function $f\colon P\rightarrow Q$ is called a \emph{derivation}\index{derivation} if it is additive,
i.e. 
\[
f(x+y)=f(x)+f(y)
\quad
\left(x, y\in P\right)
\]
and also satisfies the so-called \emph{Leibniz rule}\index{Leibniz rule}, i.e.  
\[
f(xy)=f(x)y+xf(y)
\quad
\left(x, y\in P\right). 
\]
As well as homomorphisms, derivations give a lot of information about the rings between which they act. 
Therefore, the characterization of derivations can also be applicable not only from the theory of 
functional equations but also from the point of view of some algebraic investigations. 
As it can be seen, the notion of derivations is already formulated via functional equations 
(the Cauchy equation and the Leibniz rule). 
In Chapter 3 we intend to show that derivations can be characterized by one single functional equation. 

The results presented here are based on Gselmann \cite{Gse12}. 

More exactly, we would like to examine whether the equations occurring in the definition of
derivations are independent in the following sense.

Let $Q$ be a commutative ring and let $P$ be a subring of $Q$.
Let $\lambda, \mu\in Q\setminus\left\{0\right\}$  be arbitrary,
$f\colon P\rightarrow Q$ be a function and consider the
equation
\[
\lambda\left[f(x+y)-f(x)-f(y)\right]+
\mu\left[f(xy)-xf(y)-yf(x)\right]=0
\quad
\left(x, y\in P\right). 
\]
Clearly, if the function $f$ is a derivation, then this equation holds.
We will however investigate the opposite direction, and it will be proved that
under some assumptions on the rings $P$ and $Q$, derivations can be
characterized via the above equation.
This result will be proved as a consequence of the main theorem that will be devoted
to the study of equation 
\[
f(x+y)-f(x)-f(y)=g(xy)-xg(y)-yg(x)
\quad
\left(x, y\in P\right), 
\]
where $f, g\colon P\rightarrow Q$ are unknown functions. 

Similar problems were already studied by J.~Dhombres in \cite{Dho88}\index{Dhombres, J.}. 
However, the interested reader should also consult the survey paper \cite{GerSab17} of R.~Ger and M.~Sablik\index{Sablik, M.} 
and also the papers Ger \cite{Ger98, Ger00}\index{Ger, R.} and Ger--Reich \cite{GerRei10}\index{Reich, L.}. 

Let $X$ and $Y$ be nonvoid sets and $E_{1}(f)=0$ and $E_{2}(f)=0$ be two functional equations for the function $f\colon X\to Y$. 
We say that equations $E_{1}$ and $E_{2}$ are \emph{alien}\index{alien!--- functional equation} if any solution $f\colon X\to Y$ of the functional equation 
\[
 E_{1}(f)+E_{2}(f)=0
\]
also solves the system 
\[
 \begin{array}{rcl}
  E_{1}(f)&=&0\\
  E_{2}(f)&=&0. 
 \end{array}
\]

Furthermore, equations $E_{1}$ and $E_{2}$ are \emph{strongly alien}\index{strongly alien!--- functional equation} if any pair $f, g\colon X\to Y$ of functions that solves 
\[
 E_{1}(f)+E_{1}(g)=0
\]
also yields a solution for 
\[
  \begin{array}{rcl}
  E_{1}(f)&=&0\\
  E_{2}(g)&=&0. 
 \end{array}
\]
In this setting, our main result says that 
that the (additive) Cauchy equation, i.e.  
\[
 f(x+y)=f(x)+f(y) 
 \qquad 
 \left(x, y\in \mathbb{F}\right)
\]
and the Leibniz rule, that is, 
\[
 f(xy)=xf(y)+f(x)y 
 \qquad 
 \left(x, y\in \mathbb{F}\right). 
\]
are \emph{alien}, but not \emph{strongly alien}.

Chapter 4 is devoted to
 the additive solvability of the system of 
functional equations
\[
 d_{k}(xy)=\sum_{i=0}^{k}\Gamma(i,k-i) d_{i}(x)d_{k-i}(y)
  \qquad (x,y\in \mathbb{R},\,k\in\{0,\ldots,n\}), 
\]
where $\Delta_n:=\big\{(i,j)\in\mathbb{Z}\times\mathbb{Z}\, \vert\,  0\leq i,j\mbox{ and }i+j\leq n\big\}$
and $\Gamma\colon\Delta_n\to\mathbb{R}$ is a symmetric function such that $\Gamma(i,j)=1$ 
whenever $i\cdot j=0$. 

Moreover, the linear dependence and independence of the additive 
solutions $d_{0},d_{1},\dots,d_{n} \allowbreak \colon\mathbb{R}\to\mathbb{R}$ of the above system of equations is characterized. 
As a consequence of the main result, for any nonzero real derivation $d\colon\mathbb{R}\to\mathbb{R}$, the 
iterates $d^0,d^1,\dots,d^n$ of $d$ are shown to be linearly independent, 
and the graph of the mapping $x\mapsto (x,d^1(x),\dots,d^n(x))$ to be dense in $\mathbb{R}^{n+1}$.
The results of this chapter were achieved jointly with Zs.~Páles\index{Páles, Zs.} and were published in 
\cite{GsePal16}. 

Finally, the closing chapter deals with the following problem. 
Assume that $\xi\colon \mathbb{R}\to \mathbb{R}$ 
is a given differentiable function and for the additive function $f\colon \mathbb{R}\to \mathbb{R}$, the mapping 
\[
 \varphi(x)=f\left(\xi(x)\right)-\xi'(x)f(x)
\]
fulfills some regularity condition (e.g. local boundedness, continuity, measurability etc.) on its domain. Is it true that in such a case 
$f$ admits a representation
\[
 f(x)=\chi(x)+f(1)\cdot x 
\quad 
\left(x\in \mathbb{R}\right), 
\]
where $\chi\colon \mathbb{R}\to \mathbb{R}$ is a real derivation?

In case $\varphi$ is identically zero and $\xi(x)=x^{k}$, several results are known for instance due to 
Jurkat \cite{Jur65}\index{Jurkat, W.}, Kurepa\index{Kurepa, S.} \cite{Kur64}, and
Kannappan--Kurepa \cite{KanKur70}\index{Kannappan, Pl.}.
Our investigation in this area began in a joint work with Z.~Boros\index{Boros, Z.}, see \cite{BorGse12}. 
After some preliminary results, in this chapter we will show that 
in case $n\in\mathbb{Z}\setminus\left\{0\right\}$ and 
$\left(\begin{array}{cc}
a&b\\
c&d
\end{array}
\right)\in\mathbf{GL}_{2}(\mathbb{Q})$
and the function $\xi$ is 
\[
 \xi(x)=\dfrac{ax^{n}+b}{cx^{n}+d}
\quad 
\left(x\in \mathbb{R}, cx^{n}+d\neq 0\right), 
\]
then the answer is \emph{affirmative}. These results can be found in \cite{Gse13, Gse14, Gse14b}. 
Furthermore, we will also show that the above class of functions is expandable. 
More precisely, we will show (among others) the following. 
Assume that for the additive function 
$f\colon \mathbb{R}\to \mathbb{R}$ the mapping $\varphi$ defined by 
\[
 \varphi(x)=f\left(\xi(x)\right)-\xi'(x)f(x)
\]
is regular. Then the function $f$ can be represented as 
\[
 f(x)=\chi(x)+f(1)\cdot x 
\quad 
\left(x\in \mathbb{R}\right), 
\]
where $\chi\colon \mathbb{R}\to \mathbb{R}$ is a derivation 
in any of the following cases
\begin{multicols}{4}
\begin{enumerate}[(a)]
 \item \[\xi(x)=a^{x}\]
\item \[\xi(x)=\cos(x)\]
\item \[\xi(x)=\sin(x)\]
\item \[\xi(x)=\cosh(x)\]
\item \[\xi(x)=\sinh(x) \]
 \item \[\xi(x)=\ln(x)\]
\item \[\xi(x)=\mathrm{arccos}(x)\]
\item \[\xi(x)=\mathrm{arcsin}(x)\]
\item \[\xi(x)=\mathrm{arcosh}(x)\]
\item \[\xi(x)=\mathrm{arsinh}(x). \]
\end{enumerate}
\end{multicols}
With the aid of Hyers' theorem and this result we were also able to prove 
stability type results concerning derivations.

I would like to express my gratitude to the many people who saw me through this book; 
to all those who provided support, talked things over, read, wrote, offered comments, 
allowed me to quote their remarks and assisted in the editing, proofreading. 
Most of all, I am indebted to \emph{Professors Gábor Horváth, Gyula Maksa and Zsolt Páles}. 

I would like to cordially thank the anonymous referee for reading my
manuscript carefully and for all of his/her helpful and constructive
comments. Clearly, they helped to improve the quality of this work. 

I would also express my gratitude to \emph{Professors Roman Ger and \.{Z}ywilla Fechner}, who made me 
available monograph \cite{Prz88}. 

The research of the author has been supported by the Hungarian Scientific Research Fund
(OTKA) Grant K 111651 and by the ÚNKP-4 New National Excellence Program of the Ministry of Human Capacities.
The research is also supported by the
EFOP-3.6.1-16-2016-00022 project. The project is
co-financed by the European Union and the European
Social Fund.

Above all, I want to thank my husband, \emph{Alfréd}, my daughter, \emph{Hilda} and the rest of my family, 
who supported and encouraged me in spite of all the time it took me away from them.

\chapter{Preliminaries from the theory of functional equations}\index{functional equation}

As J.~Aczél\index{Aczél, J.} wrote in his famous monograph \cite{Acz66}, 
the theory of functional equations is one of the oldest areas in mathematics. 
Already, J.~D'Alembert\index{D'Alembert, J.}, L.~Euler\index{Euler, L.} , C.F.~Gau\ss\index{Gau\ss, C.F.}, A. L.~Cauchy\index{Cauchy, A.L.}, N.H.~Abel\index{Abel, N.H.},
K.~Weierstra\ss\index{Weierstra\ss, K.}, J.G.~Darboux\index{Darboux, J.G.} and D.~Hilbert\index{Hilbert, D.} considered and also solved functional equations. 

In this section we will summarize the most important notions and results that will be used subsequently. 
For the details the reader should consult the two basic monographs Aczél \cite{Acz66} and Kuczma \cite{Kuc09}\index{Kuczma, M.}. 

\section{Additive functions}

Henceforth $\mathbb{N}, \mathbb{Z}, \mathbb{Q}$, and
$\mathbb{R}$ denote the set of the
natural (positive integer), the integer, the rational,
and the real numbers, respectively. Furthermore 
$\mathbb{R}^{n}$ denotes the $n$-dimensional Euclidean space. 

A function $f\colon \mathbb{R}^{n}\to \mathbb{R}$ is called \emph{additive}\index{additive function} if it satisfies 
\begin{equation}
 f(x+y)=f(x)+f(y)
\end{equation}
for all $x, y\in \mathbb{R}^{n}$. 
In case $n=1$, this equation was investigated by A.~M.~Legendre 
and C.F.~Gau\ss for the first time. It was however A.~L.~Cauchy who first
found its continuous solutions. Therefore this equation has been named after him 
\emph{Cauchy's equation}\index{Cauchy's equation}. 

By induction on $k$ easily follows that if $f\colon \mathbb{R}^{n}\to \mathbb{R}$ is additive, then 
\[
 f\left(\sum_{i=1}^{k}x_{i}\right)=\sum_{i=1}^{k}f(x_{i})
\]
is valid for all $k\in \mathbb{N}$ and for any $x_{1}, \ldots, x_{k}\in \mathbb{R}^{n}$. 
Furthermore, any additive function $f\colon \mathbb{R}^{n}\to \mathbb{R}$ is 
$\mathbb{Q}$-homogeneous as well, that is 
\[
 f(\lambda x)=\lambda f(x)
\]
is fulfilled for any $x\in \mathbb{R}^{n}$ and $\lambda\in \mathbb{Q}$. 

From this we get that if $n=1$ and $f$ is continuous then 
\[
 f(x)=f(1)\cdot x 
\qquad \left(x\in \mathbb{R}\right). 
\]

For many years the existence of discontinuous additive functions was an open problem. 
Finally, in 1905 Hamel\index{Hamel, G.} proved the following, see \cite{Ham05}. 

\begin{theorem}
 Let $\mathscr{H}$ be a Hamel base\index{Hamel base} of the space 
$(\mathbb{R}^{n}; \mathbb{Q}; +; \cdot)$. 
Then for every function $g\colon \mathscr{H}\to \mathbb{R}$ there exists a 
unique additive function 
$f\colon \mathbb{R}^{n}\to \mathbb{R}$ such that $f\vert_{\mathscr{H}}=g$. 
\end{theorem}

In view of the above theorem, assuming the Axiom of Choice, all additive functions 
$f\colon \mathbb{R}^{n}\to \mathbb{R}$ can be received. Indeed, every additive $f$ can be obtained 
as the unique additive extension of a certain function $g\colon \mathscr{H}\to \mathbb{R}$. 

With the aid of this theorem, the following corollary can be derived immediately. 

\begin{cor}
  Let $\mathscr{H}$ be a Hamel base of the space 
$(\mathbb{R}^{n}; \mathbb{Q}; +; \cdot)$. 
Let $g\colon \mathscr{H}\to \mathbb{R}$ be a non-identically zero function for which 
\[
 g(\mathscr{H})=\left\{g(h)\, \vert \, h\in \mathscr{H}\right\} \subset \mathbb{Q}. 
\]
Then the additive extension $f\colon \mathbb{R}^{n}\to \mathbb{R}$ of $g$ 
is a discontinuous additive function. 
\end{cor}

\section{The remaining Cauchy equations}

The following functional equations are also referred to as 
Cauchy's equations. 

A function $f\colon \mathbb{R}^{n}\to \mathbb{R}$ is called 
an \emph{exponential} function\index{exponential function} if 
\begin{equation}
 f(x+y)=f(x)f(y)
\end{equation}
is fulfilled for all $x, y\in \mathbb{R}^{n}$. 

A function $f\colon \mathbb{R}\to \mathbb{R}$ is a 
\emph{logarithmic} function\index{logarithmic function} if 
\begin{equation}
 f(xy)=f(x)+f(y)
\end{equation}
is satisfied for all $x, y\in \mathbb{R}$. 

Finally, a function $f\colon \mathbb{R}\to \mathbb{R}$ is termed to be  
\emph{multiplicative}\index{multiplicative function} if 
\begin{equation}
 f(xy)=f(x)\cdot f(y) 
\qquad 
\left(x, y\in \mathbb{R}\right). 
\end{equation}

\begin{theorem}
 Let $f\colon \mathbb{R}^{n}\to \mathbb{R}$ be an exponential function. 
Then either $f$ is identically zero or there is an additive function 
$a\colon \mathbb{R}^{n}\to \mathbb{R}$ so that 
\[
 f(x)=\exp\left(a(x)\right) 
\qquad 
\left(x\in \mathbb{R}^{n}\right). 
\]
\end{theorem}

The natural domains of definition for logarithmic and multiplicative functions, resp. are sets of the form 
\[
 \mathscr{L}
=
\left\{(x, y)\in A\times A\, \vert\,  x\cdot y\in A\right\}, 
\]
while for exponential functions, sets of the form 
\[
 \mathscr{A}
=
\left\{(x, y)\in A\times A\, \vert\,  x+y\in A\right\}, 
\]
with a certain nonempty set $A\subset \mathbb{R}$. 

\begin{lem}
 Let $A\subset \mathbb{R}$ be a nonempty set and 
$f\colon A\to \mathbb{R}$ be a function. 
Assume 
\[
 f(xy)=f(x)+f(y)
\]
is valid for all $(x, y)\in \mathscr{L}$. If 
$0\in A$, then $f$ is identically zero. 
\end{lem}

\begin{theorem}
 Let $A=]0, +\infty[$ or $A=\mathbb{R}\setminus\left\{0\right\}$. 
If a function $f\colon A \to \mathbb{R}$ fulfills the logarithmic Cauchy equation 
for all pairs $(x, y)\in \mathscr{L}$, then 
there exists an additive function 
$a\colon \mathbb{R}\to \mathbb{R}$ such that 
\[
 f(x)=a\left(\ln \left(|x|\right)\right) 
\qquad 
\left(x\in A\right). 
\]
\end{theorem}

\begin{theorem}
 Let $f\colon \mathbb{R}\to \mathbb{R}$ be a multiplicative function. 
Then there exists an additive function 
$a\colon \mathbb{R}\to \mathbb{R}$ such that $f$ is one of the following forms. 
\begin{multicols}{2}
\[
 f(x)=0 
\qquad 
\left(x\in \mathbb{R}\right), 
\]

\[
 f(x)=1
\qquad 
\left(x\in \mathbb{R}\right),
\]
\[
 f(x)=
\begin{cases}
 \exp\left(g\left(\ln(|x|)\right)\right)& \text{ if } x\neq 0\\
0& \text{ if } x=0
\end{cases}
\]

\[
 f(x)=
\begin{cases}
 \exp\left(g\left(\ln(|x|)\right)\right)& \text{ if } x> 0\\
0& \text{ if } x=0\\
-\exp\left(g\left(\ln(|x|)\right)\right)& \text{ if } x< 0
\end{cases}
\]
\end{multicols}
\end{theorem}

\section{Jensen equation and Hosszú equation}\index{Jensen equation}\index{Hosszú equation}

The equation resulting on replacing in the so-called Jensen inequality, i.e.  
\[
 f\left(\frac{x+y}{2}\right)\leq \frac{f(x)+f(y)}{2}
\]
the sign of inequality by that of equality, that is, 
\begin{equation}\label{Eq1.1.5}
 f\left(\frac{x+y}{2}\right)=\frac{f(x)+f(y)}{2}
\end{equation}
is known as the \emph{Jensen equation}. 

Usually \eqref{Eq1.1.5} is considered for functions $f\colon D\to \mathbb{R}$, where $D\subset \mathbb{R}^{n}$ is a convex set. 
If $D$ were also open, then any function $f$ fulfilling \eqref{Eq1.1.5} would also be Jensen convex, and thus all the results concerning 
convex functions (e.g. Bernstein--Doetsch theorem) would apply. Most of these result however became invalid when the set $D$ is not necessarily open. 

\begin{theorem}
 Let $D\subset \mathbb{R}^{n}$ be a convex set, such that $\mathrm{int}(D)\neq \emptyset$ further let
 $f\colon D\to \mathbb{R}$ be a solution of \eqref{Eq1.1.5}. 
 Then there exist an additive function $a\colon \mathbb{R}^{n}\to \mathbb{R}$ and a constant $\alpha\in \mathbb{R}$ such that 
 \[
  f(x)=a(x)+\alpha 
  \qquad 
  \left(x\in D\right). 
 \]
\end{theorem}

The functional equation
\begin{equation}
 f(x+y-xy)+f(xy)=f(x)+f(y)
\end{equation}
was mentioned for the first time by M.~Hosszú\index{Hosszú, M.} 
in 1967 at the International Symposium on Functional Equations held in Zakopane (Poland) and it is named after him \emph{Hosszú equation}.  
After its first appearance it was extensively studied by among others Blanu\v{s}a \cite{Bla70}\index{Blanu\v{s}a, D.}, 
Daróczy \cite{Dar71, Dar69}\index{Daróczy, Z.}, Davison \cite{Dav74, Dav74b}\index{Davison, T.~M.~K}, 
Davison--Redlin \cite{DavRed80}\index{Redlin, L.}, Fenyő \cite{Fen69}, G{\l}owacki--Kuczma \cite{GloKuc79}\index{G{\l}owacki, E.}\index{Kuczma, M.}, 
Lajkó \cite{Laj74}\index{Lajkó, K.}, \'Swiatak \cite{Swi68b, Swi68, Swi71}\index{\'Swiatak, H.}. 

In a short time turned out that for functions $f\colon \mathbb{R}\to \mathbb{R}$ Hosszú equation and Jensen equation 
are \emph{equivalent}. At the same time, even in rather simple cases it can occur that these two functional equations have a 
different set of general solutions. 
In particular, on $\mathbb{Z}$ the function $f\colon \mathbb{Z}\to \mathbb{Z}$ defined by 
\[
 f(x)=
 \left\{
 \begin{array}{rl}
  1, & \text{ if $n$ is even}\\
  0, & \text{ if $n$ is odd}
 \end{array}
\right.
\qquad 
\left(x\in \mathbb{Z}\right). 
\]
solves the Hosszú equation, however, does not solve the Jensen equation. 

For functions acting between fields and commutative groups, in Davison \cite{Dav74} the following was proved. 

\begin{theorem}
Let $(G, +)$ be an abelian group and $(\mathbb{K}, +, \cdot)$ be a field having at least five elements. 
Suppose that for the function $f\colon \mathbb{K}\to G$ Hosszú equation is satisfied, that is, we have 
\[f(x+y-xy)+f(xy)=f(x)+f(y)
\qquad 
\left(x, y\in \mathbb{K}\right). 
\] 
Then 
\[
f(x)=g(x)-g(0)
\qquad
\left(x \in \mathbb{K}\right), 
\]
where $g\colon \mathbb{K}\to G$ is a homomorphism. 
\end{theorem}

\section{Polynomial functions}\index{polynomial function}

The study of polynomial functions defined on groups is based on the notion of multiadditive functions. 
Therefore, firstly we collect some basic notions and results concerning such functions. 
Here we follow the notations and the terminology of Székelyhidi \cite{Sze91, Sze06}.\index{Székelyhidi, L.}

\subsection*{Multiadditive functions}\index{multiadditive function}

\begin{defin}
 Let $G, S$ be commutative semigroups, $n\in \mathbb{N}$ and let $A\colon G^{n}\to S$ be a function. 
 We say that $A$ is \emph{$n$-additive} if it is a homomorphism of $G$ into $S$ in each variable. 
 If $n=1$ or $n=2$ the function $A$ is simply termed to be \emph{additive}\index{additive function} 
 or \emph{biadditive}\index{biadditive function}, respectively. Further, letting 
 $G^{0}=G$, constant functions from $G$ to $S$ will be called \emph{$0$-additive}. 
 We call the function $A\colon G^{n}\to S$ \emph{multiadditive}, if there is a natural number $n$ such that 
 $A$ is $n$-additive. 
 
 The \emph{diagonalization}\index{multiadditive function!--- diagonalization} or \emph{trace} of an $n$-additive function $A\colon G^{n}\to S$ is defined as 
 \[
  A^{\ast}(x)=A\left([x]_{n}\right)=A\left(x, \ldots, x\right) 
  \qquad 
  \left(x\in G\right). 
 \]
\end{defin}

The following proposition contains the most basic properties of multiadditive functions. 

\begin{prop}
 Let $G, S$ be commutative semigroups, $n\in \mathbb{N}$. Then the set of all $n$-additive functions from 
 $G^{n}$ to $S$ forms a 
 \begin{enumerate}[(i)]
  \item commutative semigroup, if $S$ is a commutative semigroup;
  \item module over the ring $R$, if $S$ is a module over the ring $R$;
  \item linear space over the field $\mathbb{F}$, if $S$ is a linear space over $\mathbb{F}$. 
 \end{enumerate}
\end{prop}

From the definition of multiadditive functions it follows that,
if $n\in \mathbb{N}$ then each $n$-additive function $A\colon G^{n}\to S$ satisfies 
\[
 A(x_{1}, \ldots, x_{i-1}, kx_{i}, x_{i+1})=
 kA(x_{1}, \ldots, x_{i-1}, x_{i}, x_{i+1}, \ldots, x_{n})
 \qquad 
 \left(x_{1}, \ldots, x_{n}\in G\right)
\]
for all $i=1, \ldots, n$, where $k\in \mathbb{N}$ is arbitrary. 
Further, the same identity holds for any $k\in \mathbb{Z}$ if $G$ and $S$ are groups, and 
for $k\in \mathbb{Q}$, if $G$ and $S$ are linear spaces over the rationals. 
For the diagonalization of $A$ we have 
\[
 A^{\ast}(kx)=k^{n}A^{\ast}(x)
 \qquad
 \left(x\in G\right). 
\]

\begin{lem}[Binomial theorem]\index{binomial theorem}
 Let $G, S$ be commutative semigroups and $n\in \mathbb{N}$. 
 If the function $A\colon G^{n}\to S$ is symmetric and $n$-additive, then 
 \[
  A^{\ast}(x+y)
  =
  \sum_{k=0}^{n}\binom{n}{k}A\left([x]_{k}, [y]_{n-k}\right)
 \]
holds for all $x, y\in G$. 
\end{lem}

As a consequence of this result, the so-called \emph{Polynomial theorem} can also be obtained. 

\begin{lem}[Polynomial theorem]\index{polynomial theorem}
 Let $G, S$ be commutative semigroups and $n\in \mathbb{N}$. 
 If the function $A\colon G^{n}\to S$ is symmetric and $n$-additive, then 
 for any $m\in \mathbb{N}$ and for all $x_{1}, \ldots, x_{m}\in G$
 \[
  A^{\ast}(x_{1}+\cdots+x_{m})
  =
  \sum_{k_{1}+\cdots+k_{m}=n}\dfrac{n!}{k_{1}!\cdots k_{m}!}A\left([x_{1}]_{k_{1}}, \ldots, [x_{m}]_{k_{m}}\right)
 \]
is fulfilled. 
\end{lem}

With the aid of the abovementioned Binomial theorem, one can also prove the so-called 
\emph{Polarization formula}, that briefly expresses that (under some conditions on the domain as well as on the range)
every $n$-additive function is \emph{uniquely} determined by its diagonalization. 

\begin{theorem}[Polarization formula]\index{polarization formula}
 Let $G$ be a commutative semigroup and $S$ be a commutative group and $n\in \mathbb{N}$. 
 If $A\colon G^{n}\to S$ is a symmetric, $n$-additive function, then for all 
 $x, y_{1}, \ldots, y_{m}\in G$ we have 
 \[
  \Delta_{y_{1}, \ldots, y_{m}}A^{\ast}(x)=
  \left\{
  \begin{array}{rcl}
   0 & \text{ if} & m>n \\
   n!A(y_{1}, \ldots, y_{m}) & \text{ if}& m=n. 
  \end{array}
  \right.
 \]
\end{theorem}

\begin{cor}
 Let $G$ be a commutative semigroup and $S$ be a commutative group and $n\in \mathbb{N}$. 
 If $A\colon G^{n}\to S$ is a symmetric, $n$-additive function, then 
 \[
  \Delta^{n}_{y}A^{\ast}(x)=n!A^{\ast}(y)
  \qquad 
  \left(x, y\in G\right). 
 \]
\end{cor}

\begin{lem}
  Let $G$ be a commutative semigroup and $S$ be a commutative group and $n\in \mathbb{N}$. 
  Let us assume that the multiplication by $n!$ is surjective in $G$ or injective in $S$. 
  Then for any symmetric, $n$-additive function $A\colon G^{n}\to S$, $A^{\ast}=0$ implies that 
  $A$ is identically zero, as well. 
\end{lem}

\subsection*{Polynomial functions}

The theory of polynomial functions was firstly investigated by M.~Fréchet\index{Fréchet, M.}, S.~Banach\index{Banach, S.}, 
G.~Van der Lijn\index{Van der Lijn, G.}, S.~Mazur\index{Mazur, S.}, W.~Orlicz\index{Orlicz, W.}, 
who were primarily interested in polynomial operations on semigroups and linear spaces. 
The notion we will use is due to M.~Fréchet and S.~Banach. 
In this work we will however restrict ourselves only to the most basic notions and results, the interested reader 
should consult e.g. the two monographs of Székelyhidi \cite{Sze06, Sze91}, Fréchet \cite{Fre29}, Van der Lijn \cite{Lij39}, 
Mazur--Orlicz \cite{MazOrl48, MazOrl53} and also Székelyhidi \cite{Sze79, Sze82}\index{Székelyhidi, L.}. 

\begin{defin}
 Let $G, S$ be commutative semigroups, a function $p\colon G\to S$ is called a \emph{polynomial function}\index{polynomial function}, 
 if it has a representation as the sum of diagonalizations of multiadditive functions from $G$ to $S$. 
 In other words, the function $p\colon G\to S$ is a polynomial if and only if it can be written as 
 \[
  p=\sum_{k=0}^{n}A^{\ast}_{k}, 
 \]
where $n\in \mathbb{N}$, the $A_{k}\colon G^{k}\to S$ is a $k$-additive function for all 
$k=0, 1, \ldots, n$. In this case we also say that $p$ 
is a \emph{polynomial of degree at most $n$}.\index{polynomial function!--- degree at most $n$} 
\end{defin}

\begin{lem}
 Let $G$ be a commutative semigroup and $R$ be a ring. Then the set of all polynomials from $G$ 
 into the additive group of $R$ forms a 
 \begin{enumerate}[(i)]
  \item (commutative) ring, if $R$ is a (commutative) ring;
  \item (commutative) algebra over the field $\mathbb{F}$, if $R$ is a 
  (commutative) algebra over the field $\mathbb{F}$. 
 \end{enumerate}
\end{lem}

\begin{theorem}
 Let $G$ be a commutative semigroup, $S$ be a commutative group and $n\in \mathbb{N}$. 
 Let us assume further the the multiplication by $n!$ is bijective on $G$ or on $S$. 
 Then any polynomial $p\colon G\to S$ of degree at most $n$ has a \emph{unique} representation in the 
 form 
 \[
  p=\sum_{k=0}^{n}A^{\ast}_{k}, 
 \]
where $A_{k}\colon G^{k}\to S$ is a symmetric, $k$-additive function for all $k=0, 1, \ldots, n$. Further, 
$A_{n}$ is not identically zero, whenever $p$ is not of degree $n-1$. 
\end{theorem}

Under the assumptions of the previous theorem, the given representation is called the 
\emph{canonical representation}\index{polynomial function!--- canonical representation} of the polynomial $p$, 
and we call $A^{\ast}_{k}$ it \emph{homogeneous term of degree $k$}\index{polynomial function!--- homogeneous term}, 
further $A^{\ast}_{n}$ is called its \emph{leading term}\index{polynomial function!--- leading term}, 
whenever it is not identically zero.

\section{Extension theorems}\index{extension theorems}

As we will see later, it may happen that a certain functional equation does not hold on the whole possible domain 
but only on its subset. In this case the question is whether the functions appearing in the functional equation 
are extendable to the whole space so that also the functional equation is fulfilled on the whole domain. 

Let $\mathscr{G}\subset \mathbb{R}^{2n}$ be a nonempty set. 
In case $p\in \mathscr{G}$, then there are $x, y\in \mathbb{R}^{n}$ so that 
$p=(x, y)$. 
Further, let us define the following sets. 
\[
 \mathscr{G}_{1}=
\left\{x\in \mathbb{R}^{n}\, \vert \,  \text{there exists $y\in \mathbb{R}^{n}$ such that $(x, y)\in \mathscr{G}$}\right\}. 
\]

\[
 \mathscr{G}_{2}=
\left\{y\in \mathbb{R}^{n}\, \vert \,  \text{there exists $x\in \mathbb{R}^{n}$ such that $(x, y)\in \mathscr{G}$}\right\}. 
\]
and 
\[
 \mathscr{G}_{3}=\left\{x+y\in \mathbb{R}^{n}\, \vert \, (x, y)\in \mathscr{G}\right\}. 
\]

For the sake of brevity, we say that a set $A\subset \mathbb{R}^{n}$ has the $(\bullet)$ property, if 
\[
\tag{$\bullet$} \mathscr{G}_{1}\cup \mathscr{G}_{2}\cup\mathscr{G}_{3}\subset A
\]
is satisfied. 

\begin{theorem}[Daróczy--Losonczi \cite{DarLos67}]\index{Daróczy, Z.}\index{Losonczi, L.}
 Let $\mathscr{G}\subset \mathbb{R}^{2n}$ be an open and connected set and 
let $A\subset \mathbb{R}^{n}$ be a set having the $(\bullet)$ property. 
Assume that $f\colon A\to \mathbb{R}$ fulfills 
\[
 f(x+y)=f(x)+f(y)
\]
for all $(x, y)\in \mathscr{G}$. 
Then there exists a uniquely determined additive function 
$a\colon \mathbb{R}^{n}\to \mathbb{R}$ and unique constants 
$\alpha, \beta\in \mathbb{R}$
such that 
\[
f(x)=
 \begin{cases}
  a(x)+\alpha& \text{if $x\in \mathscr{G}_{1}$}\\
a(x)+\beta& \text{if $x\in \mathscr{G}_{2}$}\\
a(x)+\alpha+\beta& \text{if $x\in \mathscr{G}_{3}$}\\
 \end{cases}
\]
\end{theorem}

Let $r>0$ then 
\[
 B(x_{0}, r)=\left\{x\in \mathbb{R}^{n}\, \vert \,\left\| x-x_{0}\right\| <r\right\}
\]
will denote the open ball centered at $x_{0}$ with radius $r$.

\begin{cor}
 Let $r>0$  and 
$f\colon B(0, 2r)\to \mathbb{R}$ be a function such that 
\[
 f(x+y)=f(x)+f(y)
\]
holds for all $x, y\in B(0, r)$. 
Then there exists a uniquely determined additive function such that 
\[
 f(x)=a(x) 
\qquad 
\left(x\in B(0, r)\right). 
\]
\end{cor}

In the sequel, we will use the following extension theorem concerning the 
so-called Pexider equation, this result is a special case of 
\cite[Theorem 3]{ChuTab08} if we choose the normed space 
$X$ to be $\mathbb{R}^{n}$ and  $\mathscr{G}$ is an  open and connected subset of $X\times X$. 
We have to emphasize that here we do not have to suppose that some set 
$A\subset \mathbb{R}^{n}$ has the property $(\bullet)$. 

\begin{theorem}[Chudziak--Tabor]\index{Chudziak, J.}\index{Tabor, Jó.}
Assume that for the functions 
$k\colon \mathscr{G}_{3}\to \mathbb{R}$, 
$l\colon \allowbreak \mathscr{G}_{1}\to \mathbb{R}$ and 
$n\colon \mathscr{G}_{2}\to \mathbb{R}$

\[
 k(x+y)=l(x)+n(y)
\]
is fulfilled for any $x, y\in \mathscr{G}$. 
If the function $k$ is nonconstant, then these functions can  
\emph{uniquely} be extended to functions 
$\widetilde{k}, \widetilde{l}, \widetilde{n}\colon \mathbb{R}^{n}\to \mathbb{R}$ so that 
\[
 \widetilde{k}(x+y)=\widetilde{l}(x)+\widetilde{n}(y)  
\qquad 
\left(x, y\in \mathbb{R}^{n}\right). 
\]
\end{theorem}

Especially, if a 
function $a\colon \mathscr{G}_{3}\to \mathbb{R}$ is additive on $\mathscr{G}_{3}$, then it can  
always be \emph{uniquely} extended to an additive function $\widetilde{a}\colon\mathbb{R}^{n}\to \mathbb{R}$. 
We remark that this follows also from 
Theorem 4 of P\'{a}les \cite{Pal02}.\index{Páles, Zs.} 

In what follows, we shall present the general solution of the Pexider equation on a restricted domain. This 
theorem follows immediately from \cite[Theorem 1]{ChuTab08}, with exactly the same choice as above. 

\begin{theorem}
 Assume that, for the functions 
$k\colon \mathscr{G}_{3}\to \mathbb{R}$, 
$l\colon \mathscr{G}_{1}\to \mathbb{R}$ and 
$n\colon \mathscr{G}_{2}\to \mathbb{R}$,
\[
 k(x+y)=l(x)+n(y)
\]
is fulfilled for any $x, y\in \mathscr{G}$ and that 
the function $k$ is nonconstant. 
Then and only then there exists a uniquely determined 
additive function $a\colon \mathbb{R}^{n}\to \mathbb{R}$ and real constants 
$b$ and $c$ so that 
\[
 \begin{array}{rcl}
  k(x)&=&a(x)+b+c\\
l(x)&=&a(x)+b\\
n(x)&=&a(x)+c
 \end{array}
\qquad 
\left(x\in \mathbb{R}^{n}\right). 
\]
\end{theorem}

Finally, the following statement concerns the constant solutions of the Pexider equation.

\begin{cor}
 Let $k\in \mathbb{R}$ be fixed and 
$l\colon \mathscr{G}_{1}\to \mathbb{R}$ 
and $n\colon \mathscr{G}_{2}\to \mathbb{R}$
be functions so that 
\[
 k= l(x)+n(y)
\]
is fulfilled for all $x, y\in \mathscr{G}$. 
Then there exists $c\in \mathbb{R}$ such that 
\[
 n(x)= c 
\quad 
\text{and} 
\quad l(x)=k-c 
\qquad
\left(x\in \mathscr{G}\right). 
\]
\end{cor}

\section{Regularity theorems}\index{regularity theorems}

In this subsection a few so-called regularity results will be presented. 
Without completedness,  
here we list only those notions and statements that will be applied later. 
For further results we refer to the monograph \cite{Jar05} of A.~Járai\index{Járai, A.}  and also to the 
two monographs \cite{Sze91, Sze06} of L.~Szé\-kely\-hidi\index{Székelyhidi, L.} as well as  
Szé\-kely\-hidi \cite{Sze85}.

The first classical results in the regularity theory of functional equations concerned mostly additive functions. 
Later these results have been extended to multiadditive and polynomial functions, respectively. 

We begin with the description of continuous, respectively measurable additive functions, see Kuczma \cite{Kuc09}.

\begin{theorem}
 Let $f\colon \mathbb{R}^{n}\to \mathbb{R}$ be a continuous additive function. 
Then there exists a uniquely determined $\lambda=(\lambda_{1}, \ldots, \lambda_{n})\in \mathbb{R}^{n}$ such that 
\[
 f(x)=\langle \lambda, x \rangle 
=
\sum_{k=1}^{n} \lambda_{k}x_{k}
\qquad 
\left(x=(x_{1}, \ldots, x_{n})\in \mathbb{R}^{n}\right)
\]
holds. 
\end{theorem}

\begin{theorem}[Fréchet \cite{Fre13, Fre14}]\index{Fréchet, M.}
 Let $f\colon \mathbb{R}^{n}\to \mathbb{R}$ be an additive function 
which is measurable in the sense of Lebesgue. 
Then $f$ is continuous. 
\end{theorem}

\begin{theorem}
 Let $f\colon \mathbb{R}^{n}\to \mathbb{R}$ be an additive function. 
Assume that there exists a set $A\subset \mathbb{R}^{n}$ with positive Lebesgue measure such that 
the restriction of the function $f$ to this set $A$ is a function that is bounded above or below. 
Then $f$ is a continuous additive function. 
\end{theorem}

The following theorem (see also Kuczma \cite[Theorem 13.4.3]{Kuc09}) is about continuous mul\-ti\-additive functions.

\begin{theorem}
 Let us assume that the p-additive function $f\colon \mathbb{R}^{np}\to \mathbb{R}$ is a continuous p-additive function. 
 Then there exist constants $c_{j_{1}, \ldots, j_{p}}\in \mathbb{R}$, $j_{1}, \ldots, j_{p}=1, \ldots, N$, such that 
 \[
  f(x_{1}, \ldots, x_{p})= \sum_{j_{1}=1}^{N}\cdots \sum_{j_{p}=1}^{N}c_{j_{1}, \ldots, j_{p}}
  x_{ j_{1}}\cdots x_{j_{p}}, 
  \]
where $(i_{i})=(x_{i_{1}}, \ldots, x_{in})$, $i=1, \ldots, p$. 
\end{theorem}

The results presented below show that in general mild regularity conditions on polynomials imply their continuity. 
Here we lean on the two monographs Székelyhidi  \cite{Sze91, Sze06}. 

Due to a famous result of A.~Haar\index{Haar, A.} that can be found in \cite{Haa33}, on every locally compact topological 
group, Haar measure does exist. In the results we will present below, measurability is understood always in this sense.

\begin{theorem}
 Let $G$ be a commutative semigroup and $X$ be a locally convex topological vector space. 
 Then any bounded polynomial $p\colon G\to X$ is constant. 
\end{theorem}

\begin{theorem}
 Let $G$ be a locally compact commutative group which is generated by any neighbourhood of zero and let $X$ 
 be a linear space. 
 If a polynomial $p\colon G\to X$ vanishes on a measurable set of positive measure, then it vanishes everywhere. 
\end{theorem}

The result below is about continuous polynomials, see Theorem 3.2 of \cite{Sze91}, as well. 

\begin{theorem}\label{Tszek1}
 Let $G$ be a locally compact commutative group which is generated by any neighbourhood of zero and let $X$ 
 be a locally convex topological linear space. 
 If a polynomial $p\colon G\to X$ is continuous at a point, then it is continuous everywhere. 
\end{theorem}

The following three statements will play a key role in Chapter 5, they can also be found in 
Székelyhidi \cite{Sze91} as Theorems 3.7, 3.8 and 3.9. 

\begin{theorem}\label{Tszek2}
 Let $G$ be a locally compact commutative group which is generated by any neighbourhood of zero and let $X$ 
 be a locally convex topological linear space. 
 If a polynomial $p\colon G\to X$ is bounded on a nonvoid open set, then it is continuous everywhere. 
\end{theorem}

\begin{theorem}\label{Tszek3}
 Let $G$ be a locally compact commutative group which is generated by any neighbourhood of zero and let $X$ 
 be a locally convex topological linear space. 
 If a polynomial $p\colon G\to X$ is bounded on a measurable set of positive measure, then it is continuous everywhere. 
\end{theorem}

\begin{theorem}\label{Tszek4}
 Let $G$ be a locally compact commutative group which is generated by any neighbourhood of zero and let $X$ 
 be a locally convex and locally bounded topological linear space. 
 If a polynomial $p\colon G\to X$ is measurable on a measurable set of positive measure, then it is continuous everywhere. 
\end{theorem}

\section{Alien functional equations}

The concept of alien functional equations was introduced and developed by J.~Dhombres in the paper \cite{Dho88}\index{Dhombres, J.}. 
However, the interested reader should also consult the survey paper \cite{GerSab17} of R.~Ger and M.~Sablik\index{Sablik, M.} 
and also the papers Ger \cite{Ger98, Ger00}\index{Ger, R.} and Ger--Reich \cite{GerRei10}\index{Reich, L.}. 

Let $X$ and $Y$ be nonvoid sets and $E_{1}(f)=0$ and $E_{2}(f)=0$ be two functional equations for the function $f\colon X\to Y$. 
We say that equations $E_{1}$ and $E_{2}$ are \emph{alien}\index{alien!--- functional equation}, if any solution $f\colon X\to Y$ of the functional equation 
\[
 E_{1}(f)+E_{2}(f)=0
\]
also solves the system 
\[
 \begin{array}{rcl}
  E_{1}(f)&=&0\\
  E_{2}(f)&=&0. 
 \end{array}
\]

Furthermore, equations $E_{1}$ and $E_{2}$ are \emph{strongly alien}\index{strongly alien!--- functional equation}, if any pair $f, g\colon X\to Y$ of functions that solves 
\[
 E_{1}(f)+E_{1}(g)=0
\]
also yields a solution for 
\[
  \begin{array}{rcl}
  E_{1}(f)&=&0\\
  E_{2}(g)&=&0. 
 \end{array}
\]

Concerning the additive and the multiplicative Cauchy equation  in Ger \cite{Ger98} the following theorem was proved. 
\begin{theorem}
 Let $X$ and $Y$ are two rings, and assume that 
for all $x\in X$ there exists $e_{x}\in X$ such that $x e_{x}=x$, suppose further 
that $Y$ has no elements of order 2 and does not admit zero divisors. 
If $f$ is a solution of  
\[
 f(x+y)+f(xy)=f(x)+f(y)+f(x)f(y) 
\qquad 
\left(x, y\in X\right)
\]
such that $f(0)=0$, 
then either $3f$ is even and $3f(2x)=0$ for all $x\in X$, 
or $f$ yields a homomorphism between $X$ and $Y$. 
\end{theorem}

In a recent paper \cite{MakSab16} Gy.~Maksa\index{Maksa, Gy.} and M.~Sablik\index{Sablik, M.} 
investigated whether the exponential Cauchy equation and the Hosszú equation 
are (strongly) alien. In the abovementioned paper they proved the following. 

\begin{theorem}
 Suppose that the functions $f, g \colon \mathbb{R}\to \mathbb{R}$ satisfy functional equation 
 \[
  g(x)g(y)-g(xy)= f(x+y-xy)-f(x)-f(y)+f(xy) 
  \qquad 
  \left(x, y\in \mathbb{R}\right). 
 \]
Assuming that the function $f$ is continuous, the function $g$ is a solution of the exponential Cauchy equation and 
\[
 f(x)=\alpha x+\beta 
 \qquad 
 \left(x\in \mathbb{R}\right)
\]
is fulfilled with a certain constants $\alpha, \beta \in \mathbb{R}$. 
\end{theorem}

In view of the above notions, these results express that (under some regularity assumptions) the exponential Cauchy equation 
and the Hosszú equation are strongly alien.  

To make a more complete picture about the phenomenon of alienity of functional equations, here we also present a result from 
Ger \cite{Ger10}. 
Here the author investigates the additive and the exponential Cauchy equation in this setting. 

\begin{theorem}
Let $(S; +; 0)$ be an Abelian monoid and let $(R; +; \cdot; 1)$ stand for an integral domain. 
If functions $f, g\colon S \to R$ satisfy
\begin{equation}\label{eq:ger1}
f(x+y)+g(x+y)= f(x)+f(x)+g(x)g(y) 
\end{equation}
for all $x, y\in S$, then there exist constants 
$p, q\in R$, $q\neq 0$, additive maps $a, A\colon S\to R$ and a function $r \colon S \to R$ such that
\[
 p r(x+y)= r(x)r(y) 
 \qquad 
 \left(x, y\in S\right)
\]
so that either
\begin{equation}\label{eq:ger2}
 \begin{cases}
  q^{2}f(x)=a(x)+(p-q)r(x)+p(q-p)& \; (x\in S)\\
  qg(x)=r(x)+q-p&\; (x\in S)
 \end{cases}
\end{equation}
or
\begin{equation}\label{eq:ger3}
 \begin{cases}
 2f(x)=a(x)^2+A(x)& \; (x\in S) \\
 g(x)=1-a(x) & \; (x\in S)
 \end{cases}
\end{equation}
Conversely, each pair of functions $f , g \colon S \to R$ satisfying either of the systems
\eqref{eq:ger2}, \eqref{eq:ger3} yields a solution to equation \eqref{eq:ger1}.
\end{theorem}

Since the appearance of the result of J.~Dhombres (see the abovementioned paper \cite{Dho88}), the 
notion of alienity had been extended and generalized in several ways. In this section we restricted however ourselves 
only to those notions that will be used in the next chapters. The interested reader can found various interesting notions and results 
in the abovementioned survey paper \cite{GerSab17}.

\chapter{Derivations}

The main aim of this work is to present some characterization theorems concerning derivations. 
Thus, at first, we list some preliminary results that will be used in the sequel.
All of these statements and definitions can be found in Kuczma \cite{Kuc09}\index{Kuczma, M.} and
in Zariski--Samuel\index{Zariski, O.}\index{Samuel, P.} \cite{ZarSam75} and also in Kharchenko\index{Kharchenko, V. L.} \cite{Kha91}.

Since this work is about derivations, we endeavour to be as self-contained as it is possible concerning derivations. 
Therefore, in this section most of the statements are presented with proofs. For more results see Chapter 14 of Kuczma \cite{Kuc09}.

\section{Derivations}

\begin{defin}\label{D2.1.1}
Let $Q$ be a ring and let $P$ be a subring of $Q$.
A function $f\colon P\rightarrow Q$ is called a \emph{derivation}\index{derivation} if it is additive,
i.e. 
\[
f(x+y)=f(x)+f(y)
\quad
\left(x, y\in P\right)
\]
and also satisfies the so-called \emph{Leibniz rule}\index{Leibniz rule}, i.e.  equation
\[
f(xy)=f(x)y+xf(y)
\quad
\left(x, y\in P\right). 
\]
\end{defin}

From the above definition, the following proposition follows immediately. 

\begin{prop}\label{P2.1.1}
 Let $Q$ be a ring and let $P$ be a subring of $Q$, and let 
 $f, g\colon P\to Q$ be derivations and $\alpha, \beta\in P$. Then 
 \begin{enumerate}[(i)]
  \item the function $\alpha f+\beta g$ is also a derivation;
  \item the function $f\alpha +g\beta$ is also a derivation;
  \item assuming that $Q=P$, the \emph{bracket}\index{bracket} of $f$ and $g$, that is, 
  \[
   \left[f, g\right]= f\circ g -g\circ f
  \]
is also a derivation. 
 \end{enumerate}
\end{prop}

In connection with the composition of derivations we remark the following result of E.~C.~Posner\index{Posner, E.~C.}, see \cite{Pos57}. 

\begin{theorem}
 Let $P$ be a prime ring with $\mathrm{char}(P)\neq 2$, and assume that we are given two derivations 
 $f, g\colon P\to P$. Then the mapping $f\circ g \colon P\to P$ is a derivation, if and only if 
 $f=0$ or $g=0$. 
\end{theorem}

The set of derivations of the ring $R$ will be denoted by $\mathrm{Der}(R)$, that is, 
\[
 \mathrm{Der}(R)=\left\{f\colon R\to \mathbb{R} \, \vert \, f \text{ is a derivation}\right\}. 
\]
In view of the above proposition, this set is closed relative to the bracket operator. Therefore, 
$\mathrm{Der}(R)$ is a Lie ring, which is simultaneously a right module over the center of the ring $R$, 
provided that the multiplication by the central elements is determined by the formula 
\[
 (f\cdot z)(x)= f(x)z
 \qquad 
 \left(x, \in R, z\in Z(R)\right), 
\]
where $Z(R)$ denotes the center\index{center (of a ring)} of the ring $R$. 

At the same time, one cannot claim that $\mathrm{Der}(R)$ is an algebra over the ring 
$Z(R)$, since the definition of an algebra includes the identity 
\[
 \left[f\cdot z, g\right]= \left[f, g\right]\cdot z 
 \qquad 
 \left(f, g\in \mathrm{Der}(R), z\in Z(R)\right), 
\]
while in the ring $\mathrm{Der}(R)$ the following identity holds instead
\[
 \left[f\cdot z, g\right]= \left[f, g\right]\cdot z+f\cdot g(z)
 \qquad 
 \left(f, g\in \mathrm{Der}(R), z\in Z(R)\right). 
\]

Among derivations one can single out so-called inner derivations, similarly as in the case of automorphisms. 

\begin{defin}
 Let $R$ be a ring and $b\in R$, then the mapping $\mathrm{ad}_{b}\colon R\to R$ defined by 
 \[
  \mathrm{ad}_{b}(x)=\left[x, b\right] 
  \qquad 
  \left(x\in R\right)
 \]
is a derivation. A derivation $f\colon R\to R$ is termed to be an \emph{inner derivation}\index{derivation!--- inner} if there is a 
$b\in R$ so that $f=\mathrm{ad}_{b}$. We say that a derivation is an \emph{outer derivation}\index{derivation!--- outer} if it is not inner. 
\end{defin}

Clearly, commutative rings admit only \emph{trivial} inner derivations.

A fundamental example for derivations is the following.

\begin{ex}
Let $\mathbb{F}$ be a field, and let in the above definition $P=Q=\mathbb{F}[x]$
be the ring of polynomials with coefficients from $\mathbb{F}$. For a polynomial
$p\in\mathbb{F}[x]$, $p(x)=\sum_{k=0}^{n}a_{k}x^{k}$, define the function
$f\colon \mathbb{F}[x]\rightarrow\mathbb{F}[x]$ as
\[
f(p)=p',
\]
where $p'(x)=\sum_{k=1}^{n}ka_{k}x^{k-1}$ is the derivative of the polynomial $p$.
Then the function $f$ clearly fulfills
\[
f(p+q)=f(p)+f(q)
\]
and
\[
f(pq)=pf(q)+qf(p)
\]
for all $p, q\in\mathbb{F}[x]$. Hence $f$ is a derivation.
\end{ex}

\begin{ex}
 Let $(\mathbb{F}, +, \cdot)$ be a field, and suppose that we are given a derivation 
 $f\colon \mathbb{F}\to \mathbb{F}$. We define the mapping $f_{0}\colon \mathbb{F}[x]\to \mathbb{F}[x]$ in the following way. 
If $p\in \mathbb{F}[x]$ has the form 
\[
 p(x)=\sum_{k=0}^{n}a_{k}x^{k}, 
\]
then let 
\[
 f_{0}(p)= p^{f}(x)=\sum_{k=0}^{n}f(a_{k})x^{k}. 
\]
Then $f_{0}\colon \mathbb{F}[x]\to \mathbb{F}[x]$ is a derivation. 

Indeed, let $q\in \mathbb{F}[x]$, $q(x)=\displaystyle\sum_{k=0}^{n}b_{k}x^{k}$. Note that without the loss of generality we may assume 
that both sum goes from zero to $n$, adding to that of the smaller degree terms with coefficients zero. 
In this case 
\[
 f_{0}(p+q)=\sum_{k=0}^{n}f(a_{k}+b_{k})x^{k}= \sum_{k=0}^{n}f(a_{k})x^{k}+\sum_{k=0}^{n}f(b_{k})x^{k}=f_{0}(p)+f_{0}(q), 
\]
thus $f_{0}$ is additive. 

Since $pq=\displaystyle\sum_{k=0}^{2n}\left(\sum_{i=0}^{k}a_{i}b_{k-i}\right)x^{k}$, due to the additivity of $f_{0}$, we have 
\begin{multline*}
 f_{0}(pq)= f_{0}\left(\sum_{k=0}^{2n}\left(\sum_{i=0}^{k}a_{i}b_{k-i}\right)x^{k}\right)=
 \sum_{k=0}^{2n}f\left(\sum_{i=0}^{k}a_{i}b_{k-i}\right)x^{k}=
 \sum_{k=0}^{2n}\left(\sum_{i=0}^{k}f(a_{i}b_{k-i})\right)x^{k}\\
 =
 \sum_{k=0}^{2n}\left(\sum_{i=0}^{k}f(a_{i})b_{k-i}+a_{i}f(b_{k-i})\right)x^{k}=
 \sum_{k=0}^{2n}\left(\sum_{i=0}^{k}a_{i}f(b_{k-i})\right)x^{k}
 +\sum_{k=0}^{2n}\left(\sum_{i=0}^{k}a_{i}f(b_{k-i})\right)x^{k}
 \\= pf_{0}(q)+f_{0}(p)q. 
\end{multline*}
\end{ex}

The following lemma says that the above two examples have rather fundamental importance. 

\begin{lem}\label{L2.1.1}
 Let $(\mathbb{K}, +, \cdot)$ be a field and let $(\mathbb{F}, +, \cdot)$ be a subfield of $\mathbb{K}$. 
 If $f\colon \mathbb{F}\to \mathbb{K}$ is a derivation, then for any $a\in \mathbb{F}$ and for arbitrary 
 polynomial $p\in \mathbb{F}[x]$ we have 
 \[
  f(p(a))=  p^{f}(a)+f(a)p'(a). 
 \]
\end{lem}
\begin{biz}
 First we will show that for any $k\in \mathbb{N}$ 
 \[
  f(a^{k})=ka^{k-1}f(a). 
 \]
This is evident for $k=1$. Assuming that the above identity is valid for some $k\in \mathbb{N}$, 
\[
 f(a^{k+1})=f(a\cdot a^{k})= f(a)a^{k}+af(a^{k+1})= a^{k}f(a)+a\cdot k a^{k} f(a)= (k+1)a^{k}f(a). 
\]
Therefore, the above identity holds for any $k\in \mathbb{N}$. 

From the Leibniz rule we have 
\[
 f(1\cdot 1)= 1\cdot f(1)+f(1)\cdot 1= 2f(1), 
\]
implying $f(1)=0$. 

Let now $p(x)= \displaystyle\sum_{k=0}^{n}a_{k}x^{k}$, then 
\begin{multline*}
 f(p(a))= f\left(\sum_{k=0}^{n}a_{k}a^{k}\right)= 
 \sum_{k=0}^{n}f\left(a_{k}a^{k}\right)=
 \sum_{k=0}^{n}f(a_{k})a^{k}+\sum_{k=0}^{n}a_{k}f(a^{k})\\
 = p^{f}(a)+\sum_{k=0}^{n}a_{k}ka^{k-1}f(a)= p^{f}(a)+f(a)p'(a), 
\end{multline*}
which proves the lemma. 
\end{biz}

From the previous lemma we know that $f(1)=0$. On the other hand every real additive function is $\mathbb{Q}$-homogeneous. Thus 
we have the following. 

\begin{lem}\label{L2.1.2}
 If $f\colon \mathbb{R}\to \mathbb{R}$ is a derivation, then $f(x)=0$ for all $x\in \mathbb{Q}$. 
\end{lem}

This lemma and Lemma \ref{L2.1.1} imply the following statement. 

\begin{lem}
Any derivation on the real line vanishes at each algebraic number. 
\end{lem}
\begin{biz}
 Assume $a\in \mathbb{R}$ is an algebraic number and let $p\in \mathbb{Q}[x]$ the minimal polynomial of $a$. 
 Let $p(x)=\displaystyle\sum_{k=0}^{n}a_{k}x^{k}$. Here $a_{k}\in \mathbb{Q}$ for all $k=0, 1, \ldots, n$. 
 Applying Lemma \ref{L2.1.2}, we have 
 \[
  f(a_{k})=0 
  \qquad 
  \left(k=0, 1, \ldots, n\right). 
 \]
Whence
\[
 p^{f}(x)=\sum_{k=0}^{n}f(a_{k})x^{k}=0. 
\]
In view of Lemma \ref{L2.1.1} this means that 
\[
 f(p(a))= f(a)p'(a). 
\]
Since $p(a)=0$, from this we get that $f(p(a))=0$. So $f(a)$ or $p'(a)$ is zero. 
However, due to the minimality of $p$, inequality $\mathrm{deg}(p')< \mathrm{deg}(p)$ holds, yielding that 
$p'(a)\neq 0$. Thus $f(a)=0$. 
\end{biz}

From the additive property of derivations, we have also the following. 

\begin{theorem}
Let $f\colon \mathbb{R}\to \mathbb{R}$ be a derivation and suppose that at least one of the following is fulfilled:
 \begin{enumerate}[(i)]
  \item $f$ is measurable;
  \item $f$ is bounded above on a set of positive Lebesgue measure;
  \item $f$ is bounded below on a set of positive Lebesgue measure. 
 \end{enumerate}
 Then $f$ is identically zero. 
\end{theorem}

Let $k\in \mathbb{N}$ and $[a, b]\subset \mathbb{R}$.   
Henceforth let 
\[
\mathscr{C}([a, b])
 =
\left\{f\colon [a, b]\to \mathbb{R}\, \vert \, \text{$f$ is continuous on $[a, b]$}\right\}
\]
\begin{multline*}
 \mathscr{C}^{k}([a, b])
 =
\left\{f\colon [a, b]\to \mathbb{R}\, \vert \, \text{$f\in \mathscr{C}([a, b])$  }
\right. 
\\
\left. 
\text{
and $f$ is $k$-times continuously differentiable on $]a, b[$}\right\}
\end{multline*}
and 
\[
 \mathscr{C}^{\infty}([a, b])=\bigcap_{k=1}^{\infty}\mathscr{C}^{k}([a, b]). 
\]

If we endow the above spaces with the pointwise addition and the pointwise multiplication by scalars, then 
the above spaces are not only rings but also algebras over $\mathbb{R}$. 

Furthermore, after an easy calculation, we also have the following. 

\begin{ex}
 Define the mapping $D\colon \mathscr{C}^{1}([a, b])\to \mathscr{C}([a, b])$ through 
 \[
  D(f)= f'
  \qquad 
  \left(f\in \mathscr{C}^{1}([a, b])\right), 
 \]
where $f'$ denotes the derivative of the function $f\in \mathscr{C}^{1}([a, b])$. 
Then the mapping $D\colon \mathscr{C}^{1}([a, b])\allowbreak \to \mathscr{C}([a, b])$ is a derivation. 
\end{ex}

The above example shows that \emph{nomen est omen}, derivations imitate the  
action of the differential operator $D$ defined in the previous example. 
Therefore, in such a manner, derivations can be considered not only from the perspective of algebra but also 
from that of analysis. More precisely, following the monograph \cite{Prz88} of D.~Przeworska-Rolewicz, this topic belongs in fact to the 
area of \emph{algebraic analysis}. 
Following her, this coincides with the theory of right invertible operators 
in linear spaces (without any topology, in general). The abovementioned book is an attempt to formulate 
a common treatment of the calculus and linear differential equations with the help 
of the theory of right invertible operators. 
The research of the author has been supported by the Hungarian Scientific Research Fund
(OTKA) Grant K 111651 and by the ÚNKP-4 New National Excellence Program of the Ministry of Human Capacities.
The research is also supported by the
EFOP-3.6.1-16-2016-00022 project. The project is
co-financed by the European Union and the European
Social Fund. 
With regard to the notion of derivations, the most important part of the monograph \cite{Prz88} is Chapter 6, where the author 
introduces the notion of $D$-algebras. For the readers convenience, in what follows, we briefly report the most principal notions and examples 
from this area. 

\begin{defin}
Let $X$ be a linear space over the field $\mathbb{F}$. A linear operator $D\colon X\to X$ is said to be 
\emph{right invertible}\index{operator!--- right invertible}, if there is a linear operator 
$R\colon X\to X$ with the properties $\mathrm{dom}(R)=X$  and $R(X)\subset \mathrm{dom}(D)$ so that 
\[
 DR= \mathrm{id}, 
\]
where $\mathrm{id}\colon X\to X$ denotes the identity operator. 
The operator $R$ is termed to be a \emph{right inverse} of $D$.
The set of all right invertible operators is denoted by $R(X)$. 
\end{defin}

\begin{defin}
 Let $X$ be a commutative algebra and let $D\in R(X)$. 
 $X$ is  said to be a \emph{$D$-algebra} if 
 \begin{enumerate}[(i)]
  \item $\dim \ker D>0$;
  \item $x, y\in  \mathrm{dom}(D)$ implies that $xy\in \mathrm{dom}(D)$. 
 \end{enumerate}
\end{defin}

Let $X$ be a $D$-algebra and write 
\[
 f_{D}(x, y)= D(xy)-c_{D}\left(xDy+yDx\right) 
 \qquad 
 \left(x, y\in \mathrm{dom}(D)\right), 
\]
where 
\begin{enumerate}[(i)]
 \item $c_{D}$ is a scalar depending only on the operator $D$;
 \item $f_{D}\colon \mathrm{dom}(D)\times \mathrm{dom}(D)\to \mathrm{dom}(D)$ is a symmetric, bilinear mapping, 
 that is said to be a \emph{non-Leibniz component}\index{non-Leibniz component}. 
\end{enumerate}

\begin{ex}
 A $D$-algebra is called a \emph{Leibniz $D$-algebra}\index{Leibniz $D$-algebra} or shortly an \emph{$L$-algebra}\index{$L$-algebra}, if 
 $D$ satisfies 
 \[
  D(xy)= xD(y)+yD(x) 
  \qquad 
  \left(x, y\in \mathrm{dom}(D)\right)
 \]
In such a situation we have $c_{D}=1$ and $f_{D}=0$. 

Moreover, if $X$ has a unit $e$, then $D(e)=0$ and we also have the following identity 
\[
 f^{(n)}_{D}(x, y)=\sum_{k=1}^{n-1}\binom{n}{k} \left(D^{k}(x)\right)\left(D^{n-k}(y)\right) 
 \qquad 
 \left(x, y\in \mathrm{dom}(D^{n})\right)
\]
for any  $n\geq 2$. 
\end{ex}

\begin{ex}
 Let us consider $X=\mathscr{C}([a, b])$ with the pointwise operations. Then 
 $X$ is a Leibniz $\dfrac{d}{dt}$-algebra. 
 This algebra has a unit, namely the function 
 \[
  e(t)=1 
  \qquad 
  \left(t\in [a, b]\right), 
 \]
as well as zero divisors. 
\end{ex}

\begin{ex}
 Let $\Omega=\left\{(t, s)\in\mathbb{R}^{2}\, \vert \, t\in [a, b], s\in [c, d]\right\}$ and 
 $X=\mathscr{C}(\Omega)$ with the pointwise operations. 
 Let further 
 \[
D_{1}= \dfrac{\partial}{\partial t}   
\qquad 
\text{and}
\qquad 
D_{2}= \dfrac{\partial}{\partial s}. 
 \]
Then $X$ is a Leibniz $D_{1}$-algebra and simultaneously a Leibniz $D_{2}$-algebra. 
\end{ex}

Although, in this work we focus mainly on the algebraic nature of derivations, we have to point out, that in the abovementioned monograph among 
others the following notions are also introduced: quasi-Leibniz $D$-algebra, Duhamel $D$-algebra, simple Duhamel-algebra, almost Leibniz $D$-algebra. 

After presenting these notions, several characterization results are shown. Furthermore, the connection between these notions is also investigated.

\section{Extensions of derivations}

In the present section we study the possibility of extending a derivation
from its domain of definition onto a larger algebraic structure.
\begin{lem}\label{L14.2.1}
Let $(P, +,\cdot)$ be an integral domain, let $(\mathbb{F}, +,\cdot)$ be its field of
fractions, and let $(\mathbb{K}, +,\cdot)$ be a field so that $P\subset \mathbb{F}\subset \mathbb{K}$ is fulfilled.
If $f\colon P\to \mathbb{K}$ is a derivation, then there exists a unique derivation
$g\colon \mathbb{F}\to \mathbb{K}$ such that $g\vert_{P}=f$.
\end{lem}
\begin{biz}
Every $x\in \mathbb{F}$ can be written as $x=u/v$, where $u,v\in P$, $v\neq 0$.
For such an $x$ we put
\begin{equation}\label{E14.2.1}
    g(x)=g\left(\dfrac{u}{v}\right)=\dfrac{vf(u)-uf(v)}{v^{2}}.
\end{equation}
We must check that  this definition is unambiguous, i.e.  if
$u/v=z/w$, $v,w \neq 0$, then $g(u/v)=g(z/w)$. Now, $u/v=z/w$ means
\begin{equation}\label{E14.2.2}
    uw=vz,
\end{equation}
whence by the Leibniz rule $uf(w)+wf(u)=vf(z)+zf(v)$, i.e. 
    \[
    vf(z)-uf(w)=wf(u)-zf(v).
\]
Multiplying this by $uw$ we get
    \[
    v^{2}wf(z)-uvwf(w)=v w^{2}f(u)-zvwf(v),
\]
or, by \eqref{E14.2.2}
    \[
    v^{2}wf(z)-v^{2} zf(w)=vw^{2}f(u)-uw^{2}f(v).
\]
Dividing this by $v^{2}w^{2}$ we obtain
    \[
    \dfrac{wf(z)-zf(w)}{w^{2}}=\dfrac{vf(u)-uf(v)}{v^{2}},
\]
i.e.  \ $g(z/w)=g(u/v)$. Thus expression \eqref{E14.2.1} does not depend on the
representation of $x$ as a fraction $u/v$.

Now take arbitrary $x,y\in \mathbb{F}$, $x=u/v$, $y=z/w$, $u,v,z,w\in P$, $v,w\neq 0$.
Then using the additivity and the Leibniz rule, 
\begin{align*}
    g(x+y)&=g\left(\dfrac{u}{v}+\dfrac{z}{w}\right)=
    g\left(\dfrac{uw+zv}{vw}\right)=
    \dfrac{vwf(uw+zv)-(uw+zv)f(uw)}{v^{2}w^{2}}\\
    &=\dfrac{vw\big(uf(w)+wf(u)+zf(v)+vf(z)\big)-(uw+zv)\big(vf(w)+wf(v)\big)}
    {v^{2}w^{2}}\\
    &=\dfrac{vw^{2}f(u)+v^{2}wf(z)-uw^{2}f(v)-zv^{2}f(w)}{v^{2}w^{2}}\\
    &=\dfrac{vf(u)-uf(v)}{v^{2}}+\dfrac{wf(z)-zf(w)}{w^{2}}=
    g\left(\dfrac{u}{v}\right)+g\left(\dfrac{z}{w}\right)\\
    &=g(x)+g(y)\,,
\end{align*}
and
\begin{align*}
    g(xy)&=g\left(\dfrac{uz}{vw}\right)=
    \dfrac{vwf(uz)-uzf(vw)}{v^{2}w^{2}}\\
    &=\dfrac{vw\big(uf(z)+zf(u)\big)-uz\big(vf(w)+wf(v)\big)}{v^{2}w^{2}}\\
    &=\dfrac{uv\big(wf(z)-zf(w)\big)+zw\big(vf(u)-uf(v)\big)}{v^{2}w^{2}}\\
    &=\dfrac{u}{v}\dfrac{wf(z)-zf(w)}{w^{2}}+\dfrac{z}{w}\dfrac{vf(u)-uf(v)}{v^{2}}\\
    &=\dfrac{u}{v}g\left(\dfrac{z}{w}\right)+\dfrac{z}{w}g\left(\dfrac{u}{v}\right)=
    xg(y)+yg(x).
\end{align*}
Consequently $g$ is a derivation of $\mathbb{F}$. If $x\in P$, then $x=x/1$ and bearing in mind that $f(1)=0$, 
    \[
    g(x)=g\left(\dfrac{x}{1}\right)=\dfrac{f(x)-xf(1)}{1}=f(x),
\]
i.e.  $g\vert_{P}=f$.

Now let $g\colon \mathbb{F}\to \mathbb{K}$ be an arbitrary derivation such that $g\vert_{P}=f$.
Take an arbitrary $x\in \mathbb{F}$, $x=u/v$, $u,v\in P$, $v\neq 0$.
We have $f(u)=g(u)=g(vx)=vg(x)+xg(v)=vg(x)+xf(v)$, whence
    \[
    g(x)=\dfrac{f(u)-xf(v)}{v}=\dfrac{vf(u)-uf(v)}{v^{2}}.
\]
This means that $g$ has form \eqref{E14.2.1},
which proves the uniqueness of the extension.
\end{biz}

\begin{lem}\label{L14.2.2}
Let $(\mathbb{K}, +,\cdot)$ be a field, let $(\mathbb{F}, +,\cdot)$ be a subfield of $(\mathbb{K}, +,\cdot)$,
and let $f\colon \mathbb{F}\to \mathbb{K}$ be a derivation. Further, let $a,u\in \mathbb{K}$.
There exists a derivation $g\colon \mathbb{F}(a)\to \mathbb{K}$ such that $g\vert_{\mathbb{F}}=f$ and $g(a)=u$
if and only if
\begin{equation}\label{E14.2.3}
    r^{f}(a)+ur'(a)=0
\end{equation}
for every $r\in \mathbb{F}\left[x\right]$ such that $r(a)=0$. If it exists, the extension $g$ is unique.
\end{lem}

\begin{lem}\label{L14.2.3}
Let $(\mathbb{K}, +,\cdot)$ be a field, let $(\mathbb{F}, +,\cdot)$ be a subfield of
$(\mathbb{K}, +,\cdot)$ and let $S\subset \mathbb{K}$ algebraically independent over $\mathbb{F}$.
Let $f\colon \mathbb{F}\to \mathbb{K}$ be a derivation, and let $u\colon S\to \mathbb{K}$ be an
arbitrary function. Then there exists a unique derivation $g\colon \mathbb{F}(S)\to \mathbb{K}$
such that $g\vert_{\mathbb{F}}=f$ and $g\vert_{S}=u$.
\end{lem}
\begin{biz}
Let $\mathcal{R}$ be the collection of all couples $(S_{\alpha},g_{\alpha})$ such that
$S_{\alpha}\subset S$, $g_{\alpha}\colon \mathbb{F}(S_{\alpha})\to \mathbb{K}$ is a derivation,
and $g_{\alpha}\vert_{\mathbb{F}}=f$, $g_{\alpha}\vert_{S_{\alpha}}=u$.
The couple $(\emptyset, f)\in\mathcal{R}$ so that $\mathcal{R}\neq\emptyset$.
We order $\mathcal{R}$ as follows: $(S_{\alpha},g_{\alpha})\leqslant (S_{\beta},g_{\beta})$ if and only if
$S_{\alpha}\subset S_{\beta}$ and $g_{\beta}\vert_{\mathbb{F}(S_{\alpha})}=g_{\alpha}$.
Thus $(\mathcal{R},\leqslant)$ in an ordered set, and if $\mathcal{L}\subset\mathcal{R}$
is a chain, then the couple $(S_0,g_0)$ such that
$S_0=\bigcup\limits_{(S_{\alpha},g_{\alpha})\in\mathcal{L}}S_{\alpha}$, $g_0\vert_{F(S_{\alpha})}=g_{\alpha}$
for $(S_{\alpha},g_{\alpha})\in\mathcal{L}$, is an upper bound of $\mathcal{L}$ in $\mathcal{R}$.
By the Lemma of Kuratowski--Zorn there exists a maximal element $(S_{\max},g_{\max})$ in $\mathcal{R}$.
Thus, in particular, $S_{\max}\subset S$. Suppose that there exists an $a\in S\setminus S_{\max}$.
Since $S$ is algebraically independent over $\mathbb{F}$ and $a\in S\setminus S_{\max}$, we have
$r(a)\neq 0$ for every $r\in \mathbb{F}(S_{\max})[x]$.
So the condition in Lemma \ref{L14.2.2} is trivially fulfilled. 
Thus the derivation $g_{\max}$ can be extended onto $\mathbb{F}(S_{\max})(a)$ to a derivation
$g^{*}\colon \mathbb{F}(S_{\max})(a)\to \mathbb{K}$ such that $g^{*}\vert_{\mathbb{F}(S_{\max})}=g_{\max}$, $g^{*}(a)=u(a)$.
Hence $g^{*}$ satisfies $g^{*}\vert_{\mathbb{F}}=g^{*}\vert_{\mathbb{F}(S_{\max})}\vert_{\mathbb{F}}=g_{\max}\vert_{\mathbb{F}}=f$,
$g^{*}\vert_{S_{\max}}=g_{\max}\vert_{S_{\max}}=u$, whence
$g^{*}\vert_{\big(S_{\max}\cup\{a\}\big)}=u$.
Writing $S^{*}=S_{\max}\cup \left\{a\right\}$ we obtain hence that
$\left(S^{*},g^{*}\right)\in\mathcal{R}$ and $(S_{\max},g_{\max})<(S^{*},g^{*})$,
which contradicts the maximality of $(S_{\max},g_{\max})$.

Consequently $S_{\max}=S$, and $g=g_{\max}$ is the required extension.

It remains to prove the uniqueness. Let $x\in \mathbb{F}(S)$.
Then there exists a finite
set $S_1=\left\{a_1,\ldots,a_n\right\}\subset S$ such that $x\in \mathbb{F}(S_1)$. But
    \[
    \mathbb{F}(S_1)=\mathbb{F}(a_1)(a_2)\cdots (a_n).
\]
Since, by Lemma \ref{L14.2.2}, the extension of a derivation onto a simple extension
of its field of definition (with the prescribed value at the element by which we
extend the basic field) is unique, $g$ is uniquely determined on $\mathbb{F}(a_1)\dots (a_k)$
for every $k=1,\ldots,n$, and so, in particular, $g(x)$ is uniquely determined.
Thus $g$ is uniquely determined at every point $x\in \mathbb{F}(S)$, and so $g$ is unique.
\end{biz}
\begin{lem}\label{L14.2.4}
Let $(\mathbb{K}, +,\cdot)$ be a field of characteristic zero,
let $(\mathbb{K}_{1}, +,\cdot)$ be a subfield of $(\mathbb{K}, +,\cdot)$, and let
$(\mathbb{F}, +,\cdot)$ be a subfield of $(\mathbb{K}_{1}, +,\cdot)$ such that
$\mathbb{K}_{1}\subset \operatorname{alg\,cl}(\mathbb{F})$. Let $f\colon \mathbb{F}\to \mathbb{K}$ be a derivation.
Then there exists a unique derivation $g\colon \mathbb{K}_{1}\to \mathbb{K}$ such that
$g\vert_{\mathbb{F}}=f$.
\end{lem}
\begin{biz}
Let $\mathcal{R}$ be the collection of all couples $(\mathbb{K}_\alpha,g_\alpha)$
such that $(\mathbb{K}_\alpha, +,\cdot)$ is a subfield of $(\mathbb{K}_{1}, +,\cdot)$ ,
$\mathbb{F}\subset \mathbb{K}_\alpha$, $g_\alpha\colon \mathbb{F}_\alpha\to \mathbb{K}$ is a derivation
and $g_\alpha \vert_{\mathbb{F}}=f$. $(\mathbb{F},f)\in\mathcal{R}$, so $\mathcal{R}\neq\emptyset$.
We order $\mathcal{R}$ similarly as in the previous proof. 
$(\mathbb{K}_\alpha,g_\alpha)\leqslant (\mathbb{K}_\beta,g_\beta)$ if and only if $(\mathbb{K}_\alpha, +, \cdot)$ is a subfield of
$(\mathbb{K}_\beta, +,\cdot)$ and $g_\beta \vert_{\mathbb{K}_\alpha}=g_\alpha$. Then $(\mathcal{R},\leqslant)$
is an ordered set, and as previously we verify that every chain in $\mathcal{R}$
has an upper bound $(K_0,g_0)\in\mathcal{R}$. 
Again, be the Lemma of Kuratowski--Zorn, 
there exists in $\mathcal{R}$ a maximal element $(\mathbb{K}_{\max},g_{\max})$. 
In particular, $(\mathbb{K}_{\max}, +,\cdot)$ is a subfield of $(\mathbb{K}_{1}, +,\cdot)$. 
Suppose that there exists an 
$a\in \mathbb{K}_1\setminus \mathbb{K}_{\max}$. 
Thus $a$ is algebraic over $\mathbb{F}$, and let $p\in \mathbb{F}\left[x\right]$
be its minimal polynomial. 
We have $p'\neq 0$, since the characteristic of $\mathbb{K}$ is zero
(whence also the characteristic of $\mathbb{F}$ is zero), 
whence $p'(a)\neq 0$, since
$\operatorname{degree}p'< \operatorname{degree}p$, and $p$ is the minimal polynomial of $a$.
Put $u=-p^{g_{\max}}(a)/p'(a)$, and let $r\in \mathbb{K}_{\max}[x]$ be a polynomial such that $r(a)=0$.
There exist polynomials $q,s \in \mathbb{K}_{\max}[x]$ such that
$r=qp+s$ and $\operatorname{degree}s< \operatorname{degree}p$. 
Hence $0=r(a)=q(a)p(a)+s(a)=s(a)$, whence it
follows that $s=0$, and $r=qp$. Now
\begin{align*}
    r^{g_{\max}}(a)+ur'(a)&=
    q(a)p^{g_{\max}}(a)+p(a)q^{g_{\max}}(a)+
    up'(a)q(a)+up(a)q'(a)=\\
    &=q(a)\big(p^{g_{\max}}(a)+up'(a)\big),
\end{align*}
since $p(a)=0$. By the choice of $u$ we have $p^{g_{\max}}(a)+up'(a)=0$, whence also
$r^{g_{\max}}(a)+ur'(a)=0$.

By Lemma \ref{L14.2.2} there exists a derivation $g^{*}\colon \mathbb{K}_{\max}(a)\to \mathbb{K}$
such that $g^{*} \vert_{\mathbb{K}_{\max}}=g_{\max}$. 
Write $\mathbb{K}^{*}=\mathbb{K}_{\max}(a)$. Since $\mathbb{K}_{\max}\subset \mathbb{K}_{1}$
and $a\in \mathbb{K}_{1}$, the field $(\mathbb{K}^{*}, +,\cdot)$ is a subfield of $(\mathbb{K}_{1}, +,\cdot)$, 
and, of course, $\mathbb{F}\subset \mathbb{K}_{\max}\subset \mathbb{K}^{*}$. 
Moreover, $g^{*} \vert_{\mathbb{F}}=g_{\max}\vert_{\mathbb{F}}=f$. Thus $(\mathbb{K}^{*},g^{*})\in \mathcal{R}$,
and $(\mathbb{K}_{\max},g_{\max})<(\mathbb{K}^{*},g^{*})$, 
which contradicts the maximality of $(\mathbb{K}_{\max},g_{\max})$.
Hence $\mathbb{K}_{\max}=\mathbb{K}_{1}$, and $g=g_{\max}$ is the desired extension.

To prove uniqueness, suppose that $g\colon \mathbb{K}_{1}\to \mathbb{K}$ is a derivation such that $g \vert_{\mathbb{F}}=f$.
Take an $a\in \mathbb{K}_{1}$. 
Let $p\in \mathbb{F}\left[x\right]$ be the minimal polynomial of $a$.
As we have seen above, we have $p'(a)\neq 0$. By Lemma \ref{L2.1.1}
    \[
    p^{f}(a)+g(a)p'(a)=
    p^{g}(a)+g(a)p'(a)=
    g\big(p(a)\big)=g(0)=0.
\]
Hence $g(a)=-p^{f}(a)/p'(a)$ is uniquely determined. Thus all the values of $g$
on $K_1$ are uniquely determined, whence $g$ is unique.
\end{biz}

\begin{theorem}\label{T14.2.1}
Let $(\mathbb{K}, +,\cdot)$ be a field of characteristic zero, let $(\mathbb{F}, +,\cdot)$
be a subfield of $(\mathbb{K}, +,\cdot)$, let $S$ be an algebraic base of $\mathbb{K}$ over $\mathbb{F}$,
if it exists, and let $S=\emptyset$ otherwise.
Let $f\colon \mathbb{F}\to \mathbb{K}$ be a derivation.
Then, for every function $u\colon S\to \mathbb{K}$,
there exists a unique derivation $g\colon \mathbb{K}\to \mathbb{K}$
such that $g \vert_{\mathbb{F}}=f$ and $g \vert_{S}=u$.
\end{theorem}

\begin{biz}
If $\mathbb{K}=\operatorname{alg\,cl} (\mathbb{F})$, this results from Lemma \ref{L14.2.4}.
(Then $S=\emptyset$, and there is no $u$ involved).
If $\mathbb{K}\neq \operatorname{alg\,cl} (\mathbb{F})$, then there exists an algebraic base
$S$ of $\mathbb{K}$ over $\mathbb{F}$ so that $\operatorname{alg\,cl} (\mathbb{F}(S))=\mathbb{K}$. 
Since $S$ is algebraically independent over $\mathbb{F}$,
by Lemma \ref{L14.2.3} there exists a unique derivation $g_0\colon \mathbb{F}(S)\to \mathbb{K}$
such that $g_0 \vert_{\mathbb{F}}=f$ and $g_0 \vert_{S}=u$. By Lemma \ref{L14.2.4} the derivation
$g_0$ can be uniquely extended onto $\operatorname{alg\,cl} (\mathbb{F}(S))=\mathbb{K}$ to a derivation $g\colon \mathbb{K}\to \mathbb{K}$
such that  $g \vert_{\mathbb{F}(S)}=g_0$, whence $g \vert_{\mathbb{F}}={g \vert_{\mathbb{F}(S)}}_{\vert \mathbb{F}}= g_0 \vert_{\mathbb{F}} =f$,
and $g \vert_{S}= {g \vert_{\mathbb{F}(S)}}_{\vert S}=g_0 \vert_{S}=u$.
\end{biz}

\begin{theorem}\label{T14.2.2}
There exist non-trivial derivations of $\mathbb{R}$.
\end{theorem}

\begin{biz}
We have $\mathbb{R}\neq \operatorname{alg\,cl}(\mathbb{Q})$, so there exists an
algebraic base of $\mathbb{R}$ over $\mathbb{Q}$. Let $u\colon S \to\mathbb{R}$ be an arbitrary function,
$u\neq 0$. The characteristic of $\mathbb{R}$ is zero. Now take in Theorem \ref{T14.2.1}
$\mathbb{F}=\mathbb{Q}$ and $\mathbb{K}=\mathbb{R}$. 
The trivial derivation $f_0\colon \mathbb{Q}\to\mathbb{R}$, $f_0=0$, can, 
by Theorem \ref{T14.2.1}, be uniquely extended onto $\mathbb{R}$ to a derivation
$f\colon \mathbb{R}\to\mathbb{R}$ such that $f\vert S=u$, whence $f\neq 0$.
\end{biz}
Incidentally, in this way we have obtained a description of all the derivations of $\mathbb{R}$.
Every such derivation can be arbitrary prescribed on an algebraic base $S$ of $\mathbb{R}$ over $\mathbb{Q}$,
and then it is already uniquely determined.

\section{Relations between additive functions}

Let $f\colon \mathbb{R}\to\mathbb{R}$ be a derivation, and let $x\in\mathbb{R}\setminus \left\{0\right\}$.
We have by Lemma \ref{L2.1.2}
    \[
    0= f(1)= f\left(x \dfrac{1}{x}\right)=
    x f\left(\dfrac{1}{x}\right)+\dfrac{1}{x} f(x),
\]
whence
\begin{equation}\label{E14.3.1}
    f(x)= -x^{2} f\left(\dfrac{1}{x}\right) .
\end{equation}
We will show that this relation characterizes derivations among additive
functions
$f\colon \mathbb{R}\to\mathbb{R}$.

\begin{theorem}\label{T14.3.1}
Let $f\colon \mathbb{R}\to\mathbb{R}$ be an additive function satisfying condition \eqref{E14.3.1} for all
$x\in \mathbb{R}\setminus \left\{0\right\}$. Then $f$ is a derivation.
\end{theorem}

\begin{biz}
Take an arbitrary $x\in\mathbb{R}$, $x\neq 0, 1, -1$. Then we have by \eqref{E14.3.1} and and by the additivity of  $f$ that 
\begin{align*}
    f(x)+\dfrac{1}{x^{2}} f(x)&=
    f(x)-f\left(\dfrac{1}{x}\right)=
    f\left(x-\dfrac{1}{x}\right)=
    f\left(\dfrac{x^{2}-1}{x}\right)\\
    &=- \left(\dfrac{x^{2}-1}{x}\right)^{2} f\left(\dfrac{x}{x^{2}-1}\right)=
    -\left(\dfrac{x^{2}-1}{x}\right)^{2} f\left(\dfrac{1}{x-1}-\dfrac{1}{x^{2}-1}\right)\\
    &=- \left(\dfrac{x^{2}-1}{x}\right)^{2} f\left(\dfrac{1}{x-1}\right)+
    \left(\dfrac{x^{2}-1}{x}\right)^{2} f\left(\dfrac{1}{x^{2}-1}\right)\\
    &= \left(\dfrac{x^{2}-1}{x}\right)^{2}\left(\dfrac{1}{x-1}\right)^{2} f(x-1)-
    \left(\dfrac{x^{2}-1}{x}\right)^{2}\left(\dfrac{1}{x^2 -1}\right)^{2} f\left(x^{2}-1\right)\\
    &=\left(\dfrac{x+1}{x}\right)^{2} f(x-1)- \dfrac{1}{x^{2}}f\left(x^{2}-1\right).
\end{align*}
Setting in \eqref{E14.3.1} $x=1$ or $-1$, we obtain
\begin{equation}\label{E14.3.2}
    f(1)=f(-1)=0.
\end{equation}
Hence $f\left(u-1\right)=f\left(u\right)-f\left(1\right)=f\left(u\right)$ for every
$u\in\mathbb{R}$, and we get
    \[
    f(x)+\dfrac{1}{x^{2}}f(x)=
    \left(\dfrac{x+1}{x}\right)^{2} f(x)-\dfrac{1}{x^{2}} f\left(x^{2}\right),
\]
    \[
    x^{2}f(x)+f(x)= \left(x+1\right)^{2}f(x)-f\left(x^{2}\right),
\]
i.e. 
\begin{equation}\label{E14.3.3}
    f\left(x^{2}\right)=2xf(x).
\end{equation}
Again, due to the additivity we have $f(0)=0$. 
Thus, in view of \eqref{E14.3.2}, relation \eqref{E14.3.3}, so
far obtained for $x\neq 0, 1, -1$, is valid for all $x\in\mathbb{R}$. We have by \eqref{E14.3.3} for
every $x, y\in\mathbb{R}$
    \[
    f\big((x+y)^{2}\big)=2\left(x+y\right)f(x+y),
\]
i.e. 
    \[
    f\left(x^2\right)+2f\left(xy\right)+f\left(y^2\right)=
    2xf(x)+2xf(y)+2yf(x)+2yf(y),
\]
and by \eqref{E14.3.3}
\begin{equation}\label{E14.3.4}
    2xf(x)+2f\left(xy\right)+2yf(y)=
    2xf(x)+2xf(y)+2yf(x)+2yf(y).
\end{equation}
Relation \eqref{E14.3.4} yields $f$ satisfies the Leibniz rule, i.e.  $f$ is a derivation.
\end{biz}

\begin{lem}\label{L14.3.1}
Let $f\colon \mathbb{R}\to\mathbb{R}$ and $g\colon \mathbb{R}\to\mathbb{R}$ be additive functions, $f\neq 0$, and let
$P\colon \mathbb{R}\to\mathbb{R}$ be a continuous function. If
\begin{equation}\label{E14.3.5}
    g(x)=P(x)f\left(\dfrac{1}{x}\right)
\end{equation}
for all $x\in\mathbb{R}$, $x\neq 0$, then
\begin{equation}\label{E14.3.6}
    P(x)=P(1)x^{2}
\end{equation}
for all $x\in\mathbb{R}$.
\end{lem}

\begin{theorem}\label{T14.3.2}
Let $f\colon \mathbb{R}\to\mathbb{R}$ and $g\colon \mathbb{R}\to\mathbb{R}$ be additive functions satisfying the relation
\begin{equation}\label{E14.3.8}
    g(x)=x^{2}f\left(\dfrac{1}{x}\right)
\end{equation}
for all $x\in\mathbb{R}$, $x\neq 0$. Then the function $f+g$ is continuous, and the functions
$F, G\colon \mathbb{R}\to\mathbb{R}:$
\begin{equation}\label{E14.3.9}
    F(x)=f(x)-xf(1), \quad
    G(x)=g(x)-xg(1)
\end{equation}
are derivations.
\end{theorem}

\begin{biz}
We have by \eqref{E14.3.8} $g\left(1\right)=f\left(1\right)$, whence by
\eqref{E14.3.8} and \eqref{E14.3.9}, for $x\neq 0$,
    \[
    x^{2}F\left(\dfrac{1}{x}\right)=
    x^{2}f\left(\dfrac{1}{x}\right)-xf(1)=
    g(x)-xg(1)=G\left(x\right)
\]
so that
\begin{equation}\label{E14.3.10}
    G\left(x\right)=x^{2}F\left(\dfrac{1}{x}\right)
\end{equation}
for every $x\neq 0$. Now, $F\left(1\right)=G\left(1\right)=0$ by \eqref{E14.3.9},
whence, since $F$ and $G$ clearly are additive,
$F\left(x+1\right)=F\left(x\right)+F\left(1\right)=F\left(x\right)$ and
$G\left(x+1\right)=G\left(x\right)+G\left(1\right)=G\left(x\right)$ for every $x\in\mathbb{R}$.
Hence and by \eqref{E14.3.10} we have for every $x\neq -1$
\begin{align*}
    G\left(x\right)&= G\left(x+1\right)= \left(x+1\right)^{2} F\left(\dfrac{1}{x+1}\right)=
    \left(x+1\right)^{2} F\left(1-\dfrac{x}{x+1}\right)\\
    &=-\left(x+1\right)^{2}F\left(\dfrac{x}{x+1}\right)=
    -\left(x+1\right)^{2}\left(\dfrac{x}{x+1}\right)^{2} G\left(\dfrac{x+1}{x}\right)\\
    &=-x^{2}G\left(1+\dfrac{1}{x}\right)=-x^{2}G\left(\dfrac{1}{x}\right)=-F\left(x\right),
\end{align*}
i.e. 
\begin{equation}\label{E14.3.11}
    G\left(x\right)=-F\left(x\right).
\end{equation}
For $x=-1$ we have $G\left(-1\right)=-G\left(1\right)=0=F\left(1\right)=-F\left(-1\right)$,
so \eqref{E14.3.11} holds for $x=-1$, too, and thus \eqref{E14.3.11} is valid for all
$x \in \mathbb{R}$. \eqref{E14.3.11} yields $g(x)-xg(1)=-f(x)+xf(1)$, whence
    \[
    f(x)+g(x)=x\big(f(1)+g(1)\big),
\]
and consequently $f+g$ is a continuous function. Further, we get by \eqref{E14.3.10} and \eqref{E14.3.11}
    \[
    F\left(x\right)=-x^{2}F\left(\dfrac{1}{x}\right)
\]
for $x\neq 0$. By Theorem \ref{T14.3.1} $F$ is a derivation, and by Proposition \ref{P2.1.1} $G=-F$ also
is a derivation.
\end{biz}

\begin{theorem}\label{T14.3.3}
Let $f\colon \mathbb{R}\to\mathbb{R}$ be an additive function, and let $P\colon \mathbb{R}\to\mathbb{R}$ be a continuous function.
If
\begin{equation}\label{E14.3.12}
    f(x)=P\left(x\right)f\left(\dfrac{1}{x}\right)
\end{equation}
for all $x\in\mathbb{R}$, $x\neq 0$, then if $P\left(1\right)=-1$, $f$ is a derivation,
and if $P\left(1\right)\neq -1$, then $f$ is continuous.
\end{theorem}
\begin{biz}
By Lemma \ref{L14.3.1}
\begin{equation}\label{E14.3.13}
    P\left(x\right)=cx^{2}, \quad x\in\mathbb{R},
\end{equation}
where $c=P\left(1\right)$. If $c=0$, then $f=0$ by \eqref{E14.3.12}, and $f$ is continuous.
If $c\neq 0$, write $g(x)=c^{-1}f(x)$, $x\in\mathbb{R}$. Thus $g\colon \mathbb{R}\to\mathbb{R}$ is an additive function,
and \eqref{E14.3.12} with \eqref{E14.3.13} imply \eqref{E14.3.8}. By Theorem \ref{T14.3.2}
$f+g=\left(1+c^{-1}\right)f$ is continuous, and hence $f$ is continuous if $c\neq -1$.
If $c=-1$, then $P\left(x\right)=-x^{2}$, and by Theorem \ref{T14.3.1} $f$ is a derivation.
\end{biz}

\section{The cocycle equation}\index{cocycle equation}

The cocycle functional equation has a long history in connection with
many  areas  of  mathematics  and  its  applications,  as  discussed  for  example
in Jessen--Karpf--Thorup \cite{JesKarTho68}, where the authors used the hereunder results 
to give a simplified proof of Sydler's theorem on polyhedra. 
From our point of view it will be important that the cocycle equation plays a key role in the 
theory of derivations, as well.

Let $A$ and $X$ be commutative groups. 
Then any solution $F\colon A^{2}\to X$ of the functional equation 
\[
 F(a+b, c)+F(a, b)=F(a, b+c)+F(b, c) 
 \qquad 
 \left(a, b, c\in A\right)
\]
will be called a \emph{cocycle}\index{cocycle} on the group $A$. 
Further, this equation will be referred to as the \emph{cocycle equation}\index{cocycle equation}. 
Concerning the solutions of this equation the reader should consult the following works 
Davison--Ebanks \cite{DavEba95}\index{Davison, T.~M.~K}\index{Ebanks, B.}, Ebanks \cite{Eba79, Eba82}, Erdős \cite{Erd59}\index{Erdős, J. }. 
The results presented here are based on Jessen--Karpf--Thorup \cite{JesKarTho68}. 

\begin{theorem}\label{T1.2.1}
 Let $A$ and $X$ be commutative groups, and let 
 $f\colon A\to X$ be an arbitrary function. 
 Then the function $F\colon A^{2}\to X$ defined by 
 \begin{equation}
\tag{A} F(a, b)=f(a+b)-f(a)-f(b)
\qquad 
\left(a, b\in A\right),
\end{equation}
satisfies equations
\begin{equation}
\tag{$\alpha$} F(a, b)=F(b, a)
\qquad 
\left(a, b\in A\right)
\end{equation}
and 
\begin{equation}
\tag{$\beta$} F(a+b, c)+F(a, b)=F(a, b+c)+F(b, c)
\qquad 
\left(a, b, c\in A\right),
\end{equation}
If $A$ is free or $X$ is divisible, the function $F$ determined by the function $f$ 
through equation $(A)$ is the only function which satisfies equations $(\alpha)$ and 
$(\beta)$. 
\end{theorem}

Equation 
\[
F(a, b)=0 
\qquad 
\left(a, b\in A\right)
\]
expresses that the function $f$ is additive. Thus, 
for an arbitrary function $f$, the function $F$ may be said to measure how much 
$f$ deviates from being additive. 

\begin{theorem}\label{T1.2.2}
Let $A$ be a commutative ring and $X$ be a module over $A$,
and let $f\colon A\rightarrow X$ be an arbitrary function.
Then the functions
$F\colon A^{2}\rightarrow X$ and $G\colon A^{2}\rightarrow X$ defined by
\begin{equation}
\tag{A} F(a, b)=f(a+b)-f(a)-f(b),
\end{equation}
\begin{equation}
\tag{B} G(a, b)=f(ab)-af(b)-bf(a)
\end{equation}
satisfy
\begin{equation}
\tag{$\alpha$} F(a, b)=F(b, a)
\end{equation}
\begin{equation}
\tag{$\beta$} F(a+b, c)+F(a, b)=F(a, b+c)+F(b, c),
\end{equation}
\begin{equation}
\tag{$\gamma$} G(a, b)=G(b, a)
\end{equation}
\begin{equation}
\tag{$\delta$} cG(a, b)+G(ab, c)=aG(b, c)+G(a, bc),
\end{equation}
\begin{equation}
\tag{$\varepsilon$} F(ac, bc)-cF(a, b)=G(a+b, c)-G(a, c)-G(b, c).
\end{equation}
Furthermore, if $A$ has a unity $1$ and $X$ is unitary, then the function
$F$ satisfies the equation
\begin{equation}
\tag{$\zeta$} \sum_{i=1}^{p}F(1, i1)=0, \quad p=\mathrm{char}A.
\end{equation}
If $A$ is an integral domain and $X$ is a unitary module over $A$
which is uniquely $A$-divisible, then the pairs of functions
$F, G$ determined by means of a function $f$ through the equations
$(A), (B)$ are the only pairs of functions which satisfy the system
$(\alpha)-(\zeta)$.
\end{theorem}

Notice that  equation $(\zeta)$ is void if $p=0$. 
The equations $F=0$, $G=0$ express that $f$ is a derivation. 
Thus, for an arbitrary $f$, the pair of functions $F, G$ may be said to measure 
how much $f$ deviates from being a derivation.

\begin{lem}\label{L1.2.1}
 Let $A$ be an ordered commutative group and $X$ be a commutative group, 
 consider further 
 \[
  A_{+}=\left\{a\in A\, \vert \, a>0\right\}. 
 \]
Let $F\colon A^{2}_{+}\to X$ be a function which satisfies equations $(\alpha)$, $(\beta)$ 
for all $a, b, c\in A_{+}$. 
Then the function $F$ can be extended to a function $\widetilde{F}\colon A^{2}\to X$ which satisfies 
$(\alpha)$ and $(\beta)$ for all $a, b, c \in A$. 
\end{lem}
\begin{biz}
 Define the function $\widetilde{F}$ as zero when at least one of the elements 
 $a$, $b$, $a+b$ is zero. Otherwise we define $\widetilde{F}$ according to the following table, where 
 $+$, respectively $-$ stands for $>0$ and $<0$. 
 \begin{center}
 \begin{tabular}{|c|c|c|c|}
  \hline
 $a$ & $b$& $a+b$& $\widetilde{F}(a, b)$ 
 \\ 
 \hline
 \hline
+ &+&+& $F(a, b)$ \\
 \hline
 + &-&+& $-F(a+b, -b)$ \\
 \hline
 + &-&-& $F(-a-b, a)$ \\
 \hline
 - &+&+& $-F(a+b, -a)$ \\
 \hline
 - &+&-& $F(-a-b, b)$ \\
 \hline
 - &-&-& $-F(-a, -b)$ \\
 \hline
 \end{tabular} 
 \end{center}
 One easily verifies that the function $\widetilde{F}\colon A^{2}\to X$ thus defined satisfies the conditions. 
\end{biz}

The verification mentioned in the above proof concerning the definition of the extension 
$\widetilde{F}$ requires the consideration of a number of cases. 
In the paper Jessen--Karpf--Thorup \cite{JesKarTho68}\index{Jessen, B.}\index{Karpf, J.}\index{Thorup, A.}, 
the authors gave an alternative proof, in which there is no need to distinguish 
between cases. 

\begin{lem}\label{L1.2.2}
 Let $A$ be an ordered commutative ring and $X$ be a module over $A$, and let 
 $A_{+}=\left\{a\in A\, \vert \, a>0\right\}$. 
 Let further $F, G\colon A^{2}\to X$ be functions which satisfy equations $(\alpha)$--$(\varepsilon)$ for all 
 $a, b, c\in A_{+}$. 
 Then these functions can be extended to functions $\widetilde{F}, \widetilde{G}\colon A^{2}\to X$ 
 which satisfy equations $(\alpha)$--$(\varepsilon)$ for all $a, b, c\in A$. 
\end{lem}
\begin{biz}
 Define the function $\widetilde{F}$ as in the proof of Lemma \ref{L1.2.1}. 
 Further let $\widetilde{G}(a, b)$ be zero when at least one of the elements 
 $a$, $b$, $a+b$ is zero. Otherwise we define $\widetilde{G}$ according to the following table, where 
 $+$, respectively $-$ stands for $>0$ and $<0$. 
 \begin{center}
 \begin{tabular}{|c|c|c|}
  \hline
 $a$ & $b$&  $\widetilde{G}(a, b)$ 
 \\ 
 \hline
 \hline
+ &+& $G(a, b)$ \\
 \hline
 + &-&$-G(a, -b)$ \\
 \hline
 - &+&$-G(-a, b)$ \\
 \hline
 - &-&$G(-a, -b)$ \\
 \hline
 \end{tabular} 
 \end{center}
 One easily verifies that the functions $\widetilde{F}, \widetilde{G}\colon A^{2}\to X$ thus defined satisfies the conditions. 
\end{biz}

From Theorem \ref{T1.2.2} with the choice $F=0$ (and interchanging the roles of $f$ and $g$) 
the following statement can be
obtained immediately.

\begin{theorem}\label{T1.2.3}
Let $A$ be a commutative ring, $X$ be a module over $A$ and
$f\colon A\rightarrow X$ be a function such that
\[
f(ab)=af(b)+bf(a)
\qquad
\left(a, b\in A\right). 
\]
Then the function $F\colon A\times A\rightarrow X$ defined by
$(A)$ fulfills
\begin{equation}
\tag{$\alpha$} F(a, b)=F(b, a), 
\end{equation}
equation $(\beta)$ is satisfied and also
\begin{equation}
\tag{$\eta$} F(ac, bc)=cF(a, b)
\end{equation}
holds for any $a, b, c\in A$.

Furthermore, in case $A$ is an integral domain and $X$ is a
unitary module over $A$ which is uniquely $A$--divisible, then
the function $F$ defined by the function $f$ through equation $(A)$ is the only function
which satisfies equations $(\alpha)$, $(\beta)$ and $(\eta)$.
\end{theorem}

\chapter{Characterization of derivations through one equation}

\section{Preliminaries and former results}

The purpose of this chapter is to provide characterization theorems on derivations.

As we saw in the second chapter, the characterization of derivations has an extensive literature, the reader should consult for
instance Horinouchi--Kannappan \cite{HorKan71}, Jurkat \cite{Jur65},
Kurepa \cite{Kur64, Kur65} and also the two monographs Kuczma \cite{Kuc09} and
Zariski--Samuel \cite{ZarSam75}.

Nevertheless, to the best of the author's knowledge, all of the characterizations have the following
form: additivity and another property imply that the function in question is a derivation.
We intend to show that derivations can be characterized by one single functional equation. 
The results of this chapter are based on the paper Gselmann \cite{Gse12}. 

More precisely, we would like to examine whether the equations occurring in the definition of
derivations are independent in the following sense.

Let $Q$ be a commutative ring and let $P$ be a subring of $Q$.
Let $\lambda, \mu\in Q\setminus\left\{0\right\}$  be arbitrary,
$f\colon P\rightarrow Q$ be a function and consider the
equation
\[
\lambda\left[f(x+y)-f(x)-f(y)\right]+
\mu\left[f(xy)-xf(y)-yf(x)\right]=0.
\quad
\left(x, y\in P\right)
\]
Clearly, if the function $f$ is a derivation, then this equation holds.
In the next section we will investigate the opposite direction, and it will be proved that, 
under some assumptions on the rings $P$ and $Q$, derivations can be
characterized through the above equation.
This result will be proved as a consequence of the main theorem that will be devoted
to the equation
\[
f(x+y)-f(x)-f(y)=g(xy)-xg(y)-yg(x),
\quad
\left(x, y\in P\right)
\]
where $f, g\colon P\rightarrow Q$ are unknown functions. 

We remark that similar investigations were made by
Dhombres \cite{Dho88}, Ger \cite{Ger98, Ger00} and also by Ger--Reich \cite{GerRei10} concerning ring homomorphisms. 
For instance, in Ger \cite{Ger98} the following theorem was proved. 
\begin{theorem}
 Let $X$ and $Y$ are two rings, and assume that 
for all $x\in X$ there exists $e_{x}\in X$ such that $x e_{x}=x$, suppose further 
that $Y$ has no elements of order 2 and does not admit zero divisors. 
If $f$ is a solution of 
\[
 f(x+y)+f(xy)=f(x)+f(y)+f(x)f(y) 
\qquad 
\left(x, y\in X\right)
\]
such that $f(0)=0$, 
then either $3f$ is even and $3f(2x)=0$ for all $x\in X$, 
or $f$ yields a homomorphism between $X$ and $Y$. 
\end{theorem}

During the proof of the main result the celebrated \emph{cocycle equation}
will play a key role.
About this equation one can read e.g. in
Acz\'{e}l \cite{Acz65}, Davison--Ebanks \cite{DavEba95},
Ebanks \cite{Eba82}, Erd\H{o}s \cite{Erd59},  
Hossz\'{u} \cite{Hos63} and also in Jessen--Karpf--Thorup \cite{JesKarTho68}. 
In the next section we will however utilize only
Theorem \ref{T1.2.2}.

\section{Characterization of derivations through one equation}

Our main result in this direction is contained in the following. 

\begin{theorem}\label{T3.2.1}
Let $\mathbb{F}$ be a field, $X$ be a vector space over $\mathbb{F}$ and
$f, g\colon \mathbb{F}\rightarrow X$ be functions such that
\begin{equation}
\tag{$\mathscr{E}$}
f(x+y)-f(x)-f(y)=g(xy)-xg(y)-yg(x)
\end{equation}
holds for all $x, y\in\mathbb{F}$.
Then, and only then, there exist additive functions $\alpha, \beta\colon\mathbb{F}\to X$  and a 
function 
$\varphi\colon \mathbb{F}\to X$ with the property
\[
 \varphi(xy)=x\varphi(y)+y\varphi(x) 
\qquad 
\left(x, y\in\mathbb{F}\right), 
\]
 such that 
\[
 f(x)=\beta(x)+\frac{1}{2}\alpha(x^{2})-x\alpha(x)
\qquad 
\left(x, y\in\mathbb{F}\right)
\]
and 
\[
 g(x)=\varphi(x)+\alpha(x) 
\qquad 
\left(x, y\in\mathbb{F}\right)
\]
are satisfied. 
\end{theorem}

\begin{biz}
Define the functions $\mathscr{C}_{f}, \mathscr{D}_{f}, \mathscr{C}_{g}$ and $\mathscr{D}_{g}$
on $\mathbb{F}\times \mathbb{F}$ by
\[
\begin{array}{lcl}
\mathscr{C}_{f}(x, y)&=&f(x+y)-f(x)-f(y) \\
\mathscr{D}_{f}(x, y)&=&f(xy)-xf(y)-yf(x) \\
\mathscr{C}_{g}(x, y)&=&g(x+y)-g(x)-g(y) \\
\mathscr{D}_{g}(x, y)&=&g(xy)-xg(y)-yg(x),
\end{array}
\]
respectively.
In view of Theorem \ref{T1.2.2}, we immediately get that the pairs
$(\mathscr{C}_{f}, \mathscr{D}_{f})$ and $(\mathscr{C}_{g}, \mathscr{D}_{g})$
fulfill the system of equations $(\alpha)$--$(\varepsilon)$.
Furthermore, equation $(\mathscr{E})$ yields that
\begin{equation}
\tag{$\mathscr{E}^{\ast}$}\mathscr{C}_{f}(x, y)=\mathscr{D}_{g}(x, y)
\end{equation}
for all $x, y\in \mathbb{F}$.
Due to equation $(\varepsilon)$,
\begin{equation}\label{Eq3.2.1}
\mathscr{C}_{g}(xz, yz)-z\mathscr{C}_{g}(x, y)=\mathscr{D}_{g}(x+y, z)-\mathscr{D}_{g}(x, z)-\mathscr{D}_{g}(y, z)
\end{equation}
holds for all $x, y, z\in\mathbb{F}$.
Interchanging the role of $x$ and $z$ in the previous equation, we obtain that
\begin{equation}\label{Eq3.2.2}
\mathscr{C}_{g}(xz, xy)-x\mathscr{C}_{g}(z, y)=\mathscr{D}_{g}(y+z, x)-\mathscr{D}_{g}(z, x)-\mathscr{D}_{g}(y, x)
\end{equation}
for any $x, y, z\in\mathbb{F}$.
Let us subtract equation \eqref{Eq3.2.2} from \eqref{Eq3.2.1}, to obtain
\begin{multline*}
\mathscr{C}_{g}(xz, yz)-z\mathscr{C}_{g}(x, y)
-\mathscr{C}_{g}(xz, xy)+x\mathscr{C}_{g}(z, y)
\\
=
\mathscr{D}_{g}(x+y, z)-\mathscr{D}_{g}(x, z)-\mathscr{D}_{g}(y, z)
-\mathscr{D}_{g}(z+y, x)+\mathscr{D}_{g}(z, x)+\mathscr{D}_{g}(y, x).
\end{multline*}
Because of $(\mathscr{E}^{\ast})$, the function $\mathscr{D}_{g}$ can be replaced by
$\mathscr{C}_{f}$. This implies however that
\begin{multline*}
\mathscr{C}_{g}(xz, yz)-z\mathscr{C}_{g}(x, y)
-\mathscr{C}_{g}(xz, xy)+x\mathscr{C}_{g}(z, y)
\\
=
\mathscr{C}_{f}(x+y, z)+\mathscr{C}_{f}(x, y)
-\mathscr{C}_{f}(x, y+z)-\mathscr{C}_{f}(y, z)=0,
\end{multline*}
for all $x, y, z\in\mathbb{F}$, where we used that the function
$\mathscr{C}_{f}$ fulfills $(\alpha)$ and $(\beta)$.
This equation with $z=1$ yields that
\[
\mathscr{C}_{g}(x, xy)=x\mathscr{C}_{g}(1, y),
\]
or if we replace $y$ by $\dfrac{y}{x}$, $(x\neq 0)$
\[
\mathscr{C}_{g}(x, y)=x\mathscr{C}_{g}\left(1, \frac{y}{x}\right).
\quad \left(x, y\in\mathbb{F}, x\neq 0\right)
\]
We will show that from this identity the homogeneity of
$\mathscr{C}_{g}$ follows. Indeed, let $t, x, y\in\mathbb{F},
t, x\neq 0$ be arbitrary, then
\[
\mathscr{C}_{g}(tx, ty)=tx\mathscr{C}_{g}\left(1, \frac{ty}{tx}\right)=
tx\mathscr{C}_{g}\left(1, \frac{y}{x}\right)=t\mathscr{C}_{g}(x, y).
\]
If $x=0$, we get from equation $\left(\mathscr{E}^{\ast}\right)$ that $\mathscr{C}_{g}(0,
0)=0$, thus for arbitrary $t\in\mathbb{F}$,
\[
\mathscr{C}_{g}(t 0, t 0)= 0=t\mathscr{C}_{g}(0, 0).
\]
Furthermore, in case $t=0$, then for any $x, y\in \mathbb{F}$
\[
\mathscr{C}_{g}(tx, ty)=\mathscr{C}_{g}(0, 0)=0=t\mathscr{C}_{g}(x, y).
\]
This means that the function $\mathscr{C}_{g}$ is homogeneous and
fulfills equations $(\alpha)$ and $(\beta)$. In view of Theorem
\ref{T1.2.3}, there exists a function $\varphi\colon\mathbb{F}\to X$ such that 
\[
 \varphi(xy)=x\varphi(y)+y\varphi(x)
\qquad 
\left(x, y\in\mathbb{F}\right)
\]
and 
\[
 \mathscr{C}_{g}(x, y)=\varphi(xy)-x\varphi(y)-y\varphi(x)
\qquad 
\left(x, y\in\mathbb{F}\right)
\]
hold. 
Due to the definition of the function $\mathscr{C}_{g}$, this yields that 
\[
 g(x)=\varphi(x)+\alpha(x) 
\qquad 
\left(x\in\mathbb{F}\right), 
\]
where the function $\varphi$ fulfills the above identity and $\alpha\colon\mathbb{F}\to X$ is additive. 
Writing this representation of the function $g$ into equation ($\mathscr{E}$), we have that 
\begin{equation}\label{Eq3.2.3}
 f(x+y)-f(x)-f(y)=\alpha(xy)-x\alpha(y)-y\alpha(x) 
\qquad 
\left(x, y\in\mathbb{F}\right). 
\end{equation}
Since the function $\alpha$ is additive, the two place function 
\[
 \mathscr{D}_{\alpha}(x, y)=\alpha(xy)-x\alpha(y)-y\alpha(x) 
\qquad 
\left(x, y\in\mathbb{F}\right)
\]
is a symmetric, biadditive function. 
Therefore, $\mathscr{D}_{\alpha}$ can be written as the Cauchy difference of its trace, that is 
\begin{multline*}
 \mathscr{D}_{\alpha}(x, y)=
\frac{1}{2}\alpha((x+y)^{2})-(x+y)\alpha(x+y)
\\
-\left(\frac{1}{2}\alpha(x^{2})-x\alpha(x)\right)
-\left(\frac{1}{2}\alpha(y^{2})-y\alpha(y)\right)
\qquad 
\left(x, y\in\mathbb{F}\right). 
\end{multline*}
In view of equation \eqref{Eq3.2.3}, this yields that the function 
\[
 x\mapsto f(x)-\left(\frac{1}{2}\alpha(x^{2})-x\alpha(x)\right) 
\qquad 
\left(x\in\mathbb{F}\right)
\]
is additive. 
Thus there exists an additive function $\beta\colon\mathbb{F}\to X$ such that 
\[
 f(x)=\beta(x)+\frac{1}{2}\alpha(x^{2})-x\alpha(x)
\qquad 
\left(x\in\mathbb{F}\right). 
\]
\end{biz}

According to Lemma \ref{L1.2.2}, with the aid of the previous
result, the following corollary can immediately be obtained. 

\begin{cor}
Let $\mathbb{F}$ be an ordered field, $X$ be a vector space over $\mathbb{F}$,
$\mathbb{F}_{+}=\left\{x\in\mathbb{F}\vert x>0\right\}$ and
$f, g\colon \mathbb{F}_{+}\rightarrow X$ be functions such that
\[
f(x+y)-f(x)-f(y)=g(xy)-xg(y)-yg(x)
\]
holds for all $x, y\in\mathbb{F}_{+}$.
Then the functions $f$ and $g$ can be extended to functions
$\widetilde{f}, \widetilde{g}\colon \mathbb{F}\rightarrow X$ such that
\[
\widetilde{f}(x)=\beta(x)+\frac{1}{2}\alpha(x^{2})-x\alpha(x)
\qquad
\left(x, y\in\mathbb{F}\right)
\]
and
\[
\widetilde{g}(x)=\varphi(x)+\alpha(x), 
\qquad
\left(x, y\in\mathbb{F}\right)
\]
where $\alpha, \beta\colon\mathbb{F}\to X$ are additive function and $\varphi\colon\mathbb{F}\to X$ fulfills 
\[
 \varphi(xy)=x\varphi(y)+y\varphi(x) 
\qquad 
\left(x, y\in\mathbb{F}\right). 
\]
\end{cor}

From our main result of this section, with the choice $g(x)=-\dfrac{\mu}{\lambda}f(x)$ the following corollary can be
derived easily.

\begin{cor}\label{C2.2}
Let $\mathbb{F}$ be a field and $X$ be a vector space over $\mathbb{F}$,
$\lambda, \mu\in \mathbb{F}\setminus\left\{0\right\}$. 
Then the function $f\colon \mathbb{F}\rightarrow X$ is a derivation if and only if
\[
\lambda\left[f(x+y)-f(x)-f(y)\right]+
\mu\left[f(xy)-xf(y)-yf(x)\right]=0
\]
holds for all $x, y\in \mathbb{F}$.
\end{cor}

Let us observe that Theorem \ref{T3.2.1} heavily uses that the functions involved map a field into a vector space. 
At the same time the problem can also be formulated in a more general setting, e.g., in (commutative) rings. 
Thus, we can ask the following. 

\begin{Opp}
 Let $A$ be a ring and $X$ be a module over $A$. Assume that 
 for the functions $f, g\colon A\to X$, functional equation 
 \[
  f(x+y)-f(x)-f(y)=g(xy)-xg(y)-yg(x)
 \]
is fulfilled for any $x, y\in A$. Prove or disprove that the same conclusion holds as in Theorem \ref{T3.2.1}. 
\end{Opp}

Our results can be restated with the aid of the notion of alien functional equations.

With the help of the above notions, Theorem \ref{T3.2.1}
says that the (additive) Cauchy equation, i.e.  
\[
 f(x+y)=f(x)+f(y) 
 \qquad 
 \left(x, y\in \mathbb{F}\right)
\]
and the Leibniz rule, that is, 
\[
 f(xy)=xf(y)+f(x)y 
 \qquad 
 \left(x, y\in \mathbb{F}\right). 
\]
are \emph{alien}, but not \emph{strongly alien}.

Let $f\colon \mathbb{R}\to \mathbb{R}$ be a derivation. Then $f$ clearly solves the system of functional equations 
\[
\begin{array}{rcl}
 f(x+y-xy)&=&f(x)+f(y)-f(xy) \\
 f(xy)&=&xf(y)+f(x)y
 \end{array}
 \qquad 
 \left(x, y\in \mathbb{R}\right)
\]
as well as the single equation 
\[
 f(x+y-xy)-f(x)-f(y)+f(xy)= f(xy)-xf(y)-f(x)y 
 \qquad 
 \left(x, y\in \mathbb{R}\right)
\]

Our second open problem is about the converse. 

\begin{Opp}
 Assume that the function $f\colon \mathbb{R}\to \mathbb{R}$ fulfills functional equation 
 \[
 f(x+y-xy)-f(x)-f(y)+f(xy)= f(xy)-xf(y)-f(x)y 
\]
for all $x, y\in \mathbb{R}$. 
\begin{enumerate}[(i)]
 \item Prove or disprove that the function $f$ is a real derivation. 
 \item What can be said if we consider the problem not on the set of the real numbers, but on rings or on fields?
\end{enumerate}
\end{Opp}

To formulate our last open problem in this section, let us observe that if 
$f\colon \mathbb{R}\to \mathbb{R}$ is a derivation, then $f$ yields a solution for the system of functional equations 
\[
\begin{array}{rcl}
 f\left(\dfrac{x+y}{2}\right)&=&f(x)+f(y) \\[2mm]
 f(xy)&=&xf(y)+f(x)y
 \end{array}
 \qquad 
 \left(x, y\in \mathbb{R}\right)
\]
as well as the single equation 
\[
 f\left(\frac{x+y}{2}\right)-f(x)-f(y)= f(xy)-xf(y)-f(x)y 
 \qquad 
 \left(x, y\in \mathbb{R}\right). 
\]

\begin{Opp}
 Assume that the function $f\colon \mathbb{R}\to \mathbb{R}$ fulfills functional equation 
 \[
 f\left(\frac{x+y}{2}\right)-f(x)-f(y)= f(xy)-xf(y)-f(x)y 
\]
for all $x, y\in \mathbb{R}$. 
\begin{enumerate}[(i)]
 \item Prove or disprove that the function $f$ is a real derivation. 
 \item What can be said if we consider the problem not on the set of the real numbers, but on rings or on fields?
\end{enumerate}
\end{Opp}

\chapter{Additive solvability and linear independence of the solutions of a system of functional equations}

\section{Introduction}

In this chapter we investigate two problems concerning derivations. 
On one hand, the additive solvability of the system of 
functional equations
\[
 d_{k}(xy)=\sum_{i=0}^{k}\Gamma(i,k-i) d_{i}(x)d_{k-i}(y)
  \qquad (x,y\in \mathbb{R},\,k\in\{0,\ldots,n\})
\]
is studied, where $\Delta_n:=\big\{(i,j)\in\mathbb{Z}\times\mathbb{Z}\vert 0\leq i,j\mbox{ and }i+j\leq n\big\}$
and $\Gamma\colon\Delta_n\to\mathbb{R}$ is a symmetric function such that $\Gamma(i,j)=1$ 
whenever $i\cdot j=0$. 

On the other hand, the linear dependence and independence of the additive 
solutions $d_{0},d_{1},\dots,d_{n}\colon \mathbb{R}\to\mathbb{R}$ of the above system of equations is characterized. 
As a consequence of the main result, for any nonzero real derivation $d\colon\mathbb{R}\to\mathbb{R}$, the 
iterates $d^0,d^1,\dots,d^n$ of $d$ are shown to be linearly independent, 
and the graph of the mapping $x\mapsto (x,d^1(x),\dots,d^n(x))$ to be dense in $\mathbb{R}^{n+1}$.

Given a real derivation $d\colon \mathbb{R}\to\mathbb{R}$, one can prove by induction that the iterates 
$d^0:=\mathop{id}$, $d^1:=d$, \dots, $d^{n}:=d\circ d^{n-1}$ of $d$ satisfy the following 
higher-order Leibniz rule
\[
 d^{k}(xy)=\sum_{i=0}^{k}\binom{k}{i}d^{i}(x)d^{k-i}(y)
  \qquad (x,y\in\mathbb{R},\,k\in\{1,\ldots,n\}).
\]

Motivated by this property, Heyneman--Sweedler \cite{HeySwe69}\index{Heyneman, R.~G.}\index{Sweedler, M.~E.} introduced the notion of $n$th-order 
derivation (in the context of functions mapping rings to modules, however, we will restrict 
ourselves only to real functions). 
\begin{defin}
Given $n\in \mathbb{N}$, a sequence of additive functions 
$d_{0},d_{1},\dots,d_{n}\colon \mathbb{R}\to\mathbb{R}$ is termed a \emph{derivation of order $n$}, if 
$d_0=\mathop{id}$ and, for any $k\in\{1,\ldots,n\}$,
\begin{equation}\label{dk}
 d_{k}(xy)=\sum_{i=0}^{k}\binom{k}{i}d_{i}(x)d_{k-i}(y) \qquad (x,y\in\mathbb{R})
\end{equation}
is fulfilled. 
\end{defin}

Clearly, a pair $(\mathop{id},d)$ is a first-order derivation if and only if 
$d$ is a derivation. More generally, if $d\colon\mathbb{R}\to\mathbb{R}$ is a derivation, 
then the sequence $(d^0,d^1,\dots,d^n)$ is a derivation of order $n$. However, if 
$\widetilde{d}\colon\mathbb{R}\to\mathbb{R}$ is a nontrivial derivation and $n\geq2$, then 
$(d^{0},d^{1},\dots,d^{n-1},d^{n}+\widetilde{d})$ is also an $n$th-order derivation 
where the last element is not the $n$th iterate of the derivation $d$. 

Firstly we will study the additive solvability of 
the following system of functional equations
\begin{equation}\label{dd}
 d_{k}(xy)=\sum_{i=0}^{k}\Gamma(i,k-i) d_{i}(x)d_{k-i}(y)
  \qquad (x,y\in \mathbb{R},\,k\in\{0,\ldots,n\}),
\end{equation}
where
\begin{equation}\label{Del}
  \Delta_n:=\big\{(i,j)\in\mathbb{Z}\times\mathbb{Z}\, \vert \, 0\leq i,j\mbox{ and }i+j\leq n\big\},
\end{equation}
and $\Gamma\colon\Delta_n\to\mathbb{R}$ is a symmetric function such that $\Gamma(i,j)=1$ 
whenever $i\cdot j=0$. 
After that we shall characterize the linear dependence and 
independence of the additive solutions $d_{0},d_{1},\dots,d_{n}\colon \mathbb{R}\to\mathbb{R}$ of \eqref{dd}.

\section{On the additive solvability of the system of functional equations \eqref{dd}}

As consequence of Theorem \ref{T1.2.1}, we can characterize those 
two-variable functions that are identical to the Leibniz difference of an additive function.

\begin{theorem}\label{LD}
Let $X$ be a real linear space and $D\colon\mathbb{R}^2\to X$. Then there 
exists an additive function $f\colon\mathbb{R}\to X$ fulfilling functional equation
\begin{equation}\label{dD}
  D(x,y)=f(xy)-xf(y)-yf(x) \qquad(x,y\in\mathbb{R})
\end{equation}
if and only if $D$ satisfies
\begin{equation}\label{D}
\begin{array}{rcl}
  D(x,y)&=&D(y,x)\qquad(x,y\in\mathbb{R}),\\[1.5mm]
  D(xy,z)+zD(x,y)&=&D(x,yz)+xD(y,z) \qquad(x,y,z\in\mathbb{R}),\\[1.5mm]
  D(x+y,z)&=&D(x,z)+D(y,z) \qquad(x,y,z\in\mathbb{R}).
  \end{array}
\end{equation}
\end{theorem}

Our first main result provides a sufficient condition on the recursive additive solvability of the 
functional equation \eqref{dd}. We deduce this result by using Theorem \ref{LD}, however, we note that another proof 
could be provided applying the results of Gselmann \cite{Gse12}.

\begin{theorem}\label{Solv} Let $n\geq2$ and $\Gamma\colon\Delta_n\to\mathbb{R}$ be a symmetric function such that 
$\Gamma(i,j)=1$ whenever $i\cdot j=0$ and
\begin{equation}\label{Gam}
  \Gamma(i+j,k)\Gamma(i,j)=\Gamma(i,j+k)\Gamma(j,k) \qquad (0\leq i,j,k\text{ and }i+j+k\leq n).
\end{equation}
Let $d_0=\mathop{id}$ and let $d_1,\dots,d_{n-1}\colon\mathbb{R}\to\mathbb{R}$ be additive functions such that \eqref{dd} holds
for $k\in\{1,\dots,n-1\}$. Then there exists an additive function $d_n\colon\mathbb{R}\to\mathbb{R}$ such that
\eqref{dd} is also valid for $k=n$.
\end{theorem}

\begin{biz} Using $\Gamma(0,n)=\Gamma(n,0)=1$, the functional equation for $d_n\colon\mathbb{R}\to\mathbb{R}$ can 
be rewritten as
\begin{equation}\label{dn}
  d_n(xy)-xd_n(y)-yd_n(x)=D_n(x,y)
:=\sum_{i=1}^{n-1}\Gamma(i,n-i) d_{i}(x)d_{k-i}(y) 
  \qquad(x,y\in\mathbb{R}).
\end{equation}
Thus, in view of Theorem \ref{LD}, in order that there exist an additive function $d_n$ such that \eqref{dn}
hold, it is necessary and sufficient that $D=D_n$ satisfy the conditions in \eqref{D}. 
The symmetry of $\Gamma$ implies the symmetry, the additivity of $d_1,\dots,d_{n-1}$ results the 
biadditivity of $D_n$. Thus, it suffices to prove that $D=D_n$ also satisfies the second identity
in \eqref{D}. 
This is equivalent to showing that, for all fixed $y\in\mathbb{R}$, the mapping 
$(x,z)\mapsto D_n(xy,z)+zD_n(x,y)$ is symmetric. Using equations \eqref{dd} for $k\in\{1,\dots,n-1\}$,
we obtain
\begin{multline}\label{aaa}
  D_n(xy,z)+zD_n(x,y) \\
  =\sum_{k=1}^{n-1}\Gamma(k,n-k) d_{k}(xy)d_{n-k}(z)
    +z\sum_{i=1}^{n-1}\Gamma(i,n-i) d_{i}(x)d_{n-i}(y) \\
  =\sum_{k=1}^{n-1}\Gamma(k,n-k) \bigg(\sum_{i=0}^{k}\Gamma(i,k-i)d_{i}(x)d_{k-i}(y)\bigg)d_{n-k}(z)\\
    +z\sum_{i=1}^{n-1}\Gamma(i,n-i) d_{i}(x)d_{n-i}(y) \\
  =\sum_{k=0}^{n}\sum_{i=0}^{k}\Gamma(k,n-k)\Gamma(i,k-i)d_{i}(x)d_{k-i}(y)d_{n-k}(z)\\
        -xyd_n(z)-xzd_n(y)-yzd_n(x)\\
  =\sum_{\alpha,\beta,\gamma\geq0,\,\alpha+\beta+\gamma=n}
       \Gamma(\alpha+\beta,\gamma)\Gamma(\alpha,\beta)d_{\alpha}(x)d_{\beta}(y)d_{\gamma}(z)\\
        -xyd_n(z)-xzd_n(y)-yzd_n(x).
\end{multline}
The sum of the last three terms in the above expression is symmetric in $(x,z)$. 
The symmetry of the first term is the consequence of the symmetry of 
$(\alpha,\gamma)\mapsto \Gamma(\alpha+\beta,\gamma)\Gamma(\alpha,\beta)$ which 
follows from property \eqref{Gam}.
\end{biz}

\begin{rem}\label{Rem1}
 Let $n\geq2$ and $\Gamma\colon\Delta_n\to\mathbb{R}$ be a symmetric function such that 
$\Gamma(i,j)=1$ whenever $i\cdot j=0$. Then condition \eqref{Gam} obviously holds if $i\cdot j\cdot k=0$ or 
if $i=k$, therefore \eqref{Gam} is equivalent to the following condition
\begin{equation}\label{Gam+}
  \Gamma(i+j,k)\Gamma(i,j)=\Gamma(i,j+k)\Gamma(j,k) \qquad (1\leq i,j,k,\,i\neq k\text{ and }i+j+k\leq n),
\end{equation}
which is non-trivial only for $n\geq4$. 

On the other hand, as we have seen it in the proof of Theorem \ref{Solv}, in order that there 
exist an additive function $d_n$ such that \eqref{dn}
hold, it is not only sufficient but also necessary that $D=D_n$ satisfy the second identity in \eqref{D}, which,
using formula \eqref{aaa}, the symmetry and the marginal values of $\Gamma$, is equivalent to the condition
\begin{equation}\label{Eq*}
 \sum_{\alpha,\beta,\gamma\geq1,\,\alpha\neq\gamma,\,\alpha+\beta+\gamma=n}
\Big(\Gamma(\alpha+\beta,\gamma)\Gamma(\alpha,\beta)-\Gamma(\gamma+\beta,\alpha)\Gamma(\gamma,\beta)\Big)
 d_{\alpha } (x)d_{\beta}(y)d_{\gamma}(z)
  =0.
\end{equation}
If $n\in\{0,1,2,3\}$, then the above condition holds automatically. For $n=4$, the above condition can be 
written as
\begin{equation}\label{pp}
 \Big(\Gamma(2,2)\Gamma(1,1)-\Gamma(1,3)\Gamma(1,2)\Big)
 d_{1}(y)\Big(d_{1} (x)d_{2}(z)-d_{1} (z)d_{2}(x)\Big)=0.
\end{equation}
Provided that $d_1$ and $d_2$ are linearly independent, \eqref{pp} implies that 
\[
\Gamma(2,2)\Gamma(1,1)=\Gamma(1,3)\Gamma(1,2),
\]
which proves that \eqref{Gam+} is necessary for $n=4$. An 
analogous and careful computation yields that in the case $n=5$, provided that $d_1$ and $d_2$ as well as 
$d_1$ and $d_3$ are linearly independent then $\Gamma(3,2)\Gamma(1,2)=\Gamma(1,4)\Gamma(2,2)$ and 
$\Gamma(2,3)\Gamma(1,1)=\Gamma(1,4)\Gamma(1,3)$ are necessary conditions which shows that \eqref{Gam+} is
necessary also for $n=5$. In the case $n\geq6$, we conjecture that linear independence of the functions 
$d_1,\dots,d_{n-2}$ is sufficient to guarantee the necessity of \eqref{Gam+}.
\end{rem}

In what follows, we describe the nowhere zero solutions of \eqref{Gam}.

\begin{theorem}\label{Coc}
Let $n\geq 2$ and $\Gamma\colon \Delta_{n}\to \mathbb{R}\setminus\{0\}$ 
be a symmetric function so that $\Gamma(i, j)=1$ whenever $i\cdot j=0$. 
Then $\Gamma$ satisfies the functional equation \eqref{Gam} if and only if there exists 
a function $\gamma\colon \{0, 1, \ldots, n\}\to \mathbb{R}\setminus \{0\}$ such that 
\begin{equation}\label{Gij}
 \Gamma(i, j)=\frac{\gamma(i+j)}{\gamma(i)\gamma(j)} 
\qquad 
((i, j)\in \Delta_{n}).
\end{equation}
\end{theorem}

\begin{biz}
Define the function $\gamma\colon \{0, 1, \ldots, n\}\to \mathbb{R}\setminus \{0\}$ through 
\[
 \gamma(k)= \prod_{\ell=1}^{k-1}\Gamma(\ell, 1) 
\qquad 
 (k\in\{0, 1, \ldots, n\}). 
\]
The empty product being equal to $1$, we have that $\gamma(0)=\gamma(1)=1$.

To complete the proof, we have to show that, for any $(i, j)\in \Delta_n$,
\[
 \Gamma(i, j)=\frac{\gamma(i+j)}{\gamma(i)\gamma(j)}. 
\]
This equivalent to proving that
\begin{equation}\label{ij}
 \Gamma(i,j)\prod_{\ell=1}^{i-1}\Gamma(\ell, 1) = \prod_{\ell=j}^{i+j-1}\Gamma(\ell, 1)
  \qquad((i, j)\in \Delta_n).
\end{equation}
This identity trivially holds for $i=0$, $i=1$ and for any $j\in\{0,\dots,n-i\}$. 
Let $j\in \{0,\dots,n-2\}$ be fixed. We prove \eqref{ij} by induction on 
$i\in\{1,\dots,n-j\}$. Assume that \eqref{ij} holds for $i\in\{1,\dots,n-j-1\}$. Then,
\begin{multline}\label{ij+}
  \Gamma(i+1,j)\prod_{\ell=1}^{i}\Gamma(\ell, 1)
  = \frac{\Gamma(i+1,j)\Gamma(i,1)}{\Gamma(i,j)}\bigg(\Gamma(i,j)\prod_{\ell=1}^{i-1}\Gamma(\ell, 1)\bigg)\\
  = \frac{\Gamma(i+1,j)\Gamma(i,1)}{\Gamma(i,j)}\prod_{\ell=j}^{i+j-1}\Gamma(\ell, 1)
  = \frac{\Gamma(i+1,j)\Gamma(i,1)}{\Gamma(i,j)\Gamma(i+j,1)}\prod_{\ell=j}^{i+j}\Gamma(\ell, 1).
\end{multline}
Using \eqref{Gam}, it follows that $\Gamma(i+1,j)\Gamma(i,1)=\Gamma(i,j)\Gamma(i+j,1)$, hence \eqref{ij+} yields
\eqref{ij} for $i+1$ instead of $i$.

Conversely, suppose that there exists a function 
$\gamma\colon \{0, 1, \ldots, n\}\to \mathbb{R}\setminus\{0\}$ such that 
\[
 \Gamma(i, j)=\frac{\gamma(i+j)}{\gamma(i)\gamma(j)} 
\qquad 
((i, j)\in \Delta_{n}). 
\]
Then, for any $i, j, k\geq 0$ with $i+j+k\leq n$, we have 
\[
\Gamma(i+j, k)\Gamma(i, j)
=\frac{\gamma(i+j+k)}{\gamma(i+j)\gamma(k)}\cdot \frac{\gamma(i+j)}{\gamma(i)\gamma(j)}
=\frac{\gamma(i+j+k)}{\gamma(i)\gamma(j+k)}\cdot \frac{\gamma(j+k)}{\gamma(j)\gamma(k)}
=\Gamma(i, j+k)\Gamma(j, k), 
\]
which completes the proof. 
\end{biz}

When $\Gamma$ is of the form \eqref{Gij}, then Theorem \ref{Solv} reduces to the following statement.

\begin{cor}\label{CorSolv}   Let $n\geq2$ and $\gamma\colon\{0, 1, \ldots, n\}\to \mathbb{R}\setminus \{0\}$ with $\gamma(0)=1$.
Let $d_0=\mathop{id}$ and let $d_1,\dots,d_{n-1}\colon\mathbb{R}\to\mathbb{R}$ be additive functions such that 
\begin{equation}\label{DK}
 d_{k}(xy)=\sum_{i=0}^{k}\frac{\gamma(k)}{\gamma(i)\gamma(k-i)}d_{i}(x)d_{k-i}(y) \qquad (x,y\in\mathbb{R})
\end{equation}
holds for $k\in\{1,\dots,n-1\}$. Then there exists an additive function $d_n\colon\mathbb{R}\to\mathbb{R}$ such that
\eqref{DK} is also valid for $k=n$.
\end{cor}

We note that if in the above corollary $\gamma(k)=k!$, then \eqref{DK} is equivalent to \eqref{dk}, that is
$\mathop{id}, d_1,\dots,d_{n}$ is a derivation of order $n$.

\section{A characterization of the linear dependence of additive functions}

\begin{theorem}Let $X$ be a Hausdorff locally convex linear space and let $a\colon\mathbb{R}\to X$ be an additive function. 
Then the following statements are equivalent:
\begin{enumerate}[(i)]
 \item there exists a nonzero continuous linear functional $\varphi\in X^*$ such that $\varphi\circ a=0$;
 \item there exists an upper semicontinuous function $\Phi\colon X\to\mathbb{R}$ such that $\Phi\not\geq0$ and $\Phi\circ 
a\geq0$;
 \item the range of $a$ is not dense in $X$, i.e.  $\overline{a(\mathbb{R})}\neq X$. 
\end{enumerate}
\end{theorem}

\begin{biz}
The implication \textit{(i)$\Rightarrow$(ii)} is obvious, because $\Phi$ can be chosen as $\varphi$.

To prove \textit{(ii)$\Rightarrow$(iii)}, assume that there exists an upper semicontinuous 
function $\Phi\colon X\to\mathbb{R}$ such that $\Phi\not\geq0$ and $\Phi\circ a\geq0$. Then 
$U\colon =\{x\in X\mid\Phi(x)<0\}$ is a nonempty and open set. The inequality $\Phi\circ a\geq0$ 
implies that $U\cap a(\mathbb{R})=\emptyset$, which proves that the range of $a$ cannot be dense in $X$. 

Finally, suppose that $\overline{a(\mathbb{R})}\neq X$. By the additivity of $a$, the set $a(\mathbb{R})$ is closed under
addition and multiplication by rational numbers. Therefore, the closure of $a(\mathbb{R})$ is a proper closed linear 
subspace of $X$. Then, by the Hahn--Banach theorem, there exists a nonzero continuous linear functional 
$\varphi\in X^*$ which vanishes on $a(\mathbb{R})$, i.e.  $\varphi\circ a=0$ is satisfied.
\end{biz}

By taking $X=\mathbb{R}^{n}$, the above theorem immediately simplifies to the following consequence which 
characterizes the linear dependence of finitely many additive functions.

\begin{cor}\label{cor1}
Let $n\in\mathbb{N}$ and $a_1,\dots,a_n\colon\mathbb{R}\to\mathbb{R}$ be additive functions. 
Then the following statements are equivalent:
\begin{enumerate}[(i)]
 \item the additive functions $a_1,\dots,a_n$ are linearly dependent, i.e.  there exist 
 $c_1,\dots,c_n\in \mathbb{R}$ such that $c_1^2+\cdots+c_n^2>0$ and $c_1a_1+\dots+c_na_n=0$;
 \item there exists an upper semicontinuous function $\Phi\colon\mathbb{R}^{n}\to\mathbb{R}$ such that $\Phi\not\geq0$ and 
 \[\Phi(a_1(x),\dots,a_n(x))\geq0\qquad(x\in\mathbb{R});\]
 \item the set $\{(a_1(x),\dots,a_n(x))\mid x\in\mathbb{R}\}$ is not dense in $\mathbb{R}^{n}$. 
\end{enumerate}
\end{cor}

In the particular case of this corollary, namely when $\Phi$ is an indefinite quadratic form, 
the equivalence of statements \textit{(i)} and \textit{(ii)} is the main result of the paper 
\cite{Koc12} by Kocsis. A former result in this direction is due to Maksa and Rätz \cite{MakRat81}: If two 
additive functions $a,b\colon\mathbb{R}\to\mathbb{R}$ satisfy $a(x)b(x)\geq0$ then $a$ and $b$ are linearly dependent.

\section{Linear independence of iterates of nonzero derivations}

In this section we apply Corollary \ref{cor1} to the particular case when the additive functions are iterates
of a real derivation. However, firstly we prove the following for higher order derivations. 

\begin{theorem}
Let $n\in\mathbb{N}$, let $\Gamma\colon\Delta_n\to\mathbb{R}$ be a symmetric function such that $\Gamma(i,j)=1$
whenever $i\cdot j=0$, \eqref{Gam} is satisfied and, for all $k\in\{2,\dots,n\}$ there exists 
$i\in\{1,\dots,k-1\}$
such that $\Gamma(i,k-i)\neq0$. Assume that $d_0=\mathop{id}$ and $d_{1},\dots,d_{n}\colon\mathbb{R}\to X$ are additive 
functions satisfying \eqref{dd} for all $k\in\{1,\dots,n\}$. Then the following statements are equivalent:
\begin{enumerate}[(i)]
 \item there exist $c_0,c_1,\dots,c_n\in \mathbb{R}$ such that $c_0^2+c_1^2+\cdots+c_n^2>0$ and
 \begin{equation}\label{ddd}
  c_0x+c_1d_{1}(x)+\cdots+c_nd_n(x)=0\qquad(x\in\mathbb{R});
 \end{equation}
 \item there exists an upper semicontinuous function $\Phi\colon\mathbb{R}^{n+1}\to\mathbb{R}$ such that $\Phi\not\geq0$ and 
\[
\Phi(x,d_{1}(x),\dots,d_{n}(x))\geq 0
\qquad(x\in\mathbb{R});
\]
 \item the set $\{(x,d_{1}(x),\dots,d_{n}(x))\mid x\in\mathbb{R}\}$ is not dense in $\mathbb{R}^{n+1}$;
 \item $d_{1}=0$.
\end{enumerate}
\end{theorem}

\begin{biz}
Applying Corollary \ref{cor1} to the additive functions $a_i(x)=d_{i}(x)$ $(i\in\{0,1,\dots,n\})$, it 
follows that \textit{(i)}, \textit{(ii)} and \textit{(iii)} are equivalent. The implication 
\textit{(iv)$\Rightarrow$(i)} is obvious since if $d_1=0$, then \textit{(i)} holds with
$c_1=1$ and $c_0=c_2=\cdots=c_n=0$.

Thus, it remains to show that \textit{(i)} implies \textit{(iv)}. Assume that 
\textit{(i)} holds. Then there exist a smallest $1\leq m\leq n$ and  
$c_0,\dots,c_m\in\mathbb{R}$ such that $c_0^2+c_1^2+\cdots+c_m^2>0$ and
\begin{equation}\label{dm}
   c_0x+c_1d_{1}(x)+\cdots+c_md_m(x)=0\qquad(x\in\mathbb{R}).
\end{equation}
This means that the equality
\[
  \gamma_0x+\gamma_1d_{1}(x)+\cdots+\gamma_{m-1}d_{m-1}(x)=0
  \qquad(x\in\mathbb{R})
\]
can only hold for $\gamma_0=\cdots=\gamma_{m-1}=0$.

Observe, that $d_1(1)=\cdots=d_n(1)=0$. Indeed, $d_1(1)=0$ is a consequence of \eqref{dd} when $k=1$ because
this equation means that $d_1$ is a derivation. The rest easily follows by induction on $k$ from \eqref{dd}.

Putting $x=1$ into \eqref{dm}, it follows that $c_0=0$. If $m=1$, then $c_1$ cannot be 
zero, hence we obtain that $d_{1}=0$.  
Thus, we may assume that the minimal $m$ for which \eqref{dm} is satisfied is non-smaller than $2$. 
Replacing $x$ by $xy$ in \eqref{dm} and applying \eqref{dd}, for all $x,y\in\mathbb{R}$, we get
\begin{multline*}
   0=\sum_{k=1}^m c_kd_k(xy)
    =\sum_{k=1}^m c_k\bigg(\sum_{i=0}^k\Gamma(i,k-i)d_{i}(x)d_{k-i}(y)\bigg) \\
   =\sum_{k=2}^m c_k\bigg(\sum_{i=1}^{k-1}\Gamma(i,k-i)d_{i}(x)d_{k\!-\!i}(y)\bigg)
      \!+\!x\bigg(\sum_{k=1}^m c_kd_k(y)\bigg)\!+\!y\bigg(\sum_{k=1}^m c_kd_k(x)\bigg)\\
   =\sum_{k=2}^m \sum_{i=1}^{k-1}c_k\Gamma(i,k-i)d_{i}(x)d_{k-i}(y)
    =\sum_{i=1}^{m-1} \sum_{k=i+1}^{m}c_k\Gamma(i,k-i)d_{i}(x)d_{k-i}(y) \\
   =\sum_{i=1}^{m-1}\bigg(\sum_{j=1}^{m-i}c_{i+j}\Gamma(i,j)d_j(y)\bigg)d_{i}(x).
\end{multline*}
By the minimality of $m$, it follows from the above equality that, for all $y\in\mathbb{R}$,
\[
   \sum_{j=1}^{m-i}c_{i+j}\Gamma(i,j)d_j(y)=0\qquad (i\in\{1,\dots,m-1\}).
\]
Again, by the minimality of $m$, this implies that $c_{i+j}\Gamma(i,j)=0$ for $(i,j)\in\Delta_m$
with $i,j\geq1$. By the assumption of the theorem, for all $k\in\{2,\dots,n\}$ there exists 
$i\in\{1,\dots,k-1\}$ such that $\Gamma(i,k-i)\neq0$. Thus, $c_2=\dots=c_m=0$. Therefore, by \eqref{dm}, 
$c_1$ cannot be equal to zero. Then \eqref{dm} simplifies to $d_{1}=0$, which was to be proved.
\end{biz}

Let $n\in \mathbb{N}$ be arbitrary and $d\colon \mathbb{R}\to \mathbb{R}$ be a derivation. Then the $(n+1)$-tuple 
$(\mathop{id}, d, d^{2},\dots, \allowbreak d^{n})$ is a derivation of order $n$. Thus from the 
previous theorem we immediately get the following. 

\begin{cor} Let $n\in\mathbb{N}$ and let $d\colon\mathbb{R}\to\mathbb{R}$ be a derivation. 
Then the following statements are equivalent:
\begin{enumerate}[(i)]
 \item there exist $c_0,c_1,\dots,c_n\in \mathbb{R}$ such that $c_0^2+c_1^2+\cdots+c_n^2>0$ and
 \[
 c_0x+c_1d(x)+\cdots+c_nd^n(x)=0
 \qquad(x\in\mathbb{R});
 \]
 \item there exists an upper semicontinuous function $\Phi\colon\mathbb{R}^{n+1}\to\mathbb{R}$ such that $\Phi\not\geq0$ and 
 \[
 \Phi(x,d(x),\dots,d^n(x))\geq0
 \qquad(x\in\mathbb{R});
 \]
 \item the set $\{(x,d(x),\dots,d^n(x))\mid x\in\mathbb{R}\}$ is not dense in $\mathbb{R}^{n+1}$;
 \item $d=0$.
\end{enumerate}
\end{cor}

\chapter{Characterization of derivations by actions on certain elementary function}

\section{Introduction and preparatory statements}

It is easy to see from Definition \ref{D2.1.1} that every derivation
$f\colon \mathbb{R}\rightarrow\mathbb{R}$ satisfies 
\begin{equation}\label{Eq1.3}
f(x^{k})=kx^{k-1}f(x)
\quad
\left(x\in\mathbb{R}\setminus\left\{0\right\}\right)
\end{equation}
for any fixed $k\in\mathbb{Z}\setminus\left\{0\right\}$.
Furthermore, the converse is also true, in the following sense:
if $k\in\mathbb{Z}\setminus\left\{0, 1\right\}$ is fixed and
an additive function $f\colon \mathbb{R}\rightarrow\mathbb{R}$ satisfies
\eqref{Eq1.3}, then $f$ is a derivation, see e.g.,
Jurkat \cite{Jur65}\index{Jurkat, W.}, Kurepa\index{Kurepa, S.} \cite{Kur64}, and
Kannappan--Kurepa \cite{KanKur70}\index{Kannappan, Pl.}.

Motivated by a problem of I.~Halperin\index{Halperin, I.} (1963), Jurkat \cite{Jur65} and,
independently, Kurepa \cite{Kur64} proved that every additive function
$f\colon \mathbb{R}\rightarrow\mathbb{R}$ satisfying
\[
f\left(\frac{1}{x}\right)=\frac{1}{x^{2}}f(x)
\quad
\left(x\in\mathbb{R}\setminus\left\{0\right\}\right)
\]
has to be linear.

In \cite{NisHor68} A.~Nishiyama\index{Nishiyama, A.} and S.~Horinouchi\index{Horinouchi, S.} investigated
additive functions $f\colon \mathbb{R}\rightarrow\mathbb{R}$
satisfying the additional equation
\begin{equation}\label{Eq1.4}
f(x^{n})=cx^{k}f(x^{m})
\quad
\left(x\in\mathbb{R}\setminus\left\{0\right\}\right),
\end{equation}
where $c\in\mathbb{R}$ and $n, m, k\in\mathbb{Z}$ are arbitrarily
fixed. This approach is obviously the common generalization of the
abovementioned results.
In the second part of this chapter we will deal with the stability
of this last system of functional equations.
Our main results are generalizations of the
theorems of \cite{NisHor68}.

We remark that in \cite{Bad06} R.~Badora\index{Badora, R.} solved
a stability problem for derivations mappings
between Banach algebras.
In this chapter we replace the Leibniz rule 
with an equation in a single variable, namely, with
a member of the family of equations in the form \eqref{Eq1.4}.
On the other hand, we restrict our considerations to real functions.

In order to avoid superfluous repetitions, henceforth we will say that the
function in question is \emph{locally regular} on its domain, if at least 
one of the following statements are fulfilled. 
\begin{enumerate}[(i)]
 \item bounded on a measurable set of positive measure;
\item continuous at a point;
\item there exists a set of positive Lebesgue measure so that the restriction of the function in question
is measurable in the sense of Lebesgue. 
\end{enumerate}

Furthermore, a function will be called \emph{globally regular} if instead of (ii), \\
(ii)' \; continuous on its domain \\
holds.

First we prove a simple lemma. 

\begin{lem}\label{L1}
 Let $\alpha\in\mathbb{R}$ and let us assume that for the 
function $\phi\colon ]0, +\infty[ \allowbreak \to\mathbb{R}$ the following statements are valid. 
\begin{enumerate}[(a)]
 \item the function $\phi$ is $\mathbb{Q}$-homogeneous of order $\alpha$, that is, 
\[
 \phi(rx)=r^{\alpha}\phi(x) 
\qquad 
\left(x\in ]0, +\infty[, r\in\mathbb{Q}\cap ]0, +\infty[\right). 
\]
\item the function $\phi$ is continuous at a point. 
\end{enumerate}
Then $\phi$ is continuous everywhere. 
\end{lem}
\begin{biz}
 Let us assume that the function $\phi$ is continuous at the point $x_{0}\in ]0, +\infty[$ and let 
$\tilde{x}\in ]0, +\infty[$ be arbitrary. 
Then, there exists a sequence of positive rational numbers $(r_{n})_{n\in\mathbb{N}}$ so that 
$\lim_{n\to\infty} r_{n}= \dfrac{\tilde{x}}{x_{0}}$. In this case 
the sequence $\left(\frac{1}{r_{n}}\right)_{n\in\mathbb{N}}$ is also a sequence of positive rational numbers and it 
converges to $\dfrac{x_{0}}{\tilde{x}}$. 
Due to property (a), 
\[
\frac{1}{r_{n}^{\alpha}}\phi(\tilde{x})=\phi\left(\frac{\tilde{x}}{r_{n}}\right) 
\qquad 
\left(n\in\mathbb{N}\right). 
\]
Taking the limit $n\to \infty$, the left hand side converges to 
$\left(\frac{x_{0}}{\tilde{x}}\right)^{\alpha}\phi(\tilde{x})$. Furthermore, 
the sequence $\left(\frac{\tilde{x}}{r_{n}}\right)_{n\in\mathbb{N}}$ converges to 
$x_{0}$, therefore the continuity of the function $\phi$ implies that 
the right hand side tends to $\phi(x_{0})$ as $n\to \infty$. 
This implies that 
\[
 \left(\frac{x_{0}}{\tilde{x}}\right)^{\alpha}\phi(\tilde{x})=\phi\left(\frac{x_{0}}{\tilde{x}}\tilde{x}\right)
\]
is fulfilled. 
Since $\tilde{x}\in ]0, +\infty[$ was arbitrary, we get that 
\[
 \phi(\lambda x)=\lambda^{\alpha}\phi(x) 
\qquad 
\left(\lambda, x\in ]0, +\infty[\right), 
\]
which obviously implies the (everywhere) continuity of the function $\phi$. 

\end{biz}

Furthermore, it is important pointing out the following fact. 
Fix $\alpha\in\mathbb{R}$ and let $\phi\colon ]0, +\infty[ \allowbreak \to\mathbb{R}$ 
be a $\mathbb{Q}$-homogeneous function of order $\alpha$. 
Suppose that $\phi$ fulfills property (i) or (iii) on the set of positive Lebesgue measure 
$D$. 
Then, the $\mathbb{Q}$-homogeneity of the function $\phi$ implies that 
(i), respectively (iii) holds for the function $\phi$ on the set 
$r\cdot D$ for all $r\in\mathbb{Q}$.

During the proof of our results we will also utilize a theorem a Kannappan--Kurepa \cite{KanKur70}. 

\begin{theorem}\label{KanKur70}
 Let $f, g\colon \mathbb{R}\rightarrow\mathbb{R}$ be additive functions and
$n, m\in\mathbb{Z}\setminus\left\{0\right\}$, $n\neq m$.
Suppose that
\[
f(x^{n})=x^{n-m}g(x^{m})
\]
holds for all $x\in\mathbb{R}\setminus\left\{0\right\}$.
Then the functions $F, G\colon \mathbb{R}\rightarrow\mathbb{R}$ defined by
\[
F(x)=f(x)-f(1)x
\quad
\text{and}
\quad
G(x)=g(x)-g(1)x
\quad
\left(x\in\mathbb{R}\right)
\]
are derivations and $nF(x)=mG(x)$ is fulfilled for all $x\in\mathbb{R}$.
\end{theorem}

\section{Derivations along monomial functions}

\begin{theorem}\label{Thm2.1}
 Let $n, m\in\mathbb{Z}\setminus\left\{0\right\}$, $n\neq m$ so that 
$n=-m$ or $\mathrm{sign}(n)=\mathrm{sign}(m)$, let further 
$f, g\colon\mathbb{R}\to\mathbb{R}$ be additive functions. 
Define the function $\phi\colon \mathbb{R}\setminus\left\{0\right\}\to\mathbb{R}$ by the 
formula 
\[
 \phi(x)=f\left(x^{n}\right)-x^{n-m}g\left(x^{m}\right) 
\qquad 
\left(x\in\mathbb{R}\setminus\left\{0\right\}\right), 
\]
and assume that $\phi$ is locally regular. 
Then, the functions $F, G\colon\mathbb{R}\to\mathbb{R}$ defined by 
\[
 F(x)=f(x)-f(1)x \quad 
\text{and}
\quad 
G(x)=g(x)-g(1)x
\qquad 
\left(x\in\mathbb{R}\right)
\]
are derivations and 
\[
 nF(x)=mG(x)
\]
holds for arbitrary $x\in\mathbb{R}$. 
\end{theorem}
\begin{biz}
 Concerning the values of $n$ and $m$ we have to distinguish three cases. 
At first, let us assume that $n, m>0$.
There is no loss of generality in assuming $n>m$. 
Define the function $\Phi$ on $\mathbb{R}^{n}$ by 
\[
 \Phi(x_{1}, \ldots, x_{n})=
f(x_{1}\cdots x_{n}) - 
\frac{1}{\binom{n}{m}} 
\sum_{\mathrm{card}(I)=m}
\left( \prod_{j \in \{1,2,\dots,n\} \setminus I} x_j \right)
g \left( \prod_{i \in I} x_i \right) 
\quad 
\left(x_{1}, \ldots, x_{n}\in\mathbb{R}\right), 
\]
where the summation is considered for all subsets $I$ of cardinality $m$ of the index set 
$\left\{1, \ldots, n\right\}$. 
Due to the additivity of the functions $f$ and $g$, the function 
$\Phi$ is a symmetric and $n$--additive function. Furthermore, its trace, that is, 
\[
 \Phi(x, \ldots, x)=\phi(x)=f\left(x^{n}\right)-x^{n-m}g\left(x^{m}\right) 
\qquad 
\left(x\in\mathbb{R}\right)
\]
is a polynomial function. On the other hand $\phi$ is a locally regular function.  
In view of Theorems \ref{Tszek1}, \ref{Tszek2}, \ref{Tszek3}, \ref{Tszek4}, this means that $\phi$ is a continuous polynomial function. 
Therefore, there exists $c\in\mathbb{R}$ such that 
\[
 \Phi(x_{1}, \ldots, x_{n})=cx_{1}\cdots x_{n} 
\quad
\left(x_{1}, \ldots, x_{n}\in\mathbb{R}\right), 
\]
therefore, 
\[
 \phi(x)=cx^{n} 
\quad 
\left(x\in\mathbb{R}\right). 
\]
With the substitution $x=1$, we get $\phi(1)=c$. 
On the other hand, the definition of the function $\phi$ yields that 
$\phi(1)=f(1)-g(1)$. 
Thus, 
\[
 f(x^{n})-x^{n-m}g(x^{m})=\left[f(1)-g(1)\right]x^{n} 
\qquad 
\left(x\in\mathbb{R}\right). 
\]
Define the functions $F, G\colon \mathbb{R}\to\mathbb{R}$ by 
\[
 F(x)=f(x)-f(1)x 
\quad 
\text{and}
\quad 
G(x)=g(x)-g(1)x 
\qquad 
\left(x\in\mathbb{R}\right). 
\]
Then the above identity yields that
\[
 F\left(x^{n}\right)=x^{n-m}G(x^{m}) 
\qquad 
\left(x\in\mathbb{R}\right). 
\]
The statement of the theorem follows now from Theorem \ref{KanKur70}. 

Secondly, let us assume that $n, m<0$. In this case we get that the function 
\[
 \phi(x)=f\left(x^{n}\right)-x^{n-m}g\left(x^{m}\right) 
\qquad 
\left(x\in\mathbb{R}\setminus\left\{0\right\}\right)
\]
is locally regular on its domain. Let $u\in\mathbb{R}\setminus\left\{0\right\}$, with the substitution $x=\frac{1}{u}$ this 
yields that 
\[
 \phi\left(\frac{1}{u}\right)=
f\left(u^{-n}\right)-u^{-(n-m)}g\left(u^{-m}\right) 
\qquad 
\left(u\in\mathbb{R}\setminus\left\{0\right\}\right). 
\]
Since $-n, -m>0$, the results of the previous case can be applied for the function 
\[
 \psi(u)=\phi\left(\frac{1}{u}\right) 
\qquad 
\left(u\in \mathbb{R}\setminus \left\{0\right\}\right), 
\]
which is, due to the local regularity of $\phi$, also locally regular. 

Finally, let us assume that $n=-m$. Without the loss of generality $m>0$ can be assumed. 
In this case 
\[
 \phi(x)=f\left(x^{-m}\right)-x^{-2m}g\left(x^{m}\right) 
\qquad 
\left(x\in\mathbb{R}\setminus\left\{0\right\}\right)
\]
is locally regular, or equivalently, the mapping 
\begin{equation}\label{Eq1}
\psi(x)= \phi\left(\sqrt[m]{x}\right)=f\left(\frac{1}{x}\right)-\frac{1}{x^{2}}g(x) 
\qquad 
\left(x>0\right)
\end{equation}
has the local regularity property. 
Let $u>0$ be arbitrary and  let us substitute $u(u+1)$ in place of $x$. Then 
\begin{multline*}
\psi(u(u+1))= \phi\left(\sqrt[m]{u(u+1)}\right)
\\
=f\left(\frac{1}{u(u+1)}\right)-\frac{1}{u^{2}(u+1)^{2}}g\left(u(u+1)\right)
\\
\left(u>0\right). 
\end{multline*}
Using the additivity of the function $f$, 
\begin{multline*}
\psi(u(u+1))=\phi\left(\sqrt[m]{u(u+1)}\right)
\\
=f\left(\frac{1}{u}\right)-f\left(\frac{1}{u+1}\right) -\frac{1}{u^{2}(u+1)^{2}}g\left(u(u+1)\right)
\\
\left(u>0\right). 
\end{multline*}
On the other hand, 
\[
\psi(u)=\phi\left(\sqrt[m]{u}\right)=f\left(\frac{1}{u}\right)-\frac{1}{u^{2}}g(u) 
\qquad 
\left(u>0\right)
\]
and 
\[
\psi(u+1)=\phi\left(\sqrt[m]{u+1}\right)=f\left(\frac{1}{u+1}\right)-\frac{1}{(u+1)^{2}}g(u+1) 
\qquad 
\left(u>0\right). 
\]
Therefore, 
\begin{multline*}
\psi(u(u+1))-\psi(u)+\psi(u+1)
\\
=
\phi\left(\sqrt[m]{u(u+1)}\right)-\phi(\sqrt[m]{u})+\phi(\sqrt[m]{u+1})
\\
=
f\left(\frac{1}{u}\right)-f\left(\frac{1}{u+1}\right) -\frac{1}{u^{2}(u+1)^{2}}g\left(u(u+1)\right)
\\
-f\left(\frac{1}{u}\right)+\frac{1}{u^{2}}g(u)+
f\left(\frac{1}{u+1}\right)-\frac{1}{(u+1)^{2}}g(u+1) 
\qquad 
\left(u>0\right)
\end{multline*}
Making use of the additivity of the function g, after rearrangement, we obtain that 
\[
 \chi(u)=2ug(u)-g\left(u^{2}\right) 
\qquad 
\left(u>0\right), 
\]
where 
\[
 \chi(u)=u^{2}(u+1)^{2}\left[\psi\left(u(u+1)\right)-\psi(u)+\psi(u+1)\right]+u^{2}g(1) 
\qquad 
(u>0). 
\]
By our assumptions, the function $\phi$  is locally regular on 
$\mathbb{R}\setminus\left\{0\right\}$ and due to the additivity of $f$ and $g$, 
it is $\mathbb{Q}$-homogeneous of order $n$. Thus, by Lemma \ref{L1}, 
$\phi$ is globally regular on $\mathbb{R}\setminus \left\{0\right\}$. 
This implies that $\psi$ is globally regular on $]0, +\infty[$, which means that 
$\chi$ is locally regular. 
Due to the results of the first case this yields that the function 
$G\colon\mathbb{R}\to\mathbb{R}$ defined by 
\[
 G(x)=g(x)-g(1)x 
\qquad 
\left(x\in\mathbb{R}\right)
\]
is a derivation. 
In view of \eqref{Eq1}, this implies that 
\[
\phi\left(\sqrt[m]{x}\right)=f\left(\frac{1}{x}\right)-\frac{1}{x^{2}}\left[G(x)+g(1)x\right] 
\qquad 
\left(x>0\right), 
\]
that is , 
\[
\phi\left(\sqrt[m]{x}\right)=f\left(\frac{1}{x}\right)+G\left(\frac{1}{x}\right)+g(1)\frac{1}{x} 
\qquad 
\left(x>0\right), 
\]
since $G$ is a derivation. 
Let $u>0$, with the substitution $x=\frac{1}{u}$ we get that 
\[
\psi(u)=\phi\left(\sqrt[m]{\frac{1}{u}}\right)=f(u)+G(u)+g(1)u
\qquad 
\left(u>0\right). 
\]
Let us observe that the right hand side of this identity is an additive function, being the sum of additive functions. 
Moreover, the left hand side is locally regular, due to the local regularity of $ \phi$. 
Thus $\psi$ is a regular additive function, which means that there exists $c\in\mathbb{R}$ so that 
\[
 f(u)+G(u)+g(1)u=cu 
\qquad 
\left(u\in\mathbb{R}\right). 
\]
With $u=1$, $c=f(1)+g(1)$ can be obtained, therefore, 
\[
 f(u)=\left[f(1)+g(1)\right]u-g(1)u+G(u) 
\qquad 
\left(u\in\mathbb{R}\right), 
\]
i.e.  
\[
 f(u)=-G(u)+f(1)u
\qquad 
\left(u\in\mathbb{R}\right). 
\]
This means that the function $F\colon \mathbb{R}\to\mathbb{R}$ defined by 
\[
 F(x)=f(x)-f(1)x 
\qquad 
\left(x\in\mathbb{R}\right)
\]
is a derivation and 
\[
 F(x)=-G(x) 
\qquad 
\left(x\in\mathbb{R}\right)
\]
holds. 
\end{biz}

\begin{lem}\label{L2.1}
 Let $\kappa\in\mathbb{R}$, $n, m\in\mathbb{Z}$, $n \neq m$ and assume that 
$f\colon\mathbb{R}\to\mathbb{R}$ is an additive function. 
Define the function $\phi\colon\mathbb{R}\setminus\left\{0\right\}\to\mathbb{R}$ by 
\[
 \phi(x)=f\left(x^{n}\right)-\kappa x^{n-m}f\left(x^{m}\right) 
\qquad 
\left(x\in\mathbb{R}\setminus \left\{0\right\}\right)
\]
and assume that $\phi$ is locally regular. 
Then, the function $F\colon\mathbb{R}\to\mathbb{R}$ defined by 
\[
 F(x)=f(x)-f(1)x \quad 
\left(x\in\mathbb{R}\setminus\left\{0\right\}\right)
\]
is a derivation so that for any $x\in\mathbb{R}$ 
\[
\left(n-\kappa m\right)F(x)=0. 
\]
\end{lem}
\begin{biz}
 In view of the previous theorem, it is enough to deal with the case $\mathrm{sign}(n)\neq\mathrm{sign}(m)$ and $n\neq -m$. 
Due to the definition of the function $\phi$ 
\[
 \phi\left(x^{n}\right)=f\left(x^{n^{2}}\right)-\kappa x^{n(n-m)}f\left(x^{nm}\right)  
\qquad 
\left(x\in\mathbb{R}\setminus\left\{0\right\}\right)
\]
and 
\[
 \phi\left(x^{m}\right)=f\left(x^{nm}\right)-\kappa x^{m(n-m)}f\left(x^{m^{2}}\right) 
\qquad 
\left(x\in\mathbb{R}\setminus\left\{0\right\}\right), 
\]
therefore 
\[
 \phi\left(x^{n}\right)+\kappa x^{n(n-m)}\phi\left(x^{m}\right)=
f\left(x^{n^{2}}\right)-\kappa^{2} x^{n^{2}-m^{2}}f\left(x^{m^{2}}\right) 
\qquad 
\left(x\in\mathbb{R}\setminus\left\{0\right\}\right). 
\]
By our assumptions, $\phi$ is a locally regular function on 
$\mathbb{R}\setminus\left\{0\right\}$. 
However, the additivity of $f$ implies that 
\[
 \phi(rx)=r^{n}\phi(x) 
\qquad 
\left(x\in\mathbb{R}\setminus\left\{0\right\}, r\in\mathbb{Q}\setminus \left\{0\right\}\right). 
\]
Using Lemma \ref{L1}, we get that $\phi$ is globally regular. 
Therefore, the function
\[
 \psi(x)=\phi\left(x^{n}\right)+\kappa x^{n(n-m)}\phi\left(x^{m}\right) 
\qquad 
\left(x\in\mathbb{R}\setminus\left\{0\right\}\right)
\]
is locally regular. 
Since $n^{2}, m^{2}>0$ and $n^{2}\neq m^{2}$, the results of the previous theorem can be applied 
(with the choice $\psi(x)=\phi\left(x^{n}\right)+\kappa x^{n(n-m)}\phi\left(x^{m}\right)$ and $g(x)=\kappa^{2}f(x)$) 
to obtain that 
\[
 f(x)=F(x)+f(1)x 
\qquad 
\left(x\in\mathbb{R}\right), 
\]
where $F\colon \mathbb{R}\to\mathbb{R}$ is a derivation and 
\[
 nF(x)=m\kappa F(x)
\]
is also fulfilled for all $x\in\mathbb{R}$.

\end{biz}

From this lemma, the following corollary can be concluded immediately. 

\begin{cor}\label{C2.1}
 Let $r\in\mathbb{Q}\setminus\left\{0, 1\right\}$ and 
$f\colon\mathbb{R}\to\mathbb{R}$ be an additive function and define the function
 by
\[
 \phi(x)=f\left(x^{r}\right)-rx^{r-1}f\left(x\right)
\qquad 
\left(x\in\mathbb{R},\, x>0\right), 
\]
and assume that $\phi$ is locally regular. 
Then, the function $F\colon\mathbb{R}\to\mathbb{R}$ defined by 
\[
 F(x)=f(x)-f(1)x \quad 
\left(x\in\mathbb{R}\right)
\]
is a derivation. 
\end{cor}

\section{Derivations along rational functions}

In view of the results of the previous subsection, we are able to prove the following. 
The results presented here can be considered as a generalization that of Halter-Koch--Reich 
\cite{HalRei99, HalRei00, HalRei01}. 

\begin{theorem}\label{T2.2}
 Let $n\in\mathbb{Z}\setminus\left\{0\right\}$ and 
$\left(\begin{array}{cc}
a&b\\
c&d
\end{array}
\right)\in\mathbf{GL}_{2}(\mathbb{Q})$ be such that 
\begin{enumerate}[--]
 \item if $c=0$, then $n\neq 1$;
\item if $d=0$, then $n\neq -1$. 
\end{enumerate}
Let further $f, g\colon\mathbb{R}\to\mathbb{R}$ be additive functions and define the function 
$\phi$ by 
\[
 \phi(x)=f\left(\frac{ax^{n}+b}{cx^{n}+d}\right)-\frac{x^{n-1}g(x)}{\left(cx^{n}+d\right)^{2}} 
\qquad 
\left(x\in\mathbb{R}, \,
cx^{n}+d\neq 0\right). 
\]
Let us assume $\phi$ to be globally regular. 
Then, the functions $F, G\colon\mathbb{R}\to\mathbb{R}$ defined by 
\[
 F(x)=f(x)-f(1)x \quad \text{and} 
\quad 
G(x)=g(x)-g(1)x
\quad
\left(x\in\mathbb{R}\right)
\]
are derivations. 
\end{theorem}

\begin{biz}
 Firstly, let us suppose that  $c=0$. This means that the function 
\[
 \phi(x)=f\left(\frac{a}{d}x^{n}+\frac{b}{d}\right)-\frac{1}{d^{2}}x^{n-1}g(x)
\qquad 
\left(x\in\mathbb{R}\setminus\left\{0\right\}\right)
\]
is globally regular.  
In this case, the statement immediately follows from Theorem \ref{Thm2.1}. 

Similarly, if $d=0$, then 
\[
 \phi(x)=f\left(\frac{a}{c}+\frac{b}{c}x^{-n}\right)-x^{-n-1}g(x) 
\qquad 
\left(x\in\mathbb{R}\setminus\left\{0\right\}\right)
\]
is globally regular. Therefore, due to Theorem \ref{Thm2.1}, we obtain that the functions 
\[
 F(x)=f(x)-f(1)x \quad 
\text{and}
\quad 
G(x)=g(x)-g(1)x 
\qquad 
\left(x\in\mathbb{R}\right)
\]
are derivations. 

Thus, henceforth $cd\neq 0$ can be assumed. Furthermore, due to the $\mathbb{Q}$-homogeneity of the functions 
$f$ and $g$, $c=1$ can be supposed. That is, 
\[
 \phi(x)=f\left(\frac{ax^{n}+b}{x^{n}+d}\right)-\frac{x^{n-1}g(x)}{\left(x^{n}+d\right)^{2}} 
\qquad 
\left(x\in\mathbb{R}, \,
x^{n}+d\neq 0\right). 
\]
Since the function $f$ is additive, 
\[
 f\left(\frac{ax^{n}+b}{x^{n}+d}\right)=f(a)-f\left(\frac{D}{x^{n}+d}\right), 
\qquad 
\left(x\in\mathbb{R}, \,
x^{n}+d\neq 0\right), 
\]
therefore, 
\begin{equation}\label{Eq2}
 \phi(x)=f(a)-f\left(\frac{D}{x^{n}+d}\right)-\frac{x^{n-1}g(x)}{(x^{n}+d)^{2}}
\qquad 
\left(x\in\mathbb{R}, \,
x^{n}+d\neq 0\right), 
\end{equation}
where $D=\det\left(\begin{array}{cc}
                    a&b\\
1&d
                   \end{array}
\right)$. 
Let us observe that 
\[
 \frac{D}{x^{n}+d}=\frac{D}{d}-\frac{D}{\left(\sqrt[n]{d^{2}}\frac{1}{x}\right)^{n}+d}
\]
holds for all $x\in\mathbb{R}, \, x\neq 0, x^{n}+d\neq 0.$
Using this identity, we get 
\begin{multline}\label{Eq3}
\phi(x)=
f(a)-f\left(\frac{D}{d}\right)+f\left(\frac{D}{\left(\sqrt[n]{d^{2}}\frac{1}{x}\right)^{n}+d}\right)
-\frac{x^{n-1}g(x)}{\left(x^{n}+d\right)^{2}}
\\
\left(x\in\mathbb{R}, \, x\neq 0, x^{n}+d\neq 0\right), 
\end{multline}
where the additivity of the function $f$ was also used. 
Let us replace $x$ by $\sqrt[n]{d^{2}}\dfrac{1}{x}$ in \eqref{Eq2} to acquire 
\begin{multline*}
 \phi\left(\sqrt[n]{d^{2}}\frac{1}{x}\right)=
f(a)-f\left(\frac{D}{\left(\sqrt[n]{d^{2}}\frac{1}{x}\right)^{n}+d}\right)-
\frac{\left(\sqrt[n]{d^{2}}\frac{1}{x}\right)^{n-1}g\left(\sqrt[n]{d^{2}}\frac{1}{x}\right)}
{\left(\left(\sqrt[n]{d^{2}}\frac{1}{x}\right)^{n}+d\right)^{2}}
\\
\left(x\in\mathbb{R}, \, x\neq 0, x^{n}+d\neq 0\right). 
\end{multline*}
Since 
\[
 \left(\sqrt[n]{d^{2}}\frac{1}{x}\right)^{n}+d= d\frac{1}{x^{n}}\left(x^{n}+d\right), 
\]
the above identity yields that 
\begin{multline*}
\phi\left(\sqrt[n]{d^{2}}\frac{1}{x}\right)=
f(a)-f\left(\frac{D}{\left(\sqrt[n]{d^{2}}\frac{1}{x}\right)^{n}+d}\right)-
\frac{\left(\sqrt[n]{d^{2}}\frac{1}{x}\right)^{n-1}g\left(\sqrt[n]{d^{2}}\frac{1}{x}\right)}
{\left(d\frac{1}{x^{n}}\right)^{2}\left(\left(x^{n}+d\right)\right)^{2}} 
\\
\left(x\in\mathbb{R}, \, x\neq 0, x^{n}+d\neq 0\right).
\end{multline*}
After some rearrangement, we arrive at 
\begin{multline}\label{Eq4}
 \phi\left(\sqrt[n]{d^{2}}\frac{1}{x}\right)=
f(a)-f\left(\frac{D}{\left(\sqrt[n]{d^{2}}\frac{1}{x}\right)^{n}+d}\right)-
\frac{x^{n-1}\frac{1}{\sqrt[n]{d^{2}}}x^{2}g\left(\sqrt[n]{d^{2}}\frac{1}{x}\right)}{\left(x^{n}+d\right)^{2}}
\\
\left(x\in\mathbb{R}, \, x\neq 0, x^{n}+d\neq 0\right).
\end{multline}
In case we add \eqref{Eq3} and \eqref{Eq4} together, 
\begin{multline*}
 \phi(x)+\phi\left(\sqrt[n]{d^{2}}\frac{1}{x}\right)=
f\left(2a-\frac{D}{d}\right)-
\frac{x^{n-1}}{\left(x^{n}+d\right)^{2}} 
\left[
g(x)+x^{2}\frac{1}{\sqrt[n]{d^{2}}}g\left(\sqrt[n]{d^{2}}\frac{1}{x}\right)
\right]
\\
\left(x\in\mathbb{R}, \, x\neq 0, x^{n}+d\neq 0\right).
\end{multline*}
Let us define the functions
\[
 h(x)=\frac{1}{\sqrt[n]{d^{2}}}g\left(\sqrt[n]{d^{2}}x\right)
\qquad 
\left(x\in\mathbb{R}\right)
\]
and 
\begin{multline*}
 \psi(x)=
-\frac{1}{x^{2}}\frac{\left(x^{n}+d\right)^{2}}{x^{n-1}}
\left[
\phi(x)+\phi\left(\sqrt[n]{d^{2}}\frac{1}{x}\right)-f\left(2a-\frac{D}{d}\right)
\right]
\\
\left(x\in\mathbb{R}, \, x\neq 0, x^{n}+d\neq 0\right).
\end{multline*}
In this case
\[
 \psi(x)=h\left(\frac{1}{x}\right)+\frac{1}{x^{2}}g(x)
\quad 
\left(x\in\mathbb{R}, \, x\neq 0, x^{n}+d\neq 0\right)
\]
holds. 
By our assumptions $\phi$ is a globally regular mapping, therefore the function 
$\psi$ has the local regularity property. 
Due to Theorem \ref{Thm2.1}, this gives that the functions 
$F, G\colon\mathbb{R}\to\mathbb{R}$ defined by 
\[
 F(x)=f(x)-f(1)x
\quad 
\text{and}
\quad 
G(x)=g(x)-g(1)x 
\quad 
\left(x\in\mathbb{R}\right)
\]
are derivations.
\end{biz}

\section{A characterization of linearity}

Finally, in the last part of this chapter we present a characterization of linearity. 
Just as in the proof of Theorem \ref{T2.2}, Theorem \ref{Thm2.1} plays again an important role. 

\begin{theorem}\label{lin}
 Let $n\in\mathbb{N}, \, n\neq 1$ and $f\colon\mathbb{R}\to\mathbb{R}$ be an additive function. 
Define $\phi$ on $\mathbb{R}$ by
\[
 \phi(x)=f\left(x^{n}\right)-f(x)^{n} 
\qquad 
\left(x\in\mathbb{R}\right). 
\]
Let us assume that $\phi$ is locally regular. 
Then the function $f$ is linear, that is, 
\[
 f(x)=f(1)x
\]
holds for all $x\in\mathbb{R}$. 
\end{theorem}
\begin{biz}
 Let us define the function $\Phi\colon\mathbb{R}^{n}\to\mathbb{R}$ by 
\[
 \Phi(x_{1}, \ldots, x_{n})=f\left(x_{1}\cdots x_{n}\right)-f(x_{1})\cdots f(x_{n}) 
\qquad 
\left(x_{1}, \ldots, x_{n}\in\mathbb{R}\right). 
\]
Due to the additivity of $f$, the function $\Phi$ is a symmetric, $n$-additive function. Furthermore, 
\[
 \Phi(x, \ldots, x)=\phi(x)=f\left(x^{n}\right)-f(x)^{n}
\qquad 
\left(x\in\mathbb{R}\right). 
\]
From the local regularity of the function $\phi$  we immediately deduce that 
$\phi$ is a continuous polynomial function. 
Consequently, 
\begin{equation}\label{Eq5}
 \Phi(x_{1}, \ldots, x_{n})=cx_{1}\cdots x_{n} 
\qquad 
\left(x_{1}, \ldots, x_{n}\in\mathbb{R}\right)
\end{equation}
holds with a certain $c\in\mathbb{R}$. Due to the definition of the 
function $\phi$, we have $\phi(1)=f(1)-f(1)^{n}$. On the other hand 
\[
 \phi(1)=\Phi(1, \ldots, 1)=c. 
\]
Hence $c=f(1)-f(1)^{n}$. 
Let $u\in\mathbb{R}$, with the substitution 
\[
 x_{1}=u, \qquad x_{i}=1 \, \text{for $i\geq 2$}, 
\]
equation \eqref{Eq5} yields that 
\[
 f(u)-f(u)f(1)^{n-1}=\left(f(1)-f(1)^{n}\right)u  
\qquad 
\left(u\in\mathbb{R}\right). 
\]
In case $f(1)^{n-1}\neq 1$, this furnishes 
\[
 f(u)=f(1)u 
\qquad 
\left(u\in\mathbb{R}\right). 
\]
If $f(1)^{n-1}=1$, then 
\[
 c=f(1)-f(1)^{n}=f(1)\left[1-f(1)^{n-1}\right]=0. 
\]
Therefore equation \eqref{Eq5} with the substitutions 
\[
 x_{1}=u, \quad x_{2}=v, \quad \text{and} \quad x_{i}=1 \quad \text{for $i\geq 3$} 
\qquad 
\left(u, v\in\mathbb{R}\right)
\]
yields that 
\[
 f(uv)=f(u)f(v)f(1)^{n-2} 
\qquad 
\left(u, v\in\mathbb{R}\right), 
\]
that is, $f(1)^{n-2}\cdot f$ is a non identically zero real homomorphism. 
In view of Kuczma \cite[Theorem 14.4.1.]{Kuc09} this implies that 
\[
 f(1)^{n-2}f(u)=u 
\qquad 
\left(u\in\mathbb{R}\right). 
\]
Since $f(1)^{n-2}\cdot f(1)=f(1)^{n-1}=1$, 
\[
 \frac{f(u)}{f(1)}=u
\]
holds for all $u\in\mathbb{R}$, that is, $f$ is a linear function, indeed. 
\end{biz}

\begin{defin}
 Let $R, R'$ be rings, $n\in\mathbb{N}, n\geq 2$ be fixed. 
The function 
$\varphi\colon R\rightarrow R'$ is called an \emph{$n$-homomorphism}\index{homomorphism!--- $n$} if 
\[
\varphi(a+b)=\varphi(a)+\varphi(b) 
\qquad
\left(a, b\in R\right) 
\]
and 
\[
 \varphi(a_{1}\cdots a_{n})=\varphi(a_{1})\cdots \varphi(a_{n})
\qquad 
\left(a_{1}, \ldots, a_{n}\in R\right). 
\]

The function $\varphi\colon R\rightarrow R'$ is called an \emph{$n$-Jordan homomorphism}\index{homomorphism!--- $n$-Jordan} if 
\[
\varphi(a+b)=\varphi(a)+\varphi(b) 
\qquad
\left(a, b\in R\right) 
\]
and 
\[
 \varphi(a^{n})=\varphi(a)^{n}
\qquad 
\left(a\in R\right). 
\]
\end{defin}

It was G.~Ancochea who firstly dealt with the connection of Jordan homomorphisms and 
homomorphisms, see \cite{Anc42}. The results of G.~Ancochea\index{Ancochea, G.} were generalized and extended in several ways, 
see for instance Jacobson--Rickart \cite{JacRic50}\index{Jacobson, N.}\index{Rickart, C.~E.}, 
Kaplansky \cite{Kap47}\index{Kaplansky, I.}, {\.Z}elazko \cite{Zel68}\index{{\.Z}elazko, W.}.

In Gselmann \cite{Gse14b} we proved a generalization of of the above result, namely the following theorems.

\begin{theorem}\label{Tncomm}
 Let $n\in\mathbb{N}, n\geq 2$ $R$ be a ring, $R'$ be a locally convex algebra 
over the field $\mathbb{F}$ of characteristic zero,  
$\varphi\colon R\rightarrow R'$ be an additive function and assume that the mapping 
\[
R \ni x\longmapsto \varphi(x^{n})-\varphi(x)^{n}
\]
is bounded on $R$. Then the function $\varphi$ is an $n$-Jordan homomorphism. 
\end{theorem}

\begin{theorem}
Let $n\in\mathbb{N}, n\geq 2$, $\mathbb{F}$ be a field of characteristic zero,
$R$ be a commutative topological ring and $R'$ be a commutative topological algebra over the 
field $\mathbb{F}$. Furthermore, let us consider the additive function 
$\varphi\colon R\rightarrow R'$ and suppose that for the map $\phi$ defined on $R$ by 
\[
 \phi(x)=\varphi(x^{n})-\varphi(x)^{n}
\qquad
\left(x\in R\right)
\]  
one of the following statements hold.  
\begin{enumerate}[(i)]
\item the function $\phi$ is continuous at a point;
\item assuming that $R'$ is locally convex, the function $\phi$ is bounded on a nonvoid open set of $B$;
\item assuming that $R$ is locally compact, $R'$ is locally convex,  the function $\phi$ is bounded on a measurable set of positive measure;
\item assuming that $R$ is locally compact and $R'$ is locally bounded and locally convex, 
the function $\phi$ is measurable on a measurable set of positive measure. 
\end{enumerate}
Then and only then the function $\varphi$  is 
a continuous function or it is an $n$-homomorphism. 
\end{theorem}

\section{Derivations along elementary functions}

Roughly speaking the above presented results dealt with  particular cases of the following problem. 
Assume that $\xi$ 
is a given differentiable function and for the additive function $d\colon \mathbb{R}\to \mathbb{R}$, the mapping 
\[
 x\longmapsto d\left(\xi(x)\right)-\xi'(x)d(x)
\]
is regular on its domain. It is true that in this case 
$d$ admits a representation
\[
 d(x)=\chi(x)+d(1)\cdot x 
\quad 
\left(x\in \mathbb{R}\right), 
\]
where $\chi\colon \mathbb{R}\to \mathbb{R}$ is a real derivation?

In view of the above results, in case $n\in\mathbb{Z}\setminus\left\{0\right\}$ and 
$\left(\begin{array}{cc}
a&b\\
c&d
\end{array}
\right)\in\mathbf{GL}_{2}(\mathbb{Q})$
and the function $\xi$ is 
\[
 \xi(x)=\dfrac{ax^{n}+b}{cx^{n}+d}
\quad 
\left(x\in \mathbb{R}, cx^{n}+d\neq 0\right), 
\]
then the answer is \emph{affirmative}. 
The main aim of this section is to extend this result to other classes of elementary functions such as 
the exponential function, the logarithm function, the trigonometric functions and the hyperbolic functions. 
Concerning such type of investigations, we have to mention the paper of Gy.~Maksa (see \cite{Mak13})\index{Maksa, Gy.}, where the previous 
problem was investigated under the supposition that the mapping 
\[
 x\longmapsto d\left(\xi(x)\right)-\xi'(x)d(x)
\]
is identically zero.

Our main result in this section is contained in the following. 

\begin{theorem}\label{thm2}
 Assume that for the additive function 
$d\colon \mathbb{R}\to \mathbb{R}$ the mapping $\varphi$ defined by 
\[
 \varphi(x)=d\left(\xi(x)\right)-\xi'(x)d(x)
\]
is regular. Then the function $d$ can be represented as 
\[
 d(x)=\chi(x)+d(1)\cdot x 
\quad 
\left(x\in \mathbb{R}\right), 
\]
where $\chi\colon \mathbb{R}\to \mathbb{R}$ is a derivation, 
in any of the following cases
\begin{multicols}{2}
\begin{enumerate}[(a)]
 \item \[\xi(x)=a^{x}\]
\item \[\xi(x)=\cos(x)\]
\item \[\xi(x)=\sin(x)\]
\item \[\xi(x)=\cosh(x)\]
\item \[\xi(x)=\sinh(x). \]
\end{enumerate}
\end{multicols}
\end{theorem}

\begin{biz}
 \begin{enumerate}[{Case} (a)]
  \item Let $a\in \mathbb{R}\setminus\left\{1\right\}$ be an arbitrary positive real number and suppose that the 
mapping $\varphi$ defined by 
\[
 \varphi(x)=d\left(a^{x}\right)-a^{x}\ln(a)d(x) 
\quad 
\left(x\in \mathbb{R}\right)
\]
is regular. 
An easy calculation shows that 
\[
 \varphi(2x)-2a^{x}\varphi(x)=d\left((a^{x})^{2}\right)-2a^{x}d\left(x\right)
\quad 
\left(x\in \mathbb{R}\right), 
\]
that is 
\[
 \varphi\left(2\log_{a}(u)\right)-2u\varphi\left(\log_{a}(u)\right)=
d(u^{2})-2ud(u)
\quad 
\left(u\in ]0, +\infty[\right). 
\]
Due to the regularity of the function $\varphi$, the mapping 
\[
 ]0, +\infty[\ni u\longmapsto \varphi\left(2\log_{a}(u)\right)-2u\varphi\left(\log_{a}(u)\right)
\]
is regular, too. Thus by Theorem \ref{T2.2}, 
\[
 d(x)=\chi(x)+d(1)\cdot x 
\qquad 
\left(x\in \mathbb{R}\right), 
\]
where the function $\chi\colon \mathbb{R}\to \mathbb{R}$ is a derivation. 
\item Assume now that for the additive function $d\colon \mathbb{R}\to \mathbb{R}$, the mapping 
$\varphi$ defined on $\mathbb{R}$ by 
\[
 \varphi(x)=d\left(\cos(x)\right)+\sin(x)d(x) 
\qquad 
\left(x\in \mathbb{R}\right)
\]
is regular. 
If so, then
\[
 \dfrac{\varphi(2x)-4\cos(x)\varphi(x)+d(1)}{2}=
d\left(\cos^{2}(x)\right)-2\cos(x)d\left(x\right)
\]
holds for all $x\in \mathbb{R}$. 
Let now $u\in ]-1, 1[$ and write $\mathrm{arccos}(u)$ in place of $x$ to get 
\[
 \dfrac{\varphi(2\mathrm{arccos}(u))-4u\varphi(\mathrm{arccos}(u))+d(1)}{2}=
d(u^{2})-2ud(u). 
\]
Again, due to the regularity of the function $\varphi$, the mapping 
\[
 ]-1, 1[ \ni u\longmapsto \dfrac{\varphi(2\mathrm{arccos}(u))-4u\varphi(\mathrm{arccos}(u))+d(1)}{2}
\]
is regular, as well. Therefore, Theorem \ref{T2.2} again implies that 
\[
 d(x)=\chi(x)+d(1)\cdot x 
\qquad 
\left(x\in \mathbb{R}\right), 
\]
is fulfilled with a certain real derivation $\chi\colon \mathbb{R}\to \mathbb{R}$. 
\item Suppose that for the additive function $d$, the mapping 
\[
 \varphi(x)=d\left(\sin(x)\right)-\cos(x)d(x)
\qquad 
\left(x\in \mathbb{R}\right)
\]
is regular. 
In this case 
\begin{multline*}
 \varphi\left(x-\frac{\pi}{2}\right)
=
d\left(\sin \left(x-\frac{\pi}{2}\right)\right)-\cos\left(x-\frac{\pi}{2}\right)d\left(x-\frac{\pi}{2}\right)
\\
-d\left(\cos(x)\right)-\sin(x)d(x)+\sin(x)d\left(\frac{\pi}{2}\right), 
\end{multline*}
that is, 
\[
 -\varphi\left(x-\frac{\pi}{2}\right)+\sin(x)d\left(\frac{\pi}{2}\right)=
d\left(\cos(x)\right)+\sin(x)d(x)
\qquad
\left(x\in \mathbb{R}\right). 
\]
 In view of Case (b) this yields that the function $d$ has the desired representation as stated. 
\item Assume the $d\colon \mathbb{R}\to \mathbb{R}$ is an additive function and the mapping 
\[
 \varphi(x)=d\left(\cosh(x)\right)-\sinh(x)d(x)
\qquad 
\left(x\in \mathbb{R}\right)
\]
is regular. The additivity of $d$ and the addition formula of the $\cosh$ function furnish 
\begin{multline*}
 \dfrac{\varphi(2x)-4\cosh(x)\varphi(x)+d(1)}{2}\\=
d\left(\cosh^{2}(x)\right)-2\cosh(x)d\left(\cosh(x)\right) 
\qquad 
\left(x\in \mathbb{R}\right). 
\end{multline*}
Let now $u\in ]1, +\infty[$ arbitrary and put $x=\mathrm{arcosh}(u)$ into the previous identity to get 
\[
 \dfrac{\varphi(2\mathrm{arcosh}(u))-4u\varphi(\mathrm{arcosh}(u))+d(1)}{2}=
d(u^{2})-2ud(u). 
\]
Since the function $\varphi$ is regular, the mapping 
\[
 ]1, +\infty[\ni u\longmapsto \dfrac{\varphi(2\mathrm{arcosh}(u))-4u\varphi(\mathrm{arcosh}(u))+d(1)}{2}
\]
will also be regular. Therefore, Theorem \ref{T2.2} implies again the desired decomposition of the function 
$d$. 
\item Finally, assume the $d\colon \mathbb{R}\to \mathbb{R}$ is an additive function so that 
\[
 \varphi(x)=d\left(\sinh(x)\right)-\cosh(x)d(x)
\qquad 
\left(x\in \mathbb{R}\right)
\]
is regular. 
Let $x, y\in \mathbb{R}$ be arbitrary, then 
\begin{multline*}
\varphi(x+y)= d\left(\sinh(x+y)\right)-\cosh(x+y)d(x+y)
\\=
d\left(\sinh(x)\cosh(y)\right)+d\left(\sinh(y)\cosh(x)\right)
\\-\left[\sinh(x)\sinh(y)+\cosh(x)\cosh(y)\right]d(x+y)
\\=
d\left(\sinh(x)\cosh(y)\right)+d\left(\sinh(y)\cosh(x)\right)
-\sinh(x)\sinh(y)d(x+y)
\\-\cosh(x)d(x)\cosh(y)
-\cosh(x)\cosh(y)d(y)
\end{multline*}
If we use the definition of the function $\varphi$, after some rearrangement, we arrive at 
\begin{multline*}
 \varphi(x+y)-\varphi(x)\cosh(y)-\varphi(y)\cosh(x)
\\
=
d\left(\sinh(x)\cosh(y)\right)+d\left(\sinh(y)\cosh(x)\right)
-\sinh(x)\sinh(y)d(x+y)\\
-\cosh(y)d\left(\sinh(x)\right)-\cosh(x)d\left(\sinh(y)\right)
\end{multline*}
for all $x, y\in \mathbb{R}$. If we replace here $y$ by $-y$, 
\begin{multline*}
 \varphi(x-y)-\varphi(x)\cosh(y)-\varphi(-y)\cosh(x)
\\
=
d\left(\sinh(x)\cosh(y)\right)-d\left(\sinh(y)\cosh(x)\right)
+\sinh(x)\sinh(y)d(x-y)\\
-\cosh(y)d\left(\sinh(x)\right)+\cosh(x)d\left(\sinh(y)\right)
\end{multline*}
can be concluded, where we have also used that the function $\cosh$ is even and the function 
$\sinh$ is odd. 
Adding this two identities,  
\begin{multline*}
 \Phi(x, y)= 2d\left(\sinh(x)\cosh(y)\right)
\\+\sinh(x)\sinh(x)\left[d(x-y)-d(x+y)\right]
-2\cosh(y)d\left(\sinh(x)\right)
\end{multline*}
for any $x, y\in \mathbb{R}$, where 
\begin{multline*}
 \Phi(x, y)= 
\varphi(x+y)-\varphi(x)\cosh(y)-\varphi(y)\cosh(x)
\\+\varphi(x-y)-\varphi(x)\cosh(y)-\varphi(-y)\cosh(x)
\quad 
\left(x, y\in \mathbb{R}\right). 
\end{multline*}
If we put $x=\mathrm{arsinh}(1)$, we get that 
\begin{multline*}
 \dfrac{\Phi\left(\mathrm{arsinh}(1), y\right)+2\cosh(y)d(1)}{2}
\\=
d\left(\cosh(y)\right)-\sinh(y)d(y)
\quad 
\left(y\in \mathbb{R}\right). 
\end{multline*}
Due to the regularity of the function $\varphi$, the mapping 
\[
 \mathbb{R}\ni y\longmapsto \dfrac{\Phi\left(\mathrm{arsinh}(1), y\right)+2\cosh(y)d(1)}{2}
\]
is regular, too. Hence, Case (d) yields the desired form of the function $d$. 
\end{enumerate}
\end{biz}

In what follows, we  extend the list of the functions appearing in the previous statement. 
Therefore we prove the following.

\begin{lem}\label{lem3}
 Let $d\colon \mathbb{R}\to \mathbb{R}$ be an additive function, 
$I\subset \mathbb{R}$ be a nonvoid open interval and 
$\xi\colon I\to \mathbb{R}$ be a continuously differentiable function so that 
the derivative of the function $\xi^{-1}\colon \xi(I)\to \mathbb{R}$ is nowhere zero. 
The mapping 
\[
 I \ni x\longmapsto d(\xi(x))-\xi'(x)d(x)
\]
is regular if and only if the mapping 
\[
 \xi(I)\ni u\longmapsto d(\eta(u))-\eta'(u)d(u)
\]
is regular, where $\eta=\xi^{-1}$. 
\end{lem}
\begin{biz}
 Assume that for the additive function $d$, we have that the mapping 
\[
 \varphi(x)=d(\xi(x))-\xi'(x)d(x) 
\qquad 
\left(x\in I\right)
\]
is regular. 
Let now $u\in \xi(I)$ and put $\xi^{-1}(u)$ in place of $x$ to get 
\[
 -\left(\xi^{-1}\right)'(u)\varphi(\xi^{-1}(u))=
d\left(\xi^{-1}(u)\right)-\left(\xi^{-1}\right)'(u)d(u). 
\]
Due to the regularity of $\varphi$, the mapping appearing in the left hand side is also regular, as stated. 
\end{biz}

In view of Theorem \ref{thm2} and Lemma \ref{lem3}, we immediately obtain the following theorem. 

\begin{cor}\label{cor4}
  Assume that for the additive function 
$d\colon \mathbb{R}\to \mathbb{R}$ the mapping $\varphi$ defined by 
\[
 \varphi(x)=d\left(\xi(x)\right)-\xi'(x)d(x)
\]
is regular. Then the function $d$ can be represented as 
\[
 d(x)=\chi(x)+d(1)\cdot x 
\quad 
\left(x\in \mathbb{R}\right), 
\]
where $\chi\colon \mathbb{R}\to \mathbb{R}$ is a derivation, 
in any of the following cases
\begin{multicols}{2}
\begin{enumerate}[(a)]
 \item \[\xi(x)=\ln(x)\]
\item \[\xi(x)=\mathrm{arccos}(x)\]
\item \[\xi(x)=\mathrm{arcsin}(x)\]
\item \[\xi(x)=\mathrm{arcosh}(x)\]
\item \[\xi(x)=\mathrm{arsinh}(x). \]
\end{enumerate}
\end{multicols}
\end{cor}

Finally, we formulate the following. 

\begin{Opp}
 Let $\alpha\in \mathbb{R}\setminus \left\{0, 1\right\}$ and assume that the additive function 
 $d\colon \mathbb{R}\to \mathbb{R}$ also fulfills 
 \[
  d(x^{\alpha})=\alpha x^{\alpha-1}d(x) 
  \qquad 
  \left(x\in \mathbb{R}, x>0\right)
 \]
Prove or disprove that the function $d$ is a derivation. 
Note that in case $\alpha\in \mathbb{Q}\setminus \left\{0, 1\right\}$ then due to Theorem \ref{T2.2} the answer is 
affirmative, see the paper \cite{BorGse12}, as well. 

Furthermore, we remark that the above identity certainly \emph{does not} characterizes derivations among additive functions, since 
there exists $\alpha\in \mathbb{R}$ and a derivation $d\colon \mathbb{R}\to \mathbb{R}$ for which 
$d$ does not differentiates the function $x\mapsto x^{\alpha}$. 
Indeed, if $\alpha\in \mathbb{R}$ is an irrational algebraic number, 
then due to the Gelfond--Schneider\index{Gelfond, A. O.}\index{Schneider, T.} theorem  (see Gelfond \cite{Gel34})
$2^{\alpha}$ is transcendental. In view of Theorem \ref{T14.2.1}, there exists a derivation 
$d\colon \mathbb{R}\to \mathbb{R}$ so that $d(2^{\alpha})=1$. 
In such a situation however
\[
 1= d(2^{\alpha}) \neq \alpha 2^{\alpha -1}d(2)=0, 
\]
since $d$ has to be identically zero on the set $\mathop{algcl}(\mathbb{Q})$. 
\end{Opp}

Furthermore, we also have to emphasize that in contrast to the previous sections of this chapter, Theorem \ref{thm2} and also 
Corollary \ref{cor4} are not characterization theorems. Thus we also pose the following problem. 

\begin{Opp}
 Prove or disprove that real derivations differentiate the function $\xi$ appearing in Theorem \ref{thm2} and in 
 Corollary \ref{cor4}. 
\end{Opp}

\subsection*{Stability of derivations}

As a starting point of the proof of the main result of this section
the theorem of Hyers will be used. Originally this
statement was formulated in terms of functions that
are acting between Banach spaces, see Hyers \cite{Hye41}\index{Hyers, D.~H.}.
However, we will use this theorem only in the particular case
when the domain and
the range are the set of reals.
In this setting we have the following.

\begin{theorem}
Let $\varepsilon\geq 0$ and suppose that the function
$f\colon \mathbb{R}\rightarrow\mathbb{R}$ fulfills the inequality
\[
\left|f(x+y)-f(x)-f(y)\right|\leq \varepsilon
\]
for all $x,y\in\mathbb{R}$. Then there exists an additive
function $a\colon \mathbb{R}\rightarrow\mathbb{R}$ such that
\[
\left|f(x)-a(x)\right|\leq \varepsilon
\]
holds for arbitrary $x\in\mathbb{R}$.
\end{theorem}

In other words, Hyers' theorem states that if a function $f\colon \mathbb{R}\to \mathbb{R}$ fulfills 
the inequality appearing above, then it can be represented as 
\[
 f(x)=a(x)+b(x) 
\qquad 
\left(x\in \mathbb{R}\right), 
\]
where $a\colon \mathbb{R}\to \mathbb{R}$ is an additive and 
$b\colon \mathbb{R}\to \mathbb{R}$ is a bounded function. Moreover, for all 
$x \in \mathbb{R}$, we also have $\left|b(x)\right|\leq \varepsilon$. 

With the aid of Hyers' theorem and the results of the previous section, we will present several 
stability results in the following. 

Concerning stability properties of derivations the first result is due to R.~Badora, see \cite{Bad06}\index{Badora, R.}, 
where the following result was proved. 

\begin{theorem}
Let $\mathscr{A}_{1}$ be a subalgebra of a Banach algebra $\mathscr{A}$. Assume that the mapping
$f \colon \mathscr{A}_{1} \to \mathscr{A}$ satisfies
\[
\left\|f (x + y) - f (x) - f (y)\right\| \leq \delta 
\qquad 
\left(x, y\in \mathscr{A}_{1}\right)
\]
and
\[
\left\| f (xy) - x f (y) - f (x)  y\right\| \leq \varepsilon
\qquad 
\left(x, y\in \mathscr{A}_{1}\right)
\]
for some constants $\varepsilon, \delta \geq 0$. 
Then there exists a unique derivation, i.e.  a mapping $d \colon \mathscr{A}_{1} \to A$ for which 
\[
 \begin{array}{rcl}
  d(x+y)&=&d(x)+d(x)\\
  d(xy)&=&xd(y)+d(x)y
 \end{array}
\qquad 
\left(x, y\in \mathscr{A}_{1}\right)
\]
so that 
\[
\left\| f (x) - d(x)\right\| \leq \delta  
\qquad 
\left(x\in \mathscr{A}_{1}\right). 
\]
Moreover, we also have
\[
x \cdot (f (y) - d(y)) = 0 
\qquad \left(x, y \in  \mathscr{A}_{1} \right). 
\]
\end{theorem}

Our first result in this area can be found in the paper Boros--Gselmann \cite{BorGse12}\index{Boros, Z.}, where 
we proved the following result.

\begin{theorem}
Let $\varepsilon \geq 0$, $\kappa\in\mathbb{R}$,
$ n, m \in \mathbb{Z} \setminus \{ 0 \} $, $ n \neq m \,$, 
such that 
the function $f:\mathbb{R}\rightarrow\mathbb{R}$
fulfills 
inequality
\begin{equation}\label{Eq3.1}
\left|f(x+y)-f(x)-f(y)\right|\leq \varepsilon_{1}
\end{equation}
for all $x,y\in\mathbb{R}$ and the mapping 
\begin{equation}\label{Eq3.2}
\mathbb{R}\setminus \left\{0\right\}\ni x \longmapsto f(x^{n})-\kappa x^{n-m}f(x^{m})
\end{equation}
is locally bounded.
Then there exist a derivation $F:\mathbb{R}\rightarrow\mathbb{R}$
and $\lambda\in\mathbb{R}$ such that
\[
(n - \kappa m) F(x) = 0
\]
and
\begin{equation}\label{Eq3.3}
\left|f(x)-\left[F(x)+\lambda x\right]\right|\leq \varepsilon 
\end{equation}
are satisfied for all $x\in\mathbb{R}$.
\end{theorem}
\begin{biz}
 Due the theorem of Hyers, inequality \eqref{Eq3.1}
immediately implies that there exists an additive function
$a:\mathbb{R}\rightarrow\mathbb{R}$ satisfying
\begin{equation}\label{Eq3.4}
\left|f(x)-a(x)\right|\leq \varepsilon
\end{equation}
for all $x\in\mathbb{R}$.
In view of inequality \eqref{Eq3.2} this implies that
\begin{multline*}
\left|a(x^{n})-\kappa x^{n-m}a(x^{m})\right| \\
\leq
\left|a(x^{n})-f(x^{n})\right|+
\left|\kappa x^{n-m}\right|\cdot \left|a(x^{m})-f(x^{m})\right|+
\left|f(x^{n})-\kappa x^{n-m}f(x^{m})\right|
\\
\leq 
\varepsilon+\left|\kappa x^{n-m}\right|\varepsilon+\left|f(x^{n})-\kappa x^{n-m}f(x^{m})\right|
=
\left(1+\left|\kappa x^{n-m}\right|\right)\varepsilon_{1}+\left|f(x^{n})-\kappa x^{n-m}f(x^{m})\right| 
\end{multline*}
is fulfilled for all $x\in\mathbb{R}\setminus\left\{0\right\} \,$.
Thus the expression $\left|a(x^{n})-\kappa x^{n-m}a(x^{m})\right|$ is
locally bounded. 
Therefore Lemma \ref{L2.1}. yields that there exists a derivation
$F:\mathbb{R}\rightarrow\mathbb{R}$
such that 
\[ 
(n-\kappa m) F(x) = 0 
\] 
and
\[
a(x)=F(x)+a(1) x
\]
holds for all $x\in\mathbb{R}$.
This, together with \eqref{Eq3.4}, implies \eqref{Eq3.3}
with $ \lambda = a(1) \,$.
\end{biz}

With the aid of the results presented in the previous section we can also 
derive the following stability type statement.

\begin{theorem}
 Let $\varepsilon>0$ and  
$f\colon \mathbb{R}\to \mathbb{R}$ be a function. 
Suppose that 
\begin{enumerate}[(A)]
 \item for all $x, y\in \mathbb{R}$ we have 
\[
 \left|f(x+y)-f(x)-f(y)\right|\leq \varepsilon. 
\]
\item the mapping 
\[
 x\longmapsto f\left(\xi(x)\right)-\xi'(x)f(x)
\]
is locally bounded on its domain, where the function $\xi$ is one of the functions 
\begin{multicols}{2}
 \begin{enumerate}[(a)]
 \item \[a^{x}\]
\item \[\cos(x)\]
\item \[\sin(x)\]
\item \[\cosh(x)\]
\item \[\sinh(x)\]
\item \[\ln(x)\]
\item \[\mathrm{arccos}(x)\]
\item \[\mathrm{arcsin}(x)\]
\item \[\mathrm{arcosh}(x)\]
\item \[\mathrm{arsinh}(x). \]
 \end{enumerate}
\end{multicols}

\end{enumerate}
Then there exist $\lambda\in \mathbb{R}$ and a real derivation $\chi\colon \mathbb{R}\to \mathbb{R}$ 
such that 
\[
 \left|f(x)-\left[\chi(x)+\lambda\cdot x\right]\right|\leq \varepsilon
\]
holds for all $x\in \mathbb{R}$. 
\end{theorem}
\begin{biz}
 Due to assumption (A), we immediately have that 
\[
 f(x)=a(x)+b(x)
\qquad 
\left(x\in \mathbb{R}\right), 
\]
where $a\colon \mathbb{R}\to \mathbb{R}$ is an additive and 
$b\colon \mathbb{R}\to \mathbb{R}$ is a bounded function.
If we use supposition (B), from this we get that the mapping
\[
 x\longmapsto \left[a(\xi(x))-\xi'(x)a(x)\right] +\left[b(\xi(x))-\xi'(x)b(x)\right]
\]
is locally bounded. From this however the local boundedness of the function 
\[
 x\longmapsto a(\xi(x))-\xi'(x)a(x)
\]
can be deduced. In view of the previous statements (see Theorem \ref{thm2} and Corollary \ref{cor4}), 
\[
 a(x)=\chi(x)+a(1)\cdot x 
\qquad 
\left(x\in \mathbb{R}\right)
\]
is fulfilled for any $x\in \mathbb{R}$, where $\chi\colon \mathbb{R}\to \mathbb{R}$ is a real derivation. 
For the function $f$ this means that there exists $\lambda \in \mathbb{R}$ and a real derivation 
$\chi\colon \mathbb{R}\to \mathbb{R}$ so that 
\[
f(x)=\chi(x)+\lambda\cdot x+b(x)
\qquad 
\left(x\in \mathbb{R}\right), 
\]
or equivalently 
\[
 \left|f(x)-\left[\chi(x)+\lambda\cdot x\right]\right|\leq \varepsilon
\]
is satisfied for any $x\in \mathbb{R}$. 
\end{biz}

\clearpage
\addcontentsline{toc}{chapter}{Index}
\printindex

\addcontentsline{toc}{chapter}{Bibliography}

\end{document}